\documentclass[12pt]{article}
\usepackage{eqsection,subeqnarray,indent,amsfonts,amssymb}
\usepackage{pstricks,pst-node,pst-text,pst-3d,pst-plot}
\usepackage{amsmath}     
\usepackage{bm}          
\usepackage{graphicx,epsf}    
\usepackage{cite}
\usepackage{rotating}    

\begin{document}

\bibliographystyle{plain}

\date{December 23, 2008 \\ revised February 18, 2009}

\title{A new Kempe invariant and the (non)-ergodicity of the
       Wang--Swendsen--Koteck\'y algorithm} 

\author{
  {\small Bojan Mohar}\thanks{On leave from
   Department of Mathematics, IMFM \& FMF, University of Ljubljana,
   Ljubljana, Slovenia.}
                                               \\[-1mm]
  {\small\it Department of Mathematics}        \\[-0.2cm]
  {\small\it Simon Fraser University}          \\[-0.2cm]
  {\small\it Burnaby, B.C.~~V5A 1S6, Canada}   \\[-0.2cm]
  {\small\tt mohar@sfu.ca}                     \\[1mm]
  {\small Jes\'us Salas}                       \\[-1mm]
  {\small\it Instituto Gregorio Mill\'an}      \\[-0.2cm]
  {\small\it and}                              \\[-0.2cm]
  {\small\it Grupo de Modelizaci\'on, Simulaci\'on Num\'erica y Matem\'atica 
             Industrial}  \\[-0.2cm]
  {\small\it Universidad Carlos III de Madrid} \\[-0.2cm]
  {\small\it Avda.\  de la Universidad, 30}    \\[-0.2cm]
  {\small\it 28911 Legan\'es, SPAIN}           \\[-0.2cm]
  {\small\tt JSALAS@MATH.UC3M.ES}              \\[-0.2cm] 
  {\protect\makebox[5in]{\quad}}  
  \\
}

\newcommand{\be}{\begin{equation}}
\newcommand{\ee}{\end{equation}} 
\newcommand{\<}{\langle}
\renewcommand{\>}{\rangle}
\newcommand{\widebar}{\overline}
\def\reff#1{(\protect\ref{#1})}
\def\spose#1{\hbox to 0pt{#1\hss}}
\def\ltapprox{\mathrel{\spose{\lower 3pt\hbox{$\mathchar"218$}}
 \raise 2.0pt\hbox{$\mathchar"13C$}}}
\def\gtapprox{\mathrel{\spose{\lower 3pt\hbox{$\mathchar"218$}}
 \raise 2.0pt\hbox{$\mathchar"13E$}}}
\def\textprime{${}^\prime$} 
\def\proof{\par\medskip\noindent{\sc Proof.\ }}
\def\qed{\hbox{\hskip 6pt\vrule width6pt height7pt depth1pt \hskip1pt}\bigskip}
\def\proofof#1{\bigskip\noindent{\sc Proof of #1.\ }}
\def\half{ {1 \over 2} }
\def\third{ {1 \over 3} }
\def\twothird{ {2 \over 3} }
\def\smfrac#1#2{\textstyle{\frac{#1}{#2}}}
\def\smhalf{ \smfrac{1}{2} }
\newcommand{\real}{\mathop{\rm Re}\nolimits}
\renewcommand{\Re}{\mathop{\rm Re}\nolimits}
\newcommand{\imag}{\mathop{\rm Im}\nolimits} 
\renewcommand{\Im}{\mathop{\rm Im}\nolimits}
\newcommand{\sgn}{\mathop{\rm sgn}\nolimits}
\newcommand{\tr}{\mathop{\rm tr}\nolimits}
\newcommand{\diag}{\mathop{\rm diag}\nolimits}
\newcommand{\Gal}{\mathop{\rm Gal}\nolimits}
\newcommand{\mycup}{\mathop{\cup}}
\newcommand{\Arg}{\mathop{\rm Arg}\nolimits}
\def\hboxscript#1{ {\hbox{\scriptsize\em #1}} }
\def\zhat{ {\widehat{Z}} } 
\def\phat{ {\widehat{P}} }
\def\qtilde{ {\widetilde{q}} }
\renewcommand{\mod}{\mathop{\rm mod}\nolimits}
\renewcommand{\emptyset}{\varnothing}

\def\scra{\mathcal{A}}
\def\scrb{\mathcal{B}}
\def\scrc{\mathcal{C}}
\def\scrd{\mathcal{D}}
\def\scrf{\mathcal{F}}
\def\scrg{\mathcal{G}}
\def\scrl{\mathcal{L}}
\def\scro{\mathcal{O}}
\def\scrp{\mathcal{P}}
\def\scrq{\mathcal{Q}}
\def\scrr{\mathcal{R}}
\def\scrs{\mathcal{S}}
\def\scrt{\mathcal{T}}
\def\scrv{\mathcal{V}}
\def\scrz{\mathcal{Z}}

\def\q{{\sf q}}

\def\Z{{\mathbb Z}}
\def\R{{\mathbb R}}
\def\C{{\mathbb C}}
\def\Q{{\mathbb Q}}
\def\N{{\mathbb N}}

\def\T{{\mathsf T}}
\def\H{{\mathsf H}}
\def\V{{\mathsf V}}
\def\D{{\mathsf D}}
\def\J{{\mathsf J}}
\def\P{{\mathsf P}}
\def\QQ{{\mathsf Q}}
\def\RR{{\mathsf R}}

\def\bsigma{{\boldsymbol{\sigma}}}
\def\bone{{\mathbf 1}}
\def\vv{{\bf v}}
\def\uu{{\bf u}}
\def\w{{\bf w}}

\newtheorem{theorem}{Theorem}[section]
\newtheorem{definition}[theorem]{Definition}
\newtheorem{proposition}[theorem]{Proposition}
\newtheorem{lemma}[theorem]{Lemma}
\newtheorem{corollary}[theorem]{Corollary}
\newtheorem{conjecture}[theorem]{Conjecture}


\newenvironment{sarray}{
          \textfont0=\scriptfont0
          \scriptfont0=\scriptscriptfont0
          \textfont1=\scriptfont1
          \scriptfont1=\scriptscriptfont1
          \textfont2=\scriptfont2
          \scriptfont2=\scriptscriptfont2
          \textfont3=\scriptfont3
          \scriptfont3=\scriptscriptfont3
        \renewcommand{\arraystretch}{0.7}
        \begin{array}{l}}{\end{array}}

\newenvironment{scarray}{
          \textfont0=\scriptfont0
          \scriptfont0=\scriptscriptfont0
          \textfont1=\scriptfont1
          \scriptfont1=\scriptscriptfont1
          \textfont2=\scriptfont2
          \scriptfont2=\scriptscriptfont2
          \textfont3=\scriptfont3
          \scriptfont3=\scriptscriptfont3
        \renewcommand{\arraystretch}{0.7}
        \begin{array}{c}}{\end{array}}

\newcommand\Kc{\kappa}  
\newcommand\proofofcase[2]{\bigskip\noindent{\sc Case #1: #2.\ }}
\newcommand\algcr{\mathop{\rm algcr}}  

\maketitle

\thispagestyle{empty}   

\begin{abstract}
We prove that for the class of three-colorable triangulations of a closed
oriented surface, the degree of a four-coloring modulo $12$ is an invariant 
under Kempe changes. We use this general result to prove that for all 
triangulations $T(3L,3M)$ of the torus with $3\leq L\leq M$, there
are at least two Kempe equivalence classes. 
This result implies in particular that the Wang--Swendsen--Koteck\'y algorithm
for the zero-temperature 4--state Potts antiferromagnet on these triangulations
$T(3L,3M)$ of the torus is not ergodic.
\end{abstract}
\clearpage

%
%
\section{Introduction} 

The $q$-state Potts model \cite{Potts_52,Wu_82,Wu_84} is certainly one of 
the simplest and most studied models in Statistical Mechanics.  
However, despite many efforts over more than 50 years, its {\em exact}
solution (even in two dimensions) is still unknown. The ferromagnetic
regime is the best understood case: there are exact (albeit not always
rigorous) results for the location of the critical temperature, 
the order of the transition, etc.   
The antiferromagnetic regime is less understood, partly because 
universality is not expected to hold in general (in contrast with the
ferromagnetic regime); in particular, critical behavior may depend 
on the lattice structure of the model. 
One interesting feature of this antiferromagnetic 
regime is that zero-temperature phase transition may occur for certain values
of $q$ and certain lattices: e.g., the models with $q=2,4$ on the triangular
lattice, and $q=3$ on the square and kagom\'e lattices 
\cite[and references therein]{Salas_Sokal_97}. 

The standard $q$-state Potts model can be defined on any finite undirected 
graph $G = (V,E)$ with vertex set $V$ and edge set $E$. On each vertex
of the graph $i\in V$, we place a spin $\sigma(i)\in \{1,2,\ldots,q\}$,
where $q\ge 2$ in an integer. The spins interact via a Hamiltonian 
\begin{equation}
H(\{\sigma\}) \;=\; -J \sum\limits_{e=ij\in E} \delta_{\sigma(i),\sigma(j)} \,,
\end{equation}
where the sum is over all edges $e \in E$, $J\in\R$ is the coupling constant, 
and $\delta_{a,b}$ is the Kronecker delta.  
The {\em Boltzmann weight}\/ of a configuration is then $e^{-\beta H}$, 
where $\beta \ge 0$ is the inverse temperature.
The {\em partition function}\/ is the sum, taken over all configurations,
of their Boltzmann weights:
\begin{equation}
   Z_G^{\rm Potts}(q, \beta J)  \;=\;  
   \sum_{ \sigma \colon\, V \to \{ 1,2,\ldots,q \} } \; 
      e^{- \beta H(\{\sigma\}) } \,.
 \label{def.ZPotts}
\end{equation}
A coupling $J$ is called {\em ferromagnetic}\/ if $J \ge 0$, 
as it is then favored for adjacent spins to take the same value; and 
{\em antiferromagnetic}\/ if $-\infty \le J \le 0$,
as it is then favored for adjacent spins to take different values.
The zero-temperature ($\beta \to +\infty$) limit of the antiferromagnetic 
($J < 0$) Potts model has an interpretation as a coloring problem: the limit 
$\lim_{\beta.\to+\infty} Z_G^{\rm Potts}(q,-\beta |J|)=P_G(q)$ is the
{\em chromatic polynomial}, which gives the number of proper 
$q$-colorings of $G$. A {\em proper $q$-coloring}\/ of $G$ is a map 
$\sigma \colon\, V \to \{ 1,2,\ldots,q \}$ such that
$\sigma(i) \neq \sigma(j)$ for all pairs of adjacent vertices $ij\in E$.

For many Statistical Mechanics systems for which an exact solution is
not known, Markov Chain Monte Carlo simulations \cite{Bremaud}
have become a very valuable tool to extract physical information.  
An necessary condition for a Markov Chain Monte Carlo algorithm to work is
that it should be ergodic (or irreducible): i.e., the chain can eventually
get from each state to every other state. This condition is usually easy to
check at positive temperature; but in many cases, it becomes a 
highly non-trivial question at zero temperature in the antiferromagnetic regime.

One popular Monte Carlo algorithm for the {\em antiferromagnetic}\/ 
$q$-state Potts model is the Wang--Swendsen--Koteck\'y (WSK) {\em non-local} 
cluster dynamics \cite{WSK_89,WSK_90}. 
At zero temperature (where we expect interesting critical phenomena), it 
leaves invariant the uniform measure over proper $q$-colorings; but 
its ergodicity is a non-trivial question (and not completely 
understood).\footnote{ 
  WSK dynamics can indeed be defined for positive temperature.
  In this case, it is easy to show its ergodicity on the set of {\em all}
  $q$-colorings of the graph $G$ (i.e., proper and non-proper).
}
It is interesting to note that at zero temperature, the basic moves of 
the WSK dynamics correspond to the so-called {\em Kempe changes}, 
introduced by Kempe in his unsuccessful proof of the four-color theorem. 
This connection makes this problem interesting from a purely 
mathematical point of view.  

In this paper we will address the problem of the ergodicity of the 
WSK algorithm for the 4--state Potts antiferromagnet on the triangular 
lattice. Although the Potts model can be defined on any graph $G$, in 
Statistical Mechanics one is mainly interested in ``large'' regular graphs 
embedded on the torus (to minimize finite-size effects). Therefore, we will 
focus on certain regular triangulations of the torus, that we will denoted as 
$T(3L,3M)$ (loosely speaking the triangulation $T(3L,3M)$ is a subset of a 
triangular lattice with linear size $(3L)\times (3M)$ and fully periodic
boundary conditions. For a more detailed definition, see next section). 

The ergodicity of the WSK algorithm for the $q$-state antiferromagnetic 
on the triangular lattice embedded on a torus is only an open question
for $q=4,5,6$. For $q=2$ (the Ising model) it is trivially non-ergodic, as each
WSK move is equivalent to a global spin flip; while for $q=3$ is trivially 
ergodic, as there is a single allowed three-coloring modulo global color 
permutations. On the contrary, for $q\ge 7$ the algorithm is ergodic 
(See next section for more details). Among the unknown cases, $q=4$ is the 
most interesting one, because the system is expected to be critical at zero
temperature.  

Proper 4-colorings of triangulation of the torus are rather special, as 
they can be regarded as maps from a sphere $S^2$ 
(using the tetrahedral representation of the spin) to an orientable surface. 
Therefore, one can borrow concepts from algebraic topology; in particular, 
the degree of a four-coloring. This approach (pioneered by Fisk 
\cite{Fisk_73a,Fisk_77a,Fisk_77b}) can only deal with $q=4$, and cannot 
be extended to the other two cases $q=5,6$. 

Our first goal is to obtain a quantity that is invariant under
a Kempe change (or zero-temperature WSK move), at least for a class 
of triangulations that includes all triangulations of the type $T(3L,3M)$. 
We succeeded in proving 
that for any three-colorable triangulation of a closed orientable surface,
the degree of a four-coloring modulo $12$ is a Kempe invariant. Because
any four-coloring of a closed orientable surface has a degree 
multiple of six, and any three-coloring has degree zero, then we 
conclude that WSK with $q=4$ colors is not ergodic on any three-colorable
triangulation of a closed orientable surface which admits a four-coloring
with degree congruent with 6 modulo 12. 

The next goal is to prove that for any triangulation $T(3L,3M)$ of the 
torus, such four-coloring with degree congruent with 6 modulo 12 exists.
We first proved this statement for any symmetric triangulation 
$T(3L,3L)$ with $L\ge 2$. Then, we extended this result to any
triangulation of the form $T(3L,3M)$ with $L\ge 3$ and $M\ge L$,
and those of the form $T(6,6(2M+1))$ with $M\ge 0$. Therefore,
we conclude that WSK with $q=4$ colors is generically non-ergodic on
the triangulations $T(3L,3M)$ of the torus. 

The paper is organized as follows: In Section~\ref{sec.setup} we
introduce our basic definitions, and review what is 
known in the literature about the problem of the ergodicity of the 
Kempe dynamics. In Section~\ref{sec.4colorings}, we introduce the  
algebraic topology approach borrowed from Fisk. This section includes 
two main results: the proof that the degree modulo 12 is a Kempe 
invariant for a wide enough class of triangulations, and a complete 
proof of Fisk theorem \cite{Fisk_77b} for the class of triangulations
$T(r,s,t)$ of the torus. In the next section, we apply the new invariant
to prove that WSK is non-ergodic on any triangulation $T(3L,3L)$ with
$L\ge 2$. In Section~\ref{sec.asym} we extend the later result to
non-symetric triangulations of the torus $T(3L,3M)$ with $L\ge 3$ and
$M\ge L$ (and also to $T(6,6(2M+1))$ with $M\ge 0$). Finally,
in Section~\ref{sec.summary} we present our conclusions and discuss
prospects of future work.

%
%
\section{Basic setup} \label{sec.setup}

Let $G = (V,E)$ be a finite undirected graph with vertex set $V$ and edge set 
$E$. Then for each graph $G$ there exists 
a polynomial $P_G$ with integer coefficients such that, for each $q \in \Z_+$,
the number of proper $q$-colorings of $G$ is precisely $P_G(q)$.
This polynomial $P_G$ is called the {\em chromatic polynomial}\/ of $G$.
The set of all proper $q$-colorings of $G$ will be denoted as 
$\mathcal{C}_q = \mathcal{C}_q(G)$ (thus, $|\mathcal{C}_q(G)|=P_G(q)$). 

It is far from obvious that $Z_G^{\rm Potts}(q, \beta J)$ 
[cf. \reff{def.ZPotts}], which is defined separately for each positive 
integer $q$, is in fact the restriction to $q \in \Z_+$ of 
a {\em polynomial}\/ in $q$. But this is in fact the case, and indeed we have:

\begin{theorem}[Fortuin--Kasteleyn \protect\cite{Kasteleyn_69,Fortuin_72}
   representation of the Potts model] \label{thm.FK}
\hfill\break
\vspace*{-4mm}
\par\noindent
For every integer $q \ge 1$, we have
\begin{equation}
   Z_G^{\rm Potts}(q, v) \;=\;  
   \sum_{ A \subseteq E }  q^{k(A)} \, v^{|A|} \;,
 \label{eq.FK.identity}
\end{equation}
where $v=e^{\beta J}-1$, and $k(A)$ denotes the number of connected components
in the spanning subgraph $(V,A)$.
\end{theorem}

The foregoing considerations motivate defining the {\em Tutte polynomial}\/ 
of the graph $G$:
\begin{equation}
   Z_G(q, v)   \;=\;
   \sum_{A \subseteq E}  q^{k(A)} \,  v^{|A|} \;,
 \label{def.ZG}
\end{equation}
where $q$ and $v$ are commuting indeterminates. This polynomial is 
equivalent to the standard Tutte polynomial $T_G(x,y)$ after a simple 
change of variables. If we set $v=-1$, we obtain the 
{\em chromatic polynomial} $P_G(q) = Z_G(q,-1)$. In particular, $q$ and
$v$ can be taken as complex variables. See \cite{Sokal_bcc2005} for a
recent survey.

\begin{figure}[htb]
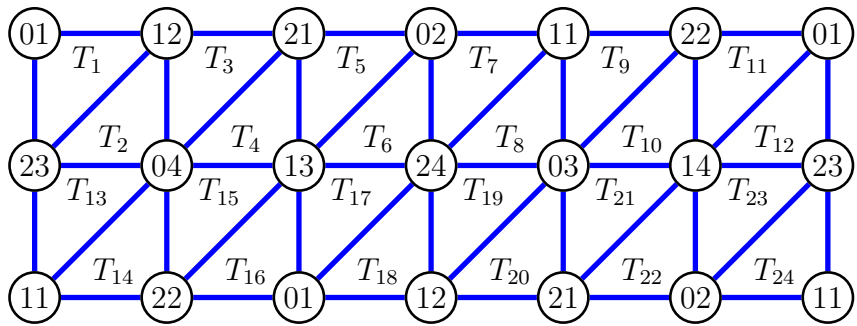

\centering
\psset{xunit=50pt}
\psset{yunit=50pt}
\psset{labelsep=10pt}
\pspicture(-0.5,-0.5)(6.5,2.5)
\multirput{0}(0,0)(0,1){3}{\psline[linewidth=2pt,linecolor=blue](0,0)(6,0)}
\multirput{0}(0,0)(1,0){7}{\psline[linewidth=2pt,linecolor=blue](0,0)(0,2)}
\multirput{0}(0,0)(1,0){5}{\psline[linewidth=2pt,linecolor=blue](0,0)(2,2)}
\psline[linewidth=2pt,linecolor=blue](0,1)(1,2)
\psline[linewidth=2pt,linecolor=blue](5,0)(6,1)
\multirput{0}(0,0)(0,1){3}{%
   \multirput{0}(0,0)(1,0){7}{%
      \pscircle*[linecolor=white]{10pt}
      \pscircle[linewidth=1pt,linecolor=black]{10pt}
   }
}
\rput{0}(0,0){${11}$}
\rput{0}(1,0){${22}$}
\rput{0}(2,0){${01}$}
\rput{0}(3,0){${12}$}
\rput{0}(4,0){${21}$}
\rput{0}(5,0){${02}$}
\rput{0}(6,0){${11}$}

\rput{0}(0,2){${01}$}
\rput{0}(1,2){${12}$}
\rput{0}(2,2){${21}$}
\rput{0}(3,2){${02}$}
\rput{0}(4,2){${11}$}
\rput{0}(5,2){${22}$}
\rput{0}(6,2){${01}$}

\rput{0}(0,1){${23}$}
\rput{0}(1,1){${04}$}
\rput{0}(2,1){${13}$}
\rput{0}(3,1){${24}$}
\rput{0}(4,1){${03}$}
\rput{0}(5,1){${14}$}
\rput{0}(6,1){${23}$}

\uput[90](0.4,1.5){${T_1}$}
\uput[90](1.4,1.5){${T_3}$}
\uput[90](2.4,1.5){${T_5}$}
\uput[90](3.4,1.5){${T_7}$}
\uput[90](4.4,1.5){${T_9}$}
\uput[90](5.4,1.5){${T_{11}}$}

\uput[90](0.4,0.5){${T_{13}}$}
\uput[90](1.4,0.5){${T_{15}}$}
\uput[90](2.4,0.5){${T_{17}}$}
\uput[90](3.4,0.5){${T_{19}}$}
\uput[90](4.4,0.5){${T_{21}}$}
\uput[90](5.4,0.5){${T_{23}}$}

\uput[270](0.6,0.5){${T_{14}}$}
\uput[270](1.6,0.5){${T_{16}}$}
\uput[270](2.6,0.5){${T_{18}}$}
\uput[270](3.6,0.5){${T_{20}}$}
\uput[270](4.6,0.5){${T_{22}}$}
\uput[270](5.6,0.5){${T_{24}}$}

\uput[270](0.6,1.5){${T_{2}}$}
\uput[270](1.6,1.5){${T_{4}}$}
\uput[270](2.6,1.5){${T_{6}}$}
\uput[270](3.6,1.5){${T_{8}}$}
\uput[270](4.6,1.5){${T_{10}}$}
\uput[270](5.6,1.5){${T_{12}}$}

\endpspicture
\caption{\label{figure_T_6_2_2}
The triangulation $T(6,2,2)=\Delta^2 \times \partial\Delta^3$ of the torus. 
Each vertex $x$ of $T(6,2,2)$ is labelled by two integers $ij$,
where  $i$ (resp.\ $j$) corresponds to the associated vertex in 
$\Delta^2$ (resp.\ $\partial\Delta^3$).
The vertices of $\Delta^2$ are labelled $\{0,1,2\}$, while the vertices of
$\partial\Delta^3$ are labelled $\{1,2,3,4\}$. The triangulation
$T(6,2,2)$ has $12$ vertices, and those in the figure with the same label
should be identified. We have also labelled the $24$ triangular faces $T_i$ 
in $T(6,2,2)$.
}
\end{figure}

As explained in the Introduction, we will focus on regular triangulations 
embedded on the torus
The class of regular triangulations of the torus with degree six is
characterized by the following theorem:

\begin{theorem}[Altschulter \protect\cite{Altschulter}] 
Let\/ $T$ be a triangulation of the torus such that all vertices have degree 
six. Then\/ $T$ is one of triangulations $T(r,s,t)$, which are obtained
from the $(r+1)\times (s+1)$ grid by adding diagonals in the squares of
the grid as shown in Figure~\ref{figure_T_6_2_2}, and then identifying
opposite sides to get a triangulation of the torus.
In $T(r,s,t)$ the top and bottom rows have $r$ edges, the left and right
sides $s$ edges. The left and right side are identified as usual; 
but the top and the bottom row are identified after (cyclically) shifting
the top row by $t$ edges to the right. 
\end{theorem}

%
%
%
%

In Figure~\ref{figure_T_6_2_2} we have displayed the triangulation $T(6,2,2)$
of the torus. We will represent these triangulations as 
embedded in a rectangular
grid with three kinds of edges: horizontal, vertical, and diagonal.

The three-colorability of the triangulations $T(r,s,t)$ is given by the
following result (whose proof is left to the reader): 

\begin{proposition}
The triangulation $T(r,s,t)$ is three-colorable if and only if 
$r\equiv 0 \pmod{3}$ and $s-t\equiv 0 \pmod{3}$.
\end{proposition}

In Monte Carlo simulations, it is usual to consider toroidal boundary conditions
with no shifting, so $t=0$. Then, the three-colorability condition reduces
to the standard result $r,s\equiv 0\pmod{3}$. In general, we will consider
the following triangulations of the torus $T(3L,3M,0)=T(3L,3M)$ with 
$L,M\geq 1$. 

The unique three-coloring $c_0$ of $T(3L,3M)$ can be described as: 
\be
c_0(x,y) \;=\; \mod(x+y-2,3) + 1 \,, 
               \quad 1\leq x \leq 3L\,, \quad 1\leq y \leq 3M \,, 
\label{def_coloring_c0}
\ee
where we have explicitly used the above-described embedding of the
triangulation $T(3L,3M)$ in a square grid. 

Finally, in most Monte Carlo simulations one usually considers tori of 
aspect ratio one: i.e., $T(3L,3L)$. This is the class of triangulations we 
are most interested in from the point of view of Statistical Mechanics.  

%
%
\subsection{Kempe changes}

Given a graph $G=(V,E)$ and $q\in\N$, we can define the following dynamics on 
$\mathcal{C}_q$: Choose uniformly at random two distinct colors 
$a,b\in\{1,2,\ldots,q\}$, and let $G_{ab}$ be the induced subgraph of $G$ 
consisting of vertices $x\in V$ for which $\sigma(x)=a$ or $b$. Then, 
independently for each connected component of $G_{ab}$, with probability 
$\smhalf$ either interchange the colors $a$ and $b$ on it, or leave 
the component unchanged. This dynamics is the zero--temperature
limit of the Wang--Swendsen--Koteck\'y (WSK) {\em non-local} cluster 
dynamics \cite{WSK_89,WSK_90} for the antiferromagnetic $q$-state Potts model. 
This zero--temperature Markov chain leaves invariant the uniform
measure over proper $q$-colorings; but its ergodicity cannot be taken
for granted.  

The basic moves of the WSK dynamics correspond to {\em Kempe changes\/} 
(or K-{\em changes}). 
In each K-change, we interchange the colors $a,b$ on
a given connected component (or K-{\em component\/}) of the induced subgraph 
$G_{ab}$.
 
Two $q$-colorings $c_1,c_2\in\mathcal{C}_q(G)$ related by a series of 
K-changes are {\em Kempe equivalent\/} (or K$_q$-{\em equivalent\/}). 
This (equivalence) relation is denoted as $c_1 \stackrel{q}{\sim} c_2$. 
The equivalence classes $\mathcal{C}_q(G)/\stackrel{q}{\sim}$ are called the
{\em Kempe classes\/} (or {\em K$_q$-classes\/}). The number of K$_q$-classes 
of $G$ is denoted by $\Kc(G,q)$. Then, if $\Kc(G,q)>1$, the 
zero-temperature WSK dynamics is not ergodic on $G$ for $q$ colors.

In this paper, we will consider two $q$-colorings related by a {\em global}
color permutation to be the same one. In other words, a $q$-coloring is
actually an equivalence class of standard $q$-colorings modulo global 
color permutations. Thus, the number of (equivalence classes of) proper
$q$-colorings is given by $P_G(q)/q!$. This convention will simplify 
the notation in the sequel. 

%
%
\subsection{The number of Kempe classes}

In this section we will briefly review what it is known in the literature 
about the number of Kempe equivalence classes for several families of graphs. 
The first result implies that WSK dynamics is ergodic on any bipartite 
graph.\footnote{All the cited authors have discovered this theorem 
independently.} 

\begin{proposition}[Burton and Henley \protect\cite{Henley_97a}, %
 Ferreira and Sokal \protect\cite{Sokal_99a}, %
 Mohar \protect\cite{Mohar_05}]
 \label{prop.bipartite}
\hfill\break
\vspace*{-4mm}
\par\noindent
Let\/ $G$ be a bipartite graph and $q\geq 2$ an integer. Then, $\Kc(G,q)=1$.
\end{proposition} 

It is worth noting that Lubin and Sokal \cite{Sokal_93} showed that 
the WSK dynamics with 3 colors is not ergodic on any square--lattice 
grid of size $3M\times 3N$ (with $M,N$ relatively prime) wrapped on a torus.
These graphs are indeed not bipartite. 

The second type of results deals with graphs of bounded maximum degree 
$\Delta$, and shows that $\Kc(G,q)=1$ whenever $q$ is large enough.

\begin{proposition}[Jerrum \protect\cite{Jerrum_private} and %
 Mohar \protect\cite{Mohar_05}]
\label{prop.deltamax}
Let $\Delta$ be the maximum degree of a graph $G$ and let $q\geq \Delta+1$ be
an integer. Then $\Kc(G,q)=1$. If $G$ is connected and contains a vertex 
of degree $<\Delta$, then also $\Kc(G,\Delta)=1$.
\end{proposition}

This result implies that for any 6--regular triangulation $T=T(r,s,t)$,
$\Kc(T,q)=1$ for any $q\geq \Delta+1=7$. However, the cases $q=4,5,6$ 
are not covered by the above proposition. The case $q=3$ is not covered 
either; but this one is trivial if the triangulation is {\em three-colorable}: 
the three-coloring is unique and therefore, $\Kc(T,3)=1$. 

Finally, if we consider planar graphs the situation is better 
understood. Fisk \cite{Fisk_77a} and Moore and Newman \protect\cite{Moore_00} 
showed that $\Kc(T,4)=1$ for planar 3-colorable triangulations. 
Moore and Newman's goal was to establish a height representation
of the corresponding zero-temperature antiferromagnetic Potts model.
One of the authors extended this result as follows:

\begin{theorem}[Mohar \protect\cite{Mohar_05}, Theorem~4.4]
Let $G$ be a three-colorable planar graph. Then $\Kc(G,4)=1$.
\end{theorem}

\begin{corollary}[Mohar \cite{Mohar_05}, Corollary~4.5]
Let $G$ be a planar graph and $q > \chi(G)$. Then $\Kc(G,q)=1$.
\end{corollary}

Indeed, none of our graphs $T(3L,3M)$ is planar. Thus, the above results do 
not apply to our case. 
The main theorem for triangulations appears in \cite{Fisk_77b}.
It involves the notion of the degree of a four-coloring, whose definition
is deferred to the next section.

\begin{theorem}[Fisk \protect\cite{Fisk_77b}] \label{theo_Fisk}
Suppose that\/ $T$ is a triangulation of the sphere, projective plane, or torus. 
If\/ $T$ has a three-coloring, then all four-colorings with degree divisible 
by $12$ are Kempe equivalent.
\end{theorem}

In Section \ref{sect.Fisk_Trst},
we provide a complete self-contained proof of Fisk's result when restricted
to the 6-regular triangulations of the torus treated in this paper.

%
%
\section{Four-colorings of triangulations of the torus} \label{sec.4colorings}

In this section we will consider four-colorings of triangulations of the
torus. Most of the known results concerning this section were obtained
by Fisk \cite{Fisk_73a,Fisk_77a,Fisk_77b}. We will follow his
notation hereafter.

%
%
\subsection{An alternative approach to four-colorings}

Fisk \cite{Fisk_73a,Fisk_77a} considered a definition of a four-coloring
that allows to borrow concepts and results from algebraic topology.
A (proper) four-coloring $f$ of a triangulation $T$ is a non-degenerate 
simplicial map 
\be
f \;\colon\; T \; \longrightarrow \; \partial \Delta^3 \,, 
\ee
where $\partial \Delta^3$ is the surface of a tetrahedron (thus, it can also
be considered as a triangulation of the sphere $S^2$).\footnote{
  A map $f \colon T \to \partial \Delta^3$ is non-degenerate if the image
  of every triangle of $T$ under $f$ is a triangle of $\partial \Delta^3$.
} 
From algebraic topology \cite{Fisk_77a}, if $T$ is the triangulation of an
orientable closed surface (e.g., a sphere or a torus), there is an 
integer-valued function $\deg(f)$ determined up to a sign by $f$. 
In any practical computation, we should choose orientations
for the triangulation $T$ and the tetrahedron $\partial \Delta^3$. Then,
given any triangle $t$ of $\partial \Delta^3$ (i.e., a particular 
three-coloring of a triangular face), we can compute the number $p$ 
(resp.\ $n$) of triangles of $T$ mapping to $t$ which have their 
orientation preserved (resp.\ reversed) by $f$. 
Then, the degree of the four-coloring $f$ is defined as 
\be
\deg(f) \;=\; p-n \,, 
\label{def_deg}
\ee
and it is independent of the choice of the triangle $t$. For instance,
the three-coloring of any triangulation has zero degree, as 
there are no vertices colored $4$, so for $t=124$ we have
$p=n=0$. As we are interested in equivalence classes 
of four-colorings modulo global color permutations, in practical computations 
it only makes sense to consider the absolute value of the degree: i.e.,  
$|\deg(f)|$. 

Tutte \cite{Tutte_69} proved a formula for the degree of a four-coloring 
modulo 2 (the parity of a four-coloring) in terms of the degrees 
of all vertices colored with a specific color. 
We write $\rho(x)$ for the degree of a vertex $x\in V$.
A vertex is {\em even\/} (resp.\ {\em odd\/}) 
if its degree is even (resp.\ odd). 

\begin{lemma}[Tutte \protect\cite{Tutte_69}]
\label{lemma.Tutte}
Given a triangulation $T$ of a closed orientable surface, the degree of a 
four-coloring\/ $f$ of\/ $T$ satisfies
\be
\deg(f) \;\equiv\; \sum\limits_{f(x) = a} \rho(x) \pmod{2} 
\label{eq.lemma.Tutte}
\ee
for $a=1,2,3,4$.
\end{lemma}
 
\proof
By definition, the degree of a four-coloring is modulo 2 equal to the 
number $N$ of triangles of $T$ mapping to a given triangle of 
$\partial \Delta^3$: $\deg(f) \equiv p + n \pmod{2}$ and $N=p+n$. 
If we take a color $a$, which is a vertex of $\partial \Delta^3$, then 
there are three triangular faces of $\partial \Delta^3$ sharing this vertex $a$:
i.e., $t_1$, $t_2$, and $t_3$. For each of these triangles $t_i$, there are 
$N_i$ triangles of $T$ mapping to $t_i$. Then,
\begin{eqnarray}
\deg(f) &\equiv& 3 \deg(f) \pmod{2} \nonumber \\
        &\equiv& N_1 + N_2 + N_3 \pmod{2}  
\end{eqnarray} 
which is equal to the number of triangles of $T$ with a vertex colored $a$.
This number can indeed be written as the r.h.s. of \reff{eq.lemma.Tutte}. \qed 

\bigskip

Lemma \ref{lemma.Tutte} implies that any Eulerian triangulation, in particular,
any triangulation $T(r,s,t)$, can only have four-colorings with even degree, 
as every vertex $x\in V$ has even degree [i.e., $\rho(x)=6$ for any vertex $x$
of $T(r,s,t)$].

A natural question is how many possible values the degree of a 
four-coloring $f$ can take. An answer for a restricted class of triangulations
is given by the following proposition: 

\begin{proposition}[Fisk \protect\cite{Fisk_73a}, %
Problem I.6.6 in \protect\cite{Fisk_77a}] 
\label{prop_Fisk}
Let\/ $T$ be a triangulation of a closed orientable surface, and let\/
$f$ be a four-coloring of\/ $T$. If\/ $T$ admits a three-coloring, 
then $\deg(f) \equiv 0 \pmod{6}$.
\end{proposition}

\proof 
The idea is to mimic the proof of Theorem~4 in \cite{Fisk_77a}.
If $T$ has a three-coloring $h$, and $f$ is a 4-coloring of $T$, then we
can combine these two maps and give
\be
h\times f \;\colon \; T \; \longrightarrow \; 
\Delta^2 \times \partial\Delta^3 \,, 
\ee
where $\Delta^2 \times \partial\Delta^3 =T(6,2,2)$ 
(see Figure~\ref{figure_T_6_2_2}). We have the following diagram
$$
\psset{xunit=1cm}
\psset{yunit=1cm}
\pspicture(-0.5,-0.5)(4.5,2.5)
\rput(0,2){\rnode{G}{$T$}}
\rput(0,0){\rnode{D3}{$\partial \Delta^3$}}
\rput(4,0){\rnode{Prod}{$\Delta^2 \times \partial \Delta^3$}}
\rput(4,2){\rnode{D2}{$\Delta^2$}}
\ncline[nodesep=5pt,linewidth=1pt]{->}{G}{D2}
\ncline[nodesep=5pt,linewidth=1pt]{->}{G}{D3}
\ncline[nodesep=5pt,linewidth=1pt]{->}{Prod}{D2}
\ncline[nodesep=5pt,linewidth=1pt]{->}{Prod}{D3}
\ncline[nodesep=5pt,linewidth=1pt]{->}{G}{Prod}
\uput[180](0,1){$f$}
\uput[0](2,1.1){$h\times f$}
\uput[270](2,0){$g$}
\uput[90](2,2){$h$}
\endpspicture 
$$
where $g$ is the projection of $\Delta^2 \times \partial \Delta^3$ onto its 
second factor $\partial\Delta^3$. By commutativity, 
$\deg(f)=\deg(h\times f)\deg(g)$. 
As the degree of $g$ is $6$, then $\deg(f)= 6\deg(h\times f)\equiv 0 \pmod{6}$.
\qed 

\medskip

In this geometric approach to four-colorings, it is useful to introduce
the concept of a Kempe region \cite{Fisk_77a}.
Suppose that $D$ is a region of the triangulation $T$ (i.e., a union of 
triangles of $T$), and that the four-coloring $f$ uses only two colors on the 
boundary $\partial D$ of $D$. We define a new coloring $g$ of $T$ that is 
equal to $f$ on $T\setminus D$, and equal to $\pi(f)$ on $D$, where 
$\pi$ is the permutation which interchanges the two colors {\em not} on 
$\partial D$. Fisk calls $D$ a Kempe region of $f$, and $\partial D$ a
Kempe cycle. The coloring is {\em not} changed on $\partial D$ itself.
Indeed, inside a Kempe region $D$ we find one or more Kempe components
of the two colors not on $\partial D$. So, the new coloring is K-equivalent
to $f$. Conversely, every K-change can be described as a change on the
region consisting of all triangles containing an edge affected by the
K-change.

Finally it is worth noting that Lemma~\ref{lemma.Tutte} implies that
the parity of a four-coloring [i.e., $\deg(f) \pmod{2}$]
is a Kempe invariant: 

\begin{corollary} \label{cor.Tutte}
Given a triangulation $T$ of a closed orientable surface, then 
the parity of a four-coloring of\/ $T$ is a Kempe invariant.
\end{corollary}

\proof 
If we consider a K-change on a region $D$, we take $a$ to be one of the 
colors on the boundary $\partial D$ (or one of the colors not on the Kempe 
component $T_{bc}$). Then, the parity given by \reff{eq.lemma.Tutte} 
is not affected by the K-change, and therefore, it is an invariant. \qed

\medskip

Unfortunately, the parity is not useful for our purposes, as we are interested
in $6$-regular triangulations of the torus $T(r,s,t)$. Thus, all 
four-colorings have even parity. 
In addition, in the class of three-colorable triangulations of any orientable 
surface, Proposition~\ref{prop_Fisk} ensures that all four-colorings have 
$\deg(f) \equiv 0 \pmod{6}$.

%
%
\subsection{A new Kempe invariant for a class of triangulations}

In this section we shall consider a special class of triangulations in 
which every vertex is of even degree. Such a triangulation is said
to be {\em even} (or Eulerian). 
Observe that every 3-colorable triangulation is even.

Tutte's lemma \ref{lemma.Tutte} implies that if we have a four-coloring $f$
of a triangulation $T$ and we perform a Kempe change to obtain a new 
four-coloring $g$, then 
\be
\deg(g) \;\equiv\; \deg(f) \pmod{2} \,.
\ee
For even triangulations this result has no useful consequences, as all 
four-colorings have even degree. However, for the restricted class of 
three-colorable triangulations of orientable surfaces we can do better. 

\begin{theorem}
\label{main.theo}
Let\/ $T$ be a three-colorable triangulation of a closed orientable surface.
If $f$ and $g$ are two four-colorings of\/ $T$ related by a Kempe change on a 
region $R$, then
\be
\deg(g) \;\equiv \;\deg(f) \pmod{12} \,. 
\label{main.eq}
\ee 
\end{theorem}

\proof
We begin by noting that if $T$ is three-colorable, then it is an even
triangulation. 
Proposition~\ref{prop_Fisk} ensures that $\deg(f),\deg(g)\equiv 0 \pmod{6}$.
As in the proof of Proposition~\ref{prop_Fisk}, we can combine the 
three-color map $h$ with both four-colorings to define the following maps 
\begin{subeqnarray}
 F &=& h\times f \\
 G &=& h\times g 
\end{subeqnarray}
from $T$ onto $\Delta^2 \times \partial \Delta^3 = T(6,6,2)$, where $h$ is
the 3-coloring of $T$. Let us consider the following commutative diagram:
$$
\psset{xunit=1cm}
\psset{yunit=1cm}
\pspicture(-0.5,-0.5)(8.5,2.5)
\rput(0,0){\rnode{D31}{$\partial \Delta^3$}}
\rput(4,0){\rnode{T}{$T$}}
\rput(8,0){\rnode{D32}{$\partial \Delta^3$}}

\rput(0,2){\rnode{T6621}{$\Delta^2 \times \partial \Delta^3$}}
\rput(4,2){\rnode{D2}{$\Delta^2$}}
\rput(8,2){\rnode{T6622}{$\Delta^2 \times \partial \Delta^3$}}

\ncline[nodesep=5pt,linewidth=1pt]{->}{T}{D2}
\ncline[nodesep=5pt,linewidth=1pt]{->}{T}{D31}
\ncline[nodesep=5pt,linewidth=1pt]{->}{T}{D32}
\ncline[nodesep=5pt,linewidth=1pt]{->}{T6621}{D2}
\ncline[nodesep=5pt,linewidth=1pt]{->}{T6622}{D2}
\ncline[nodesep=5pt,linewidth=1pt]{->}{T6621}{D31}
\ncline[nodesep=5pt,linewidth=1pt]{->}{T6622}{D32}
\ncline[nodesep=5pt,linewidth=1pt]{->}{T}{T6621}
\ncline[nodesep=5pt,linewidth=1pt]{->}{T}{T6622}
\uput[180](0,1){$p_2$}
\uput[0]  (8,1){$p_2$}
\uput[90] (2,2){$p_1$}
\uput[90] (6,2){$p_1$}
\uput[270](6,0){$f$}
\uput[270](2,0){$g$}
\uput[0](2,1.1){$G$}
\uput[0](6,0.9){$F$}
\endpspicture 
$$
Since $\deg(f) = \deg(F)\deg(p_2) = 6\deg(F)$ and $\deg(g) = 6\deg(G)$,
our claim is equivalent to $\deg G \equiv \deg F \pmod{2}$. 

For simplicity, let us suppose that there is a Kempe region $R$ such that
its boundary $\partial R$ is colored $3$ and $4$. Then, the Kempe change
on $R$ consists in swapping colors $1$ and $2$ on $R$.
Let us see in detail what happens after this K-change. Consider
Figure~\ref{figure_T_6_2_2} for notation. Triangles in
Figure~\ref{figure_T_6_2_2} are labeled $T_1,\dots,T_{24}$. We say that
a triangle $t$ in $T$ is of {\em type $i$\/} with respect to the coloring
$f$ if it is mapped to $T_i$ by the mapping $F$. Similarly, we consider
types of triangles under $g$.

A triangle of type $T_1$ with 
positive (resp.\ negative) orientation is mapped on a
triangle of type $T_{24}$ with negative (resp.\ positive) orientation after
we swap colors $1$ and $2$. We represent this correspondence
as $\pm T_1 \leftrightarrow \mp T_{24}$.
In fact, this K-change induces a bijection from the set of triangular faces of 
$T(6,2,2)$ onto itself of the form
\begin{subeqnarray}
 \pm T_1     &\leftrightarrow& \mp T_{24}   \\ 
 \pm T_{1+k} &\leftrightarrow& \mp T_{12+k} \,, \qquad 1\le k\le 11. 
\end{subeqnarray}
This correspondence can be written shortly as 
\be
\pm T_k \;\leftrightarrow\; \mp T_{\gamma(k)}
\ee
where $\gamma$ is an appropriate permutation. 
After the K-change, the number of triangles of a given type outside $R$ 
is not changed, so we have to count only the changes inside $R$.
Let us introduce some useful notation: the total number of triangles of a 
given type $k\in\{1,\ldots,24\}$ inside a region $A$ of the triangulation 
$T$ is denoted by $N_k^{(A)}$. Let 
$P_k^{(A)}$ (resp.\ $M_k^{(A)}$) denote the number of triangles of 
type $k$ inside region $A$ with positive (resp.\ negative) orientation. Hence,
\be
N_k^{(A)} \;=\; P_k^{(A)} + M_k^{(A)} \,, \quad k=1,2,\ldots,24 \,, 
                                  \quad  A\subseteq T \,. \nonumber
\ee
If we split the triangulation $T$ into two
regions $R$ and $T\setminus R$, we get  
\be
\deg F \;=\; P^{(T\setminus R)}_k -  M^{(T\setminus R)}_k +
             P^{(R)}_k            -  M^{(R)}_k 
       \,, \quad k=1,2,\ldots,24 \,.  \nonumber
\ee
After the K-change we obtain a new four-coloring $g$. The composite 
coloring $G$ is identical to $F$ outside $R$. The differences can only occur 
inside $R$. The degree of $G$ is given by:
\be
\deg G \;=\; P^{(T\setminus R)}_k -  M^{(T\setminus R)}_k -
             P^{(R)}_{\gamma(k)}       +  M^{(R)}_{\gamma(k)} 
       \,, \quad k=1,2,\ldots,24 \,.  \nonumber
\ee
Let $\Delta\deg = \deg F - \deg G$. Then
\begin{displaymath}
\Delta\deg  \;=\; P^{(R)}_k + P^{(R)}_{\gamma(k)} -  
                     (M^{(R)}_k + M^{(R)}_{\gamma(k)}) \,, 
                \quad k=1,2,\ldots,24 \,.  \nonumber
\end{displaymath}
But this is equivalent to
\begin{eqnarray*}
\Delta\deg      &\equiv&  P^{(R)}_k + P^{(R)}_{\gamma(k)} + 
                          M^{(R)}_k + M^{(R)}_{\gamma(k)} \pmod{2} \nonumber\\
                &\equiv&  N^{(R)}_k + N^{(R)}_{\gamma(k)} \pmod{2} 
   \,, \quad k=1,2,\ldots,24.   \nonumber
\end{eqnarray*} 
In particular, we have that for $k=1,5,9$:
\begin{eqnarray*}
\Delta\deg &\equiv& N^{(R)}_1 + N^{(R)}_{24} \pmod{2}  \nonumber \\ 
\Delta\deg &\equiv& N^{(R)}_5 + N^{(R)}_{16} \pmod{2}  \nonumber \\
\Delta\deg &\equiv& N^{(R)}_9 + N^{(R)}_{20} \pmod{2}  \nonumber 
\end{eqnarray*}
Summing these three equations we arrive at the formula
\begin{eqnarray}
\Delta\deg &\equiv& N^{(R)}_1 + N^{(R)}_{24} + 
                    N^{(R)}_5 + N^{(R)}_{16} +
                    N^{(R)}_9 + N^{(R)}_{20} \pmod{2} \nonumber \\
           &\equiv& \text{\# of triangles inside $R$ with no vertex colored $4$}
               \pmod{2} \nonumber \\
           &\equiv& \text{\# of triangles inside $R$ colored $123$} \pmod{2}
\label{eq.1}
\end{eqnarray}
Note that if we repeat this procedure with $k=3,7,11$ we obtain a 
similar equation and conclude that $\Delta\deg$ has the same parity
as the number of triangles inside $R$ colored $124$.
On the other hand, we cannot obtain a similar formula for the triangles 
colored $134$ or $234$. 

Let us go back to Eq.~\reff{eq.1}. All vertices colored $1$ inside $R$ 
belong to the interior of $R$ (i.e., none of them lies on its boundary,
as $\partial R$ is colored $3,4$). 
In addition, because the triangulation $T$ is even, each interior vertex 
colored $1$ belongs to an even number of triangular faces; all of them 
belonging to $R$. 
Let us consider one of these interior vertices colored $1$, say $x$. 
If none of its neighbors is colored $4$, $x$ contributes $\rho(x)$ 
to $\Delta\deg$ in Eq.~\reff{eq.1}, which is an even number.  
For any neighboring vertex of $x$ colored $4$, this  
contribution is reduced by two.
Thus, for each interior vertex colored $1$, there is an even 
number of triangles belonging to $R$ and colored $123$. This implies that
$\Delta\deg = \deg F - \deg G \equiv 0 \pmod{2}$, and therefore  
\be
\deg f - \deg g \;=\; 6 (\deg F - \deg G) \;\equiv\; 0 \pmod{12} \,, \nonumber
\ee
as claimed. \qed
 
\medskip

Theorem~\ref{main.theo} implies that a four-coloring $f$ with degree 
$\deg f \equiv 6 \pmod{12}$ cannot be K-equivalent to the three-coloring 
$h$, whose degree is zero. This proves the following corollary:

\begin{corollary}
\label{main.corollary}
Let\/ $T$ be a three-colorable triangulation of the torus.
Then $\Kc(T,4)>1$ if and only if there exists a four-coloring $f$ with 
$\deg(f)\equiv 6 \pmod{12}$.
\end{corollary}

\proof
Fisk's Theorem~\ref{theo_Fisk} together with Theorem~\ref{main.theo}
imply the existence of a Kempe equivalence class characterized by 
$\deg(g)\equiv 0 \pmod{12}$. This class includes the three-coloring.
Thus, $\Kc(T,4)>1$ if and only if there is a four-coloring $f$ with
$\deg(f)\equiv 6 \pmod{12}$. \qed 

\medskip

By Theorem \ref{main.theo}, the ``if'' part of Corollary \ref{main.corollary}
holds on arbitrary closed orientable surfaces.

The question of the ergodicity of the WSK dynamics on triangulations
$T(3L,3M)$ reduces to the existence of four-colorings of degree 
$\equiv 6 \pmod{12}$. If there are no such four-colorings, WSK dynamics is
ergodic, while if such four-colorings exist, then WSK dynamics is 
non-ergodic, and the corresponding Markov chain will not converge to the
uniform measure over $\mathcal{C}_4(T)$. 

%
%
\subsection{A complete proof of Fisk's theorem for $\bm{T(r,s,t)}$}

\label{sect.Fisk_Trst}

The proof of Theorem \ref{theo_Fisk} in \cite{Fisk_77b} 
seems to be missing some minor details, as reported in \cite{Mohar_84}. 
However, as far as the authors can see, Fisk's proof is complete 
and correct apart from these minor issues. Nevertheless,
in this section we provide a self-contained proof of Fisk's result 
when restricted to the 6-regular triangulations of the torus treated 
in this paper.
Another advantage of our proof is that it gives a closer insight into Kempe 
equivalence between 4-colorings of triangulations $T(r,s,t)$.

\begin{theorem}
\label{thm_Fisk_for_Trst}
If the triangulation $T(r,s,t)$ admits a\/ $3$-coloring, then every\/
$4$-coloring of $T$, whose degree is divisible by 12, is K-equivalent to the
$3$-coloring.
\end{theorem}

For the proof we shall consider the ``non-singular structure'' of 4-colorings
and show that we can eliminate the ``non-singular'' part completely
by applying K-changes and thus arrive to the 3-coloring. This will be done
by a series of lemmas. But first we need some definitions.

Let $f$ be a 4-coloring of a triangulation $T$. Let $xy\in E(T)$ and let 
$xyz$ and $xyw$ be the two triangles of $T$ containing the edge $xy$. We say that
the edge $xy$ is {\em singular\/} (for the coloring $f$) if $f(z)=f(w)$,
and is {\em non-singular\/} if $f(z)\ne f(w)$. 
Let $N(f)$ be the set of all non-singular edges, and 
for any distinct colors $i,j$, let $N_{ij}=N_{ij}(f)$ be the set of
non-singular edges $xy\in N(f)$ for which $\{f(x),f(y)\}=\{i,j\}$.
For a vertex $x$, let $N_{ij}^x$ be the set of edges in $N_{ij}$ that are
incident with~$x$.

{}From now on we assume that $T=T(r,s,t)$ is a fixed triangulation of the
torus and that $f$ is a 4-coloring of $T$. We also let $i,j\in \{1,\dots,4\}$
be distinct colors used by the 4-coloring~$f$. 

\begin{lemma}
\label{lem:L1}
If\/ $x$ is a vertex of color $f(x)=i$, and $N^x_{ij}\neq \emptyset$,
then $|N^x_{ij}|=2$. Therefore, each $N_{ij}$ is a union of disjoint 
cycles in\/ $T$. If two such cycles, $C\subseteq N_{ij}$ and 
$C'\subseteq N_{il}$ $(j\ne l)$, cross each other at the
vertex $x$, then there is a third cycle $C''\subseteq N_{ik}$ $(k\ne j,l)$
passing through $x$ and crossing both $C$ and $C'$ at $x$.
\end{lemma}

\begin{figure}[htb]
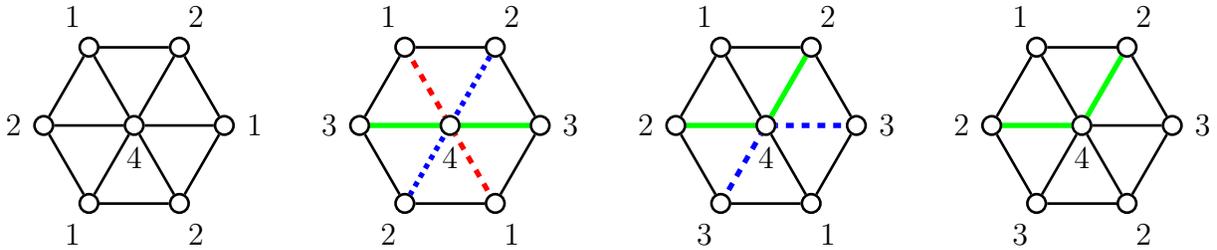

\centering
\psset{xunit=1.2cm}
\psset{yunit=1.2cm}
\psset{labelsep=0.3cm}
\pspicture(-1.5,-1.5)(12,1.5)
\rput{0}(0,0){%
   \psline[linecolor=black,linewidth=1pt](0,0)(1,0)
   \psline[linecolor=black,linewidth=1pt](0,0)(0.500,0.866)
   \psline[linecolor=black,linewidth=1pt](0,0)(-0.500,0.866)
   \psline[linecolor=black,linewidth=1pt](0,0)(-0.500,-0.866)
   \psline[linecolor=black,linewidth=1pt](0,0)(0.500,-0.866)
   \psline[linecolor=black,linewidth=1pt](0,0)(-1,0)
   \psline[linecolor=black,linewidth=1pt](1,0)(0.500,0.866)(-0.500,0.866)%
      (-1,0)(-0.500,-0.866)(0.500,-0.866)(1,0)
   \rput{0}(0.000,0.000){%
      \pscircle*[linecolor=white](0,0){4pt}
      \pscircle[linecolor=black,linewidth=1pt](0,0){4pt}
      \uput[270](0,0){4}
   }
   \rput{0}(1.000,0.000){%
      \pscircle*[linecolor=white](0,0){4pt}
      \pscircle[linecolor=black,linewidth=1pt](0,0){4pt}
      \uput[0.0](0,0){1}
   }
   \rput{0}(0.500,0.866){%
      \pscircle*[linecolor=white](0,0){4pt}
      \pscircle[linecolor=black,linewidth=1pt](0,0){4pt}
      \uput[60.0](0,0){2}
   }
   \rput{0}(-0.500,0.866){%
      \pscircle*[linecolor=white](0,0){4pt}
      \pscircle[linecolor=black,linewidth=1pt](0,0){4pt}
      \uput[120.0](0,0){1}
   }
   \rput{0}(-1.000,0.000){%
      \pscircle*[linecolor=white](0,0){4pt}
      \pscircle[linecolor=black,linewidth=1pt](0,0){4pt}
      \uput[180.0](0,0){2}
   }
   \rput{0}(-0.500,-0.866){%
      \pscircle*[linecolor=white](0,0){4pt}
      \pscircle[linecolor=black,linewidth=1pt](0,0){4pt}
      \uput[240.0](0,0){1}
   }
   \rput{0}(0.500,-0.866){%
      \pscircle*[linecolor=white](0,0){4pt}
      \pscircle[linecolor=black,linewidth=1pt](0,0){4pt}
      \uput[300.0](0,0){2}
   }
}
\rput{0}(3.5,0){%
   \psline[linecolor=green,linewidth=2pt](-1,0)(1,0)
   \psline[linecolor=blue,linewidth=2pt,linestyle=dashed,dash=2pt 2pt]%
           (0.500,0.866)(-0.500,-0.866)
   \psline[linecolor=red,linewidth=2pt,linestyle=dashed,dash=3pt 3pt]%
           (-0.500,0.866)(0.500,-0.866)
   \psline[linecolor=black,linewidth=1pt](1,0)(0.500,0.866)(-0.500,0.866)%
      (-1,0)(-0.500,-0.866)(0.500,-0.866)(1,0)
   \rput{0}(0.000,0.000){%
      \pscircle*[linecolor=white](0,0){4pt}
      \pscircle[linecolor=black,linewidth=1pt](0,0){4pt}
      \uput[270](0,0){4}
   }
   \rput{0}(1.000,0.000){%
      \pscircle*[linecolor=white](0,0){4pt}
      \pscircle[linecolor=black,linewidth=1pt](0,0){4pt}
      \uput[0.0](0,0){3}
   }
   \rput{0}(0.500,0.866){%
      \pscircle*[linecolor=white](0,0){4pt}
      \pscircle[linecolor=black,linewidth=1pt](0,0){4pt}
      \uput[60.0](0,0){2}
   }
   \rput{0}(-0.500,0.866){%
      \pscircle*[linecolor=white](0,0){4pt}
      \pscircle[linecolor=black,linewidth=1pt](0,0){4pt}
      \uput[120.0](0,0){1}
   }
   \rput{0}(-1.000,0.000){%
      \pscircle*[linecolor=white](0,0){4pt}
      \pscircle[linecolor=black,linewidth=1pt](0,0){4pt}
      \uput[180.0](0,0){3}
   }
   \rput{0}(-0.500,-0.866){%
      \pscircle*[linecolor=white](0,0){4pt}
      \pscircle[linecolor=black,linewidth=1pt](0,0){4pt}
      \uput[240.0](0,0){2}
   }
   \rput{0}(0.500,-0.866){%
      \pscircle*[linecolor=white](0,0){4pt}
      \pscircle[linecolor=black,linewidth=1pt](0,0){4pt}
      \uput[300.0](0,0){1}
   }
}
\rput{0}(7,0){%
   \psline[linecolor=green,linewidth=2pt](-1,0)(0,0)(0.500,0.866)
   \psline[linecolor=blue,linewidth=2pt,linestyle=dashed,dash=3pt 3pt]%
           (-0.500,-0.866)(0,0)(1,0)
   \psline[linecolor=black,linewidth=1pt](-0.500,0.866)(0.500,-0.866)
   \psline[linecolor=black,linewidth=1pt](1,0)(0.500,0.866)(-0.500,0.866)%
      (-1,0)(-0.500,-0.866)(0.500,-0.866)(1,0)
   \rput{0}(0.000,0.000){%
      \pscircle*[linecolor=white](0,0){4pt}
      \pscircle[linecolor=black,linewidth=1pt](0,0){4pt}
      \uput[270](0,0){4}
   }
   \rput{0}(1.000,0.000){%
      \pscircle*[linecolor=white](0,0){4pt}
      \pscircle[linecolor=black,linewidth=1pt](0,0){4pt}
      \uput[0.0](0,0){3}
   }
   \rput{0}(0.500,0.866){%
      \pscircle*[linecolor=white](0,0){4pt}
      \pscircle[linecolor=black,linewidth=1pt](0,0){4pt}
      \uput[60.0](0,0){2}
   }
   \rput{0}(-0.500,0.866){%
      \pscircle*[linecolor=white](0,0){4pt}
      \pscircle[linecolor=black,linewidth=1pt](0,0){4pt}
      \uput[120.0](0,0){1}
   }
   \rput{0}(-1.000,0.000){%
      \pscircle*[linecolor=white](0,0){4pt}
      \pscircle[linecolor=black,linewidth=1pt](0,0){4pt}
      \uput[180.0](0,0){2}
   }
   \rput{0}(-0.500,-0.866){%
      \pscircle*[linecolor=white](0,0){4pt}
      \pscircle[linecolor=black,linewidth=1pt](0,0){4pt}
      \uput[240.0](0,0){3}
   }
   \rput{0}(0.500,-0.866){%
      \pscircle*[linecolor=white](0,0){4pt}
      \pscircle[linecolor=black,linewidth=1pt](0,0){4pt}
      \uput[300.0](0,0){1}
   }
}
\rput{0}(10.5,0){%
   \psline[linecolor=green,linewidth=2pt](-1,0)(0,0)(0.500,0.866)
   \psline[linecolor=black,linewidth=1pt](-0.500,0.866)(0.500,-0.866)
   \psline[linecolor=black,linewidth=1pt](-0.500,-0.866)(0,0)(1,0)
   \psline[linecolor=black,linewidth=1pt](1,0)(0.500,0.866)(-0.500,0.866)%
      (-1,0)(-0.500,-0.866)(0.500,-0.866)(1,0)
   \rput{0}(0.000,0.000){%
      \pscircle*[linecolor=white](0,0){4pt}
      \pscircle[linecolor=black,linewidth=1pt](0,0){4pt}
      \uput[270](0,0){4}
   }
   \rput{0}(1.000,0.000){%
      \pscircle*[linecolor=white](0,0){4pt}
      \pscircle[linecolor=black,linewidth=1pt](0,0){4pt}
      \uput[0.0](0,0){3}
   }
   \rput{0}(0.500,0.866){%
      \pscircle*[linecolor=white](0,0){4pt}
      \pscircle[linecolor=black,linewidth=1pt](0,0){4pt}
      \uput[60.0](0,0){2}
   }
   \rput{0}(-0.500,0.866){%
      \pscircle*[linecolor=white](0,0){4pt}
      \pscircle[linecolor=black,linewidth=1pt](0,0){4pt}
      \uput[120.0](0,0){1}
   }
   \rput{0}(-1.000,0.000){%
      \pscircle*[linecolor=white](0,0){4pt}
      \pscircle[linecolor=black,linewidth=1pt](0,0){4pt}
      \uput[180.0](0,0){2}
   }
   \rput{0}(-0.500,-0.866){%
      \pscircle*[linecolor=white](0,0){4pt}
      \pscircle[linecolor=black,linewidth=1pt](0,0){4pt}
      \uput[240.0](0,0){3}
   }
   \rput{0}(0.500,-0.866){%
      \pscircle*[linecolor=white](0,0){4pt}
      \pscircle[linecolor=black,linewidth=1pt](0,0){4pt}
      \uput[300.0](0,0){2}
   }
}
\endpspicture
\caption{Non-singular edges around a vertex.}
\label{fig:f1}
\end{figure}

\proof
Let us consider the 
possible 4-colorings around $x$. Up to symmetries (permutations 
of the colors and the dihedral symmetries of the 6-cycle), there are precisely
four possibilities that are shown in Figure \ref{fig:f1}. The non-singular edges
are drawn by bold solid or broken lines, and a brief inspection shows that the
claims of the lemma hold.
\qed 

%
%
%
%


A 4-coloring $f$ of $T$ is said to be {\em non-singularly minimal\/}
({\em NS-minimal\/} for short) if for any two distinct colors $i,j$,
the non-singular set $N_{ij}$ is either empty or forms a single non-contractible
cycle. The next lemma and its proof explain why such colorings are called
``minimal''.

\begin{lemma}
\label{lem:L2}
Let $f$ be a $4$-coloring of\/ $T$. Then there exists an NS-minimal\/ 
$4$-coloring 
$f'$ of $T$ that is K-equivalent with $f$ and $N(f')\subseteq N(f)$. 
\end{lemma}

\proof
Let $f'$ be a 4-coloring of $T$ that is K-equivalent to $f$, such that
$N(f')\subseteq N(f)$, and $f'$ has minimum number of non-singular edges 
subject to these requirements. Since $f$ has the stated conditions, $f'$ exists.

Let us now consider an arbitrary pair of colors, say 1 and 2.
If $C\subseteq N_{12}(f')$ is a contractible cycle, let $R$ be the disk region
bounded by $C$. By exchanging colors 3 and 4 on $R$ (which keeps us in the 
same K-class), all the change in nonsingular edges is that $C$ becomes singular.
(However, note that particular sets $N_{ij}$ may be changed.)
This contradicts the minimality of $N(f')$. 
Therefore, every non-singular cycle in $N_{12}(f')$ is non-contractible.

Suppose that $N_{12}(f')$ contains distinct cycles $C,C'$. 
As proved above, $C$ and $C'$ are non-contractible. 
By Lemma \ref{lem:L1}, $C$ and $C'$ are disjoint, so 
they are homotopic and therefore together bound a cylinder region $R$.
As above, by exchanging colors 3 and 4 on $R$, we get a contradiction to the
minimality assumption. This completes the proof.
\qed 

%
%
\begin{figure}[htb]
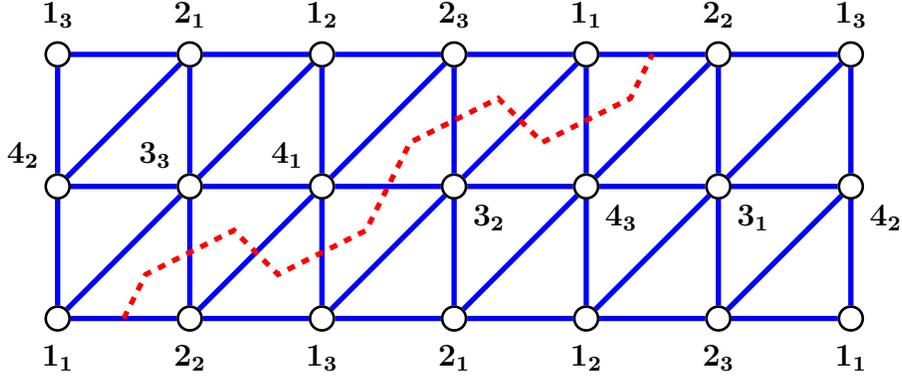

\centering
\psset{xunit=50pt}
\psset{yunit=50pt}
\psset{labelsep=10pt}
\pspicture(-0.5,-0.5)(6.5,2.5)
\multirput{0}(0,0)(0,1){3}{\psline[linewidth=2pt,linecolor=blue](0,0)(6,0)}
\multirput{0}(0,0)(1,0){7}{\psline[linewidth=2pt,linecolor=blue](0,0)(0,2)}
\multirput{0}(0,0)(1,0){5}{\psline[linewidth=2pt,linecolor=blue](0,0)(2,2)}
\psline[linewidth=2pt,linecolor=blue](0,1)(1,2)
\psline[linewidth=2pt,linecolor=blue](5,0)(6,1)
\psline[linewidth=2pt,linecolor=red,linestyle=dashed, dash=3pt 3pt]%
(0.5,0)(0.6667,0.3333)(1.3333,0.6667)(1.6667,0.3333)(2.3333,0.6667)%
(2.6667,1.3333)(3.3333,1.6667)(3.6667,1.3333)(4.3333,1.6667)(4.5,2)
\multirput{0}(0,0)(0,1){3}{%
   \multirput{0}(0,0)(1,0){7}{%
      \pscircle*[linecolor=white]{5pt}
      \pscircle[linewidth=1pt,linecolor=black]{5pt}
   }
}
\uput[270](0,0){$\bm{1_1}$}
\uput[270](1,0){$\bm{2_2}$}
\uput[270](2,0){$\bm{1_3}$}
\uput[270](3,0){$\bm{2_1}$}
\uput[270](4,0){$\bm{1_2}$}
\uput[270](5,0){$\bm{2_3}$}
\uput[270](6,0){$\bm{1_1}$}

\uput[90](0,2){$\bm{1_3}$}
\uput[90](1,2){$\bm{2_1}$}
\uput[90](2,2){$\bm{1_2}$}
\uput[90](3,2){$\bm{2_3}$}
\uput[90](4,2){$\bm{1_1}$}
\uput[90](5,2){$\bm{2_2}$}
\uput[90](6,2){$\bm{1_3}$}

\uput[135](0,1){$\bm{4_2}$}
\uput[135](1,1){$\bm{3_3}$}
\uput[135](2,1){$\bm{4_1}$}
\uput[315](3,1){$\bm{3_2}$}
\uput[315](4,1){$\bm{4_3}$}
\uput[315](5,1){$\bm{3_1}$}
\uput[315](6,1){$\bm{4_2}$}
\endpspicture
\caption{\label{fig:T622}
The triangulation $T_0 = \Delta^2\times\partial\Delta^3 
\approx T(6,2,2)$. The dashed line shows the sequence of triangles
$(g\times f)(\gamma)$ (see text).
}
\end{figure}
%
%

As defined earlier, let $T_0 = \Delta^2\times\partial\Delta^3 \approx T(6,2,2)$
be the 6-regular triangulation of the torus shown in Figure~\ref{fig:T622}.
Note that $T_0$ admits a 3-coloring and a non-singular 4-coloring. Its vertices 
can be labeled by pairs of colors, written as $i_j$, where $i\in\{1,2,3,4\}$
is the color of the non-singular 4-coloring, and $j\in\{1,2,3\}$ is its color
under the 3-coloring; see Figure \ref{fig:T622}. 
If the triangulation $T$ has a 3-coloring $g$ and a 4-coloring $f$, then we
define a simplicial map $g\times f: T\to T_0$ by setting 
$(g\times f)(x) = f(x)_{g(x)}\in V(T_0)$ for every vertex $x$ of $T$.
If $\gamma$ is a closed curve on the torus $T$ that does not pass through
the vertices of $T$, then $\gamma$ can be described (up to homotopy)
 by specifying
the sequence of triangles of $T$ traversed by it. This closed sequence of 
triangles, $A_1,A_2,\dots,A_N,A_1$, is uniquely determined if we cancel out 
possible immediate backtracking, i.e., subsequences of the form $A,B,A$.
The mapping $g\times f$ then determines a closed sequence 
$B_1,B_2,\dots,B_N,B_1$ of triangles in $T_0$, 
where $B_i = (g\times f)(A_i)$ for $i=1,\dots,N$.
This sequence will be denoted by $(g\times f)(\gamma)$ 
(See Figure~\ref{fig:T622}).
The main property of this correspondence is that $B_i=B_{i+1}$ if and only if
the edge common to $A_i$ and $A_{i+1}$ is singular with respect to the
4-coloring $f$ of $T$, i.e. $\gamma$ crosses a singular edge of $f$ when 
passing from $A_i$ to $A_{i+1}$.

\begin{lemma}
\label{lem:L3}
Let\/ $T=T(r,s,t)$ be a $3$-colorable triangulation of the torus, and let
$f$ be an NS-minimal\/ $4$-coloring of\/ $T$. If $f$ is not 
the $3$-coloring of\/ $T$, then all non-singular cycles $N_{ij}$
$(1\le i<j\le 4)$ exist. Two such cycles $N_{ij}$ and $N_{kl}$ 
$(\{i,j\}\ne\{k,l\})$ are homotopic if and only if\/ 
$\{i,j\}\cap\{k,l\}=\emptyset$.
\end{lemma}

\proof
We shall use the notation introduced above. Since $f$ is not the 3-coloring  
(which is unique, up to global permutations of colors), 
we may assume that $N_{12}\ne\emptyset$. Let $\gamma$ be a simple closed curve 
in the torus that crosses $N_{ij}$ precisely once and is given by the sequence
of triangles $A_1,\dots,A_N,A_1$. Let us consider the corresponding sequence
$\gamma' = (g\times f)(\gamma)=B_1,B_2,\dots,B_N,B_1$ of triangles in $T_0$.

Let $K_{ij}$ be the non-singular cycle in $T_0$ passing through all vertices
$i_l$ and $j_l$, $l=1,2,3$. Since $\gamma$ crosses $N_{12}$ precisely once,
$\gamma'$ crosses $K_{12}$ exactly once. We may assume that
it crosses $K_{12}$ through the edge $e=1_12_2$ as shown in 
Figure~\ref{fig:T622}.

For a cycle $K_{ij}$ we define the {\em algebraic crossing number\/} with 
$\gamma'$ by first counting the number of consecutive triangles 
$B_l,B_{l+1}$ in 
$\gamma'$ such that $B_l$ is ``on the left'' of $K_{ij}$, while $B_{l+1}$ is
``on the right'' of it, and then subtracting the number of such pairs,
where $B_l$ is ``on the right'' and $B_{l+1}$ is ``on the left''. (For the
two ``horizontal'' cycles $K_{12}$ and $K_{34}$ we replace ``left'' 
by ``bottom'' and ``right'' by ``top''. All of these directions of 
course refer to Figure~\ref{fig:T622}.) 
We denote this number by $\algcr(\gamma',K_{ij})$.

For an arbitrary edge-set $F\subseteq E(K_{ij})$, we define $\algcr(\gamma',F)$
in the same way, except that we only consider consecutive
triangles $B_l,B_{l+1}$ sharing the edges in $F$. 
Let $k = \algcr(\gamma',\{1_14_2,4_21_3\})$.
This number can be viewed as the ``winding number'' around the cylinder 
obtained from $T_0$ by cutting along the cycle $K_{12}$, 
cf.~Figure~\ref{fig:T622}. 
Using the fact that $\gamma'$ is contained in this cylinder except for its
crossing of the edge $1_12_2$, it is easy to see that 
$\algcr(\gamma',K_{13}) = 3k+1$, $\algcr(\gamma',K_{24}) = 3k+1$,
$\algcr(\gamma',K_{14}) = 3k+2$, and $\algcr(\gamma',K_{23}) = 3k+2$.
Moreover, $\algcr(\gamma',K_{12}) = \algcr(\gamma',K_{34}) = 1$.
In particular, none of these numbers is zero (modulo 3). 

Let us recall that $B_i\ne B_{i+1}$ if and only if the edge common to $A_i$ and 
$A_{i+1}$ is non-singular with respect to $f$. Therefore, $\gamma'$ crosses
an edge of $K_{ij}$ precisely when $\gamma$ crosses an edge in $N_{ij}(f)$.
Therefore $\algcr(\gamma',K_{ij}) = \algcr(\gamma,N_{ij}) \ne 0$.
This shows that none of the sets $N_{ij}$ is empty. 

If $\{i,j\}\cap\{k,l\}=\emptyset$, the two cycles $N_{ij}$ and $N_{kl}$ are
disjoint. Since they are non-contractible and the surface is the torus, they
are homotopic to each other. On the other hand, since
$\algcr(\gamma,N_{13}) = \algcr(\gamma,N_{14}) - 1$, cycles $N_{13}$ and
$N_{14}$ cannot be homotopic. Similarly, by starting the above proof with other 
cycles instead of $N_{12}$, we conclude that cycles $N_{ij}$ and $N_{kl}$ 
cannot be homotopic if $\{i,j\}\cap\{k,l\}\ne\emptyset$.
\qed 

\medskip

Note that in the proof of Lemma \ref{lem:L3}, we did not use any assumption
on the degree of the 4-coloring $f$. On the other hand, in our last lemma,
when arguing about the degree of a 4-coloring, we will not need the existence 
of the 3-coloring. 

\begin{lemma}
\label{lem:L4}
Let $f$ be an NS-minimal\/ $4$-coloring of\/ $T$ such that all non-singular
cycles $N_{ij}(f)$ exist and such that two such cycles $N_{ij}$ and $N_{kl}$ 
$(\{i,j\}\ne\{k,l\})$ are homotopic if and only 
if\/ $\{i,j\}\cap\{k,l\}=\emptyset$.
Then the degree of $f$ is congruent to $2$ modulo $4$.
In particular, it is not divisible by\/ $12$.
\end{lemma}

\proof
Let us consider cycles $N_{12}$ and $N_{13}$. Since they are not homotopic, 
they cross at least once, and this happens at vertices of color 1. 
By Lemma \ref{lem:L1}, both these cycles are crossed by $N_{14}$ at each
such crossing point. Let us fix an orientation on the torus $T$ and let
$x\in V(T)$ be a vertex of color 1 at which $N_{12},N_{13},N_{14}$ cross
each other. If the local clockwise order around $x$ is
$N_{12},N_{13},N_{14},N_{12},N_{13},N_{14}$, then we say that $x$ is 
a {\em positive crossing point\/} (of color 1); if the local clockwise order 
is $N_{12},N_{14},N_{13},N_{12},N_{14},N_{13}$, then $x$ is 
a {\em negative crossing point}. 

We claim that the difference of the number of positive minus the number 
of negative crossing points of color 1 is equal (in absolute value) 
to the algebraic crossing number $\algcr(N_{12},N_{13})$. 
This is a consequence of the fact that color 4 changes 
sides from left to right side of $N_{13}$, or vice versa,
every time when the curve $N_{13}$ passes through a crossing point of 
color 1 or through a crossing point of color 3 (thus crossing
the cycle $N_{34}$ which is homotopic to $N_{12}$). 
We leave the details to the reader.

Since the numbers of positive and negative crossing points of color 
1 are also the same for other pairs of non-singular cycles that involve 
color 1, we conclude that
\begin{equation}
|\algcr(N_{12},N_{13})| \;=\; |\algcr(N_{12},N_{14})| \;=\; 
|\algcr(N_{13},N_{14})|\,.
\label{eq:algcr1}
\end{equation}

Let us fix two simple closed curves $\gamma,\nu$ on the torus $T$, 
where $\nu$ is the curve corresponding to the cycle $N_{12}(f)$ and $\gamma$ 
crosses $\nu$ precisely once. Then every closed curve $\alpha$ on $T$ 
is homotopic to the curve which winds $a$ times around $\nu$, and then 
winds $b$ times around $\gamma$, where $a$ and $b$ are
integers. We say that $\alpha$ has {\em homotopy type\/} $(a,b)$.
The homotopy type of $N_{12}$ is clearly (1,0). Let $(a,b)$ and $(c,d)$ be the
homotopy types of $N_{13}$ and $N_{14}$, respectively. The algebraic 
crossing number between closed curves is a (free) homotopy invariant
and can be expressed as the determinant of the $2\times 2$ matrix 
whose rows are the homotopy types of the curves (see, e.g.~\cite{Zieschang}). 
In particular,
\begin{eqnarray}
 \algcr(N_{12},N_{13}) &=& \pm\det\begin{pmatrix}1&0\\a&b\end{pmatrix} 
 \;=\; \pm b \,, \label{eq:algcr2}\\
 \algcr(N_{12},N_{14}) &=& \pm\det\begin{pmatrix}1&0\\c&d\end{pmatrix} 
 \;=\; \pm d\,, \label{eq:algcr3}\\
 \algcr(N_{13},N_{14}) &=& \pm\det\begin{pmatrix}a&b\\c&d\end{pmatrix} 
 \;=\; \pm(ad-bc) \,. \label{eq:algcr4}
\end{eqnarray}
By (\ref{eq:algcr1}), all three algebraic crossing numbers in 
(\ref{eq:algcr2})--(\ref{eq:algcr4}) are equal up to the sign, so
$|b| = |d| = |ad-bc|$. It follows that either $|a-c|=1$ or $|a+c|=1$. 
Here we have used the fact that $b\ne0$, and this is true since $N_{13}$ 
is not homotopic to $N_{12}$. A particular consequence of the above conclusion
is that either $a$ or $c$ is even.

Suppose first that $a$ is even. Since $N_{13}$ is a simple curve, 
its homotopy type $(a,b)$ satisfies $\gcd(a,b)=1$ (cf.~\cite{Zieschang}).
Therefore $b$ and henceforth also $d$ are odd. 

The other case is when $c$ is even. In that case, we derive the same 
conclusion as above. From this it follows that the
total number of crossing points of color 1 is odd. 
Of course, we can repeat the same proof for crossing points of color 2
to conclude that their number is odd as well.

We are ready for the second part of the proof, where we will relate the number
of crossing points and the degree of the coloring $f$. 
Let us traverse the cycle $N_{12}$ and consider the (cyclic) sequence of 
all crossing points of colors 1 and 2 as they appear on $N_{12}$. 
We shall see that one can determine the degree of $f$ just from
this sequence.

Let us recall that $\deg(f)$ is equal to the difference between the number
of triangles colored $123$, whose orientation on the surface is $123$, 
minus the number of such triangles whose orientation is $132$.
If $t$ is such a triangle and its edge colored $12$ is not in $N_{12}$,
then there is another triangle colored 123 sharing that edge with $t$ and
having opposite orientation. The contribution of all such triangles towards
the degree of $f$ thus cancels out. 
On the other hand, each edge of $N_{12}$ is contained in precisely one triangle
colored $123$. Consider two consecutive edges $xy$ and $yz$ on $N_{12}$. If $y$
is not a crossing point with other non-singular curves, then one of the two
triangles colored $123$ and incident with these edges is oriented positively,
the other one negatively, and so their contributions will cancel out.
On the other hand, if $y$ is a crossing point, then they have the 
same orientation. If two consecutive crossing points on $N_{12}$ are of 
the same color, then the pair at one of these two crossing points is 
positively oriented, while the pair at the other crossing point is 
negatively oriented, and hence they cancel out. This has the same 
effect as removing two consecutive 1's or two consecutive 2's from the 
cyclic sequence of crossing points on $N_{12}$. Therefore, we may 
assume that the sequence of crossing points is alternating, 
$1212\dots 12$. The number of 1's
is an odd integer, say $2k+1$, as shown in the first part of the proof.
This implies that all triangles at crossing points have positive (or all have
negative) orientation. Therefore, $\deg(f) = \pm 2(2k+1) \equiv 2 \pmod{4}$, 
which we were to prove.
\qed 

\medskip

\proofof{Theorem~\ref{thm_Fisk_for_Trst}} 
Let $f$ be a 4-coloring of $T=T(r,s,t)$. By Lemma \ref{lem:L2}
there is an NS-minimal coloring $f'$ that is K-equivalent to $f$ and 
has $N(f')\subseteq N(f)$. If $f'$ is not the 3-coloring, then
by Lemma \ref{lem:L3}, all six non-singular curves $N_{ij}(f')$ exist 
and their homotopy is as stated in the lemma. But then Lemma \ref{lem:L4}
implies that $\deg(f')\equiv 2 \pmod{4}$. Since the K-equivalence
preserves the value of the degree modulo 12 (cf.~Theorem~\ref{main.theo}), 
this yields a contradiction 
to the assumption that the degree of $f$ is divisible by 12.
\qed 

%
%
\section{Consequences for the triangulations $\bm{T(3L,3L)}$} \label{sec.sym}

A simple corollary of Proposition~\ref{prop_Fisk} and Theorem~\ref{theo_Fisk}
shows that all 4-colorings of $T(3,3)$ are K-equivalent: 

\begin{corollary} \label{theo_L=3}
$\Kc(T(3,3),4)=1$.
\end{corollary}

\proof
The smallest (in modulus) non-zero degree for a four-coloring of an even  
three-colorable triangulation is $6$ by Proposition~\ref{prop_Fisk}. 
But in order to have a four-coloring $f$ with such degree, we would need 
at least $6\times 4=24$ triangular faces. 
However, the triangulation $T(3,3)$ only has $3^2 \times 2=18$ such faces. 
Then, $\deg(f)=0$ for all four-colorings of $T(3,3)$, and 
Theorem~\ref{theo_Fisk} implies that $\Kc(T(3,3),4)=1$. \qed 

\medskip

A four-coloring $f$ is said to be {\em non-singular\/} if all edges 
are non-singular with respect to $f$. Fisk \cite{Fisk_77b} showed that 
the triangulation
$T(r,s,t)$ has a non-singular four-coloring $c_\text{ns}$ if and only if
$r,s,t$ are all even. In this non-singular coloring, each horizontal 
row uses exactly two colors. This also holds for all vertical and diagonal
``straight-ahead cycles''. For the triangulation $T(3L,3M)$, the
non-singular coloring is given by 
\be
c_\text{ns}(x,y) \;=\; \begin{cases} 
                1 &  \text{if $x,y\equiv 1 \bmod{2}$} \\ 
                2 &  \text{if $x  \equiv 1$ and 
                              $y  \equiv 0 \bmod{2}$} \\ 
                3 &  \text{if $x  \equiv 0$  and 
                              $y  \equiv 1 \bmod{2}$ }\\ 
                4 &  \text{if $x,y\equiv 0 \bmod{2}$} 
\end{cases} \,, \quad 1\leq x \leq 3L\,, \quad 1\leq y \leq 3M 
\label{def_coloring_ns}
\ee

\begin{proposition}
\label{prop.non-singular}
The triangulation $T(3L,3M)$ has a non-singular four-coloring $c_\text{ns}$
if and only $L=2\ell$ and $M=2m$ are both even. If so, then
$|\deg c_\text{ns}|=18\ell m$. In particular, $\Kc(T(6\ell,6m),4)\ge2$
if $\ell$ and $m$ are both odd.
\end{proposition}

\proof
Under the non-singular coloring, all triangles are mapped to
$\partial\Delta^3$ with the same orientation. Thus,
$|\deg c_\text{ns}| = \frac{1}{4}(\#\text{triangles of }T(3L,3M)) = 18\ell m$.
If $\ell$ and $m$ are both odd, the degree is $\equiv 6 \bmod{12}$,
and now Corollary~\ref{main.corollary} applies.
\qed

The next non-trivial result shows that $\Kc(T(6,6),4)=2$; hence
WSK dynamics is non-ergodic on this triangulation. 

\begin{theorem}[with Alan Sokal]\label{theo_L=6}
$\Kc(T(6,6),4)=2$.
\end{theorem}

\proof
Proposition \ref{prop.non-singular} shows that the non-singular
four-coloring of $T(6,6)$ has $\deg(c_\text{ns})\equiv 6 \pmod{12}$
and that there are at least two Kempe equivalence classes for this 
triangulation. One class $\mathcal{C}_4^{(0)}$
corresponds to all colorings whose degree is a multiple of $12$. The other
classes contain colorings with degree $\equiv 6 \pmod{12}$. 

The fact that the number of Kempe classes is exactly two can be derived 
as follows. Let us first observe that the maximum degree of a four-coloring
of the triangulation $T(3L,3L)$ is $\lfloor 9L^2/2 \rfloor$; therefore,
for $T(6,6)$ this maximum degree is $18$. Thus, we should focus on all
four-colorings $f$ with $|\deg(f)|=6,18$, and show that they form a unique
Kempe equivalence class. 

\begin{figure}[htb]
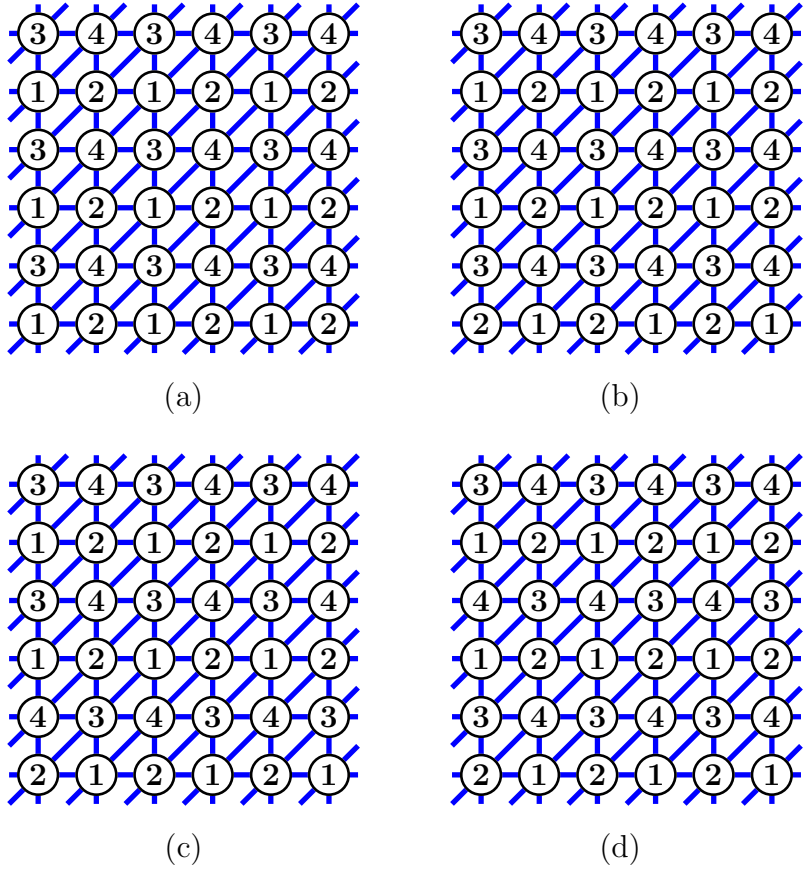

\centering
  \begin{tabular}{cc}
%
%
\psset{xunit=22pt}
\psset{yunit=22pt}
\pspicture(-0.5,-0.5)(5.5,5.5)
\psline[linewidth=2pt,linecolor=blue](-0.5,0)(5.5,0)
\psline[linewidth=2pt,linecolor=blue](-0.5,1)(5.5,1)
\psline[linewidth=2pt,linecolor=blue](-0.5,2)(5.5,2)
\psline[linewidth=2pt,linecolor=blue](-0.5,3)(5.5,3)
\psline[linewidth=2pt,linecolor=blue](-0.5,4)(5.5,4)
\psline[linewidth=2pt,linecolor=blue](-0.5,5)(5.5,5)
\psline[linewidth=2pt,linecolor=blue](0,-0.5)(0,5.5)
\psline[linewidth=2pt,linecolor=blue](1,-0.5)(1,5.5)
\psline[linewidth=2pt,linecolor=blue](2,-0.5)(2,5.5)
\psline[linewidth=2pt,linecolor=blue](3,-0.5)(3,5.5)
\psline[linewidth=2pt,linecolor=blue](4,-0.5)(4,5.5)
\psline[linewidth=2pt,linecolor=blue](5,-0.5)(5,5.5)
\psline[linewidth=2pt,linecolor=blue](-0.5,-0.5)(5.5,5.5)
\psline[linewidth=2pt,linecolor=blue](0.5,-0.5)(5.5,4.5)
\psline[linewidth=2pt,linecolor=blue](1.5,-0.5)(5.5,3.5)
\psline[linewidth=2pt,linecolor=blue](2.5,-0.5)(5.5,2.5)
\psline[linewidth=2pt,linecolor=blue](3.5,-0.5)(5.5,1.5)
\psline[linewidth=2pt,linecolor=blue](4.5,-0.5)(5.5,0.5)
\psline[linewidth=2pt,linecolor=blue](-0.5,0.5)(4.5,5.5)
\psline[linewidth=2pt,linecolor=blue](-0.5,1.5)(3.5,5.5)
\psline[linewidth=2pt,linecolor=blue](-0.5,2.5)(2.5,5.5)
\psline[linewidth=2pt,linecolor=blue](-0.5,3.5)(1.5,5.5)
\psline[linewidth=2pt,linecolor=blue](-0.5,4.5)(0.5,5.5)
\multirput{0}(0,0)(0,1){6}{%
  \multirput{0}(0,0)(1,0){6}{%
     \pscircle*[linecolor=white]{8pt}
     \pscircle[linewidth=1pt,linecolor=black] {8pt}
   }
}
\rput{0}(0,0){{\bf 1}}
\rput{0}(0,1){{\bf 3}}
\rput{0}(0,2){{\bf 1}}
\rput{0}(0,3){{\bf 3}}
\rput{0}(0,4){{\bf 1}}
\rput{0}(0,5){{\bf 3}}
\rput{0}(1,0){{\bf 2}}
\rput{0}(1,1){{\bf 4}}
\rput{0}(1,2){{\bf 2}}
\rput{0}(1,3){{\bf 4}}
\rput{0}(1,4){{\bf 2}}
\rput{0}(1,5){{\bf 4}}
\rput{0}(2,0){{\bf 1}}
\rput{0}(2,1){{\bf 3}}
\rput{0}(2,2){{\bf 1}}
\rput{0}(2,3){{\bf 3}}
\rput{0}(2,4){{\bf 1}}
\rput{0}(2,5){{\bf 3}}
\rput{0}(3,0){{\bf 2}}
\rput{0}(3,1){{\bf 4}}
\rput{0}(3,2){{\bf 2}}
\rput{0}(3,3){{\bf 4}}
\rput{0}(3,4){{\bf 2}}
\rput{0}(3,5){{\bf 4}}
\rput{0}(4,0){{\bf 1}}
\rput{0}(4,1){{\bf 3}}
\rput{0}(4,2){{\bf 1}}
\rput{0}(4,3){{\bf 3}}
\rput{0}(4,4){{\bf 1}}
\rput{0}(4,5){{\bf 3}}
\rput{0}(5,0){{\bf 2}}
\rput{0}(5,1){{\bf 4}}
\rput{0}(5,2){{\bf 2}}
\rput{0}(5,3){{\bf 4}}
\rput{0}(5,4){{\bf 2}}
\rput{0}(5,5){{\bf 4}}
\endpspicture
\qquad 
& 
\qquad 
%
%
\psset{xunit=22pt}
\psset{yunit=22pt}
\pspicture(-0.5,-0.5)(5.5,5.5)
\psline[linewidth=2pt,linecolor=blue](-0.5,0)(5.5,0)
\psline[linewidth=2pt,linecolor=blue](-0.5,1)(5.5,1)
\psline[linewidth=2pt,linecolor=blue](-0.5,2)(5.5,2)
\psline[linewidth=2pt,linecolor=blue](-0.5,3)(5.5,3)
\psline[linewidth=2pt,linecolor=blue](-0.5,4)(5.5,4)
\psline[linewidth=2pt,linecolor=blue](-0.5,5)(5.5,5)
\psline[linewidth=2pt,linecolor=blue](0,-0.5)(0,5.5)
\psline[linewidth=2pt,linecolor=blue](1,-0.5)(1,5.5)
\psline[linewidth=2pt,linecolor=blue](2,-0.5)(2,5.5)
\psline[linewidth=2pt,linecolor=blue](3,-0.5)(3,5.5)
\psline[linewidth=2pt,linecolor=blue](4,-0.5)(4,5.5)
\psline[linewidth=2pt,linecolor=blue](5,-0.5)(5,5.5)
\psline[linewidth=2pt,linecolor=blue](-0.5,-0.5)(5.5,5.5)
\psline[linewidth=2pt,linecolor=blue](0.5,-0.5)(5.5,4.5)
\psline[linewidth=2pt,linecolor=blue](1.5,-0.5)(5.5,3.5)
\psline[linewidth=2pt,linecolor=blue](2.5,-0.5)(5.5,2.5)
\psline[linewidth=2pt,linecolor=blue](3.5,-0.5)(5.5,1.5)
\psline[linewidth=2pt,linecolor=blue](4.5,-0.5)(5.5,0.5)
\psline[linewidth=2pt,linecolor=blue](-0.5,0.5)(4.5,5.5)
\psline[linewidth=2pt,linecolor=blue](-0.5,1.5)(3.5,5.5)
\psline[linewidth=2pt,linecolor=blue](-0.5,2.5)(2.5,5.5)
\psline[linewidth=2pt,linecolor=blue](-0.5,3.5)(1.5,5.5)
\psline[linewidth=2pt,linecolor=blue](-0.5,4.5)(0.5,5.5)
\multirput{0}(0,0)(0,1){6}{%
  \multirput{0}(0,0)(1,0){6}{%
     \pscircle*[linecolor=white]{8pt}
     \pscircle[linewidth=1pt,linecolor=black] {8pt}
   }
}
\rput{0}(0,0){{\bf 2}}
\rput{0}(0,1){{\bf 3}}
\rput{0}(0,2){{\bf 1}}
\rput{0}(0,3){{\bf 3}}
\rput{0}(0,4){{\bf 1}}
\rput{0}(0,5){{\bf 3}}
\rput{0}(1,0){{\bf 1}}
\rput{0}(1,1){{\bf 4}}
\rput{0}(1,2){{\bf 2}}
\rput{0}(1,3){{\bf 4}}
\rput{0}(1,4){{\bf 2}}
\rput{0}(1,5){{\bf 4}}
\rput{0}(2,0){{\bf 2}}
\rput{0}(2,1){{\bf 3}}
\rput{0}(2,2){{\bf 1}}
\rput{0}(2,3){{\bf 3}}
\rput{0}(2,4){{\bf 1}}
\rput{0}(2,5){{\bf 3}}
\rput{0}(3,0){{\bf 1}}
\rput{0}(3,1){{\bf 4}}
\rput{0}(3,2){{\bf 2}}
\rput{0}(3,3){{\bf 4}}
\rput{0}(3,4){{\bf 2}}
\rput{0}(3,5){{\bf 4}}
\rput{0}(4,0){{\bf 2}}
\rput{0}(4,1){{\bf 3}}
\rput{0}(4,2){{\bf 1}}
\rput{0}(4,3){{\bf 3}}
\rput{0}(4,4){{\bf 1}}
\rput{0}(4,5){{\bf 3}}
\rput{0}(5,0){{\bf 1}}
\rput{0}(5,1){{\bf 4}}
\rput{0}(5,2){{\bf 2}}
\rput{0}(5,3){{\bf 4}}
\rput{0}(5,4){{\bf 2}}
\rput{0}(5,5){{\bf 4}}
\endpspicture
\\[2mm]
   (a) &\phantom{(a)} (b) \\[5mm]
%
%
\psset{xunit=22pt}
\psset{yunit=22pt}
\pspicture(-0.5,-0.5)(5.5,5.5)
\psline[linewidth=2pt,linecolor=blue](-0.5,0)(5.5,0)
\psline[linewidth=2pt,linecolor=blue](-0.5,1)(5.5,1)
\psline[linewidth=2pt,linecolor=blue](-0.5,2)(5.5,2)
\psline[linewidth=2pt,linecolor=blue](-0.5,3)(5.5,3)
\psline[linewidth=2pt,linecolor=blue](-0.5,4)(5.5,4)
\psline[linewidth=2pt,linecolor=blue](-0.5,5)(5.5,5)
\psline[linewidth=2pt,linecolor=blue](0,-0.5)(0,5.5)
\psline[linewidth=2pt,linecolor=blue](1,-0.5)(1,5.5)
\psline[linewidth=2pt,linecolor=blue](2,-0.5)(2,5.5)
\psline[linewidth=2pt,linecolor=blue](3,-0.5)(3,5.5)
\psline[linewidth=2pt,linecolor=blue](4,-0.5)(4,5.5)
\psline[linewidth=2pt,linecolor=blue](5,-0.5)(5,5.5)
\psline[linewidth=2pt,linecolor=blue](-0.5,-0.5)(5.5,5.5)
\psline[linewidth=2pt,linecolor=blue](0.5,-0.5)(5.5,4.5)
\psline[linewidth=2pt,linecolor=blue](1.5,-0.5)(5.5,3.5)
\psline[linewidth=2pt,linecolor=blue](2.5,-0.5)(5.5,2.5)
\psline[linewidth=2pt,linecolor=blue](3.5,-0.5)(5.5,1.5)
\psline[linewidth=2pt,linecolor=blue](4.5,-0.5)(5.5,0.5)
\psline[linewidth=2pt,linecolor=blue](-0.5,0.5)(4.5,5.5)
\psline[linewidth=2pt,linecolor=blue](-0.5,1.5)(3.5,5.5)
\psline[linewidth=2pt,linecolor=blue](-0.5,2.5)(2.5,5.5)
\psline[linewidth=2pt,linecolor=blue](-0.5,3.5)(1.5,5.5)
\psline[linewidth=2pt,linecolor=blue](-0.5,4.5)(0.5,5.5)
\multirput{0}(0,0)(0,1){6}{%
  \multirput{0}(0,0)(1,0){6}{%
     \pscircle*[linecolor=white]{8pt}
     \pscircle[linewidth=1pt,linecolor=black] {8pt}
   }
}
\rput{0}(0,0){{\bf 2}}
\rput{0}(0,1){{\bf 4}}
\rput{0}(0,2){{\bf 1}}
\rput{0}(0,3){{\bf 3}}
\rput{0}(0,4){{\bf 1}}
\rput{0}(0,5){{\bf 3}}
\rput{0}(1,0){{\bf 1}}
\rput{0}(1,1){{\bf 3}}
\rput{0}(1,2){{\bf 2}}
\rput{0}(1,3){{\bf 4}}
\rput{0}(1,4){{\bf 2}}
\rput{0}(1,5){{\bf 4}}
\rput{0}(2,0){{\bf 2}}
\rput{0}(2,1){{\bf 4}}
\rput{0}(2,2){{\bf 1}}
\rput{0}(2,3){{\bf 3}}
\rput{0}(2,4){{\bf 1}}
\rput{0}(2,5){{\bf 3}}
\rput{0}(3,0){{\bf 1}}
\rput{0}(3,1){{\bf 3}}
\rput{0}(3,2){{\bf 2}}
\rput{0}(3,3){{\bf 4}}
\rput{0}(3,4){{\bf 2}}
\rput{0}(3,5){{\bf 4}}
\rput{0}(4,0){{\bf 2}}
\rput{0}(4,1){{\bf 4}}
\rput{0}(4,2){{\bf 1}}
\rput{0}(4,3){{\bf 3}}
\rput{0}(4,4){{\bf 1}}
\rput{0}(4,5){{\bf 3}}
\rput{0}(5,0){{\bf 1}}
\rput{0}(5,1){{\bf 3}}
\rput{0}(5,2){{\bf 2}}
\rput{0}(5,3){{\bf 4}}
\rput{0}(5,4){{\bf 2}}
\rput{0}(5,5){{\bf 4}}
\endpspicture
\qquad 
&
\qquad 
%
%
\psset{xunit=22pt}
\psset{yunit=22pt}
\pspicture(-0.5,-0.5)(5.5,5.5)
\psline[linewidth=2pt,linecolor=blue](-0.5,0)(5.5,0)
\psline[linewidth=2pt,linecolor=blue](-0.5,1)(5.5,1)
\psline[linewidth=2pt,linecolor=blue](-0.5,2)(5.5,2)
\psline[linewidth=2pt,linecolor=blue](-0.5,3)(5.5,3)
\psline[linewidth=2pt,linecolor=blue](-0.5,4)(5.5,4)
\psline[linewidth=2pt,linecolor=blue](-0.5,5)(5.5,5)
\psline[linewidth=2pt,linecolor=blue](0,-0.5)(0,5.5)
\psline[linewidth=2pt,linecolor=blue](1,-0.5)(1,5.5)
\psline[linewidth=2pt,linecolor=blue](2,-0.5)(2,5.5)
\psline[linewidth=2pt,linecolor=blue](3,-0.5)(3,5.5)
\psline[linewidth=2pt,linecolor=blue](4,-0.5)(4,5.5)
\psline[linewidth=2pt,linecolor=blue](5,-0.5)(5,5.5)
\psline[linewidth=2pt,linecolor=blue](-0.5,-0.5)(5.5,5.5)
\psline[linewidth=2pt,linecolor=blue](0.5,-0.5)(5.5,4.5)
\psline[linewidth=2pt,linecolor=blue](1.5,-0.5)(5.5,3.5)
\psline[linewidth=2pt,linecolor=blue](2.5,-0.5)(5.5,2.5)
\psline[linewidth=2pt,linecolor=blue](3.5,-0.5)(5.5,1.5)
\psline[linewidth=2pt,linecolor=blue](4.5,-0.5)(5.5,0.5)
\psline[linewidth=2pt,linecolor=blue](-0.5,0.5)(4.5,5.5)
\psline[linewidth=2pt,linecolor=blue](-0.5,1.5)(3.5,5.5)
\psline[linewidth=2pt,linecolor=blue](-0.5,2.5)(2.5,5.5)
\psline[linewidth=2pt,linecolor=blue](-0.5,3.5)(1.5,5.5)
\psline[linewidth=2pt,linecolor=blue](-0.5,4.5)(0.5,5.5)
\multirput{0}(0,0)(0,1){6}{%
  \multirput{0}(0,0)(1,0){6}{%
     \pscircle*[linecolor=white]{8pt}
     \pscircle[linewidth=1pt,linecolor=black] {8pt}
   }
}
\rput{0}(0,0){{\bf 2}}
\rput{0}(0,1){{\bf 3}}
\rput{0}(0,2){{\bf 1}}
\rput{0}(0,3){{\bf 4}}
\rput{0}(0,4){{\bf 1}}
\rput{0}(0,5){{\bf 3}}
\rput{0}(1,0){{\bf 1}}
\rput{0}(1,1){{\bf 4}}
\rput{0}(1,2){{\bf 2}}
\rput{0}(1,3){{\bf 3}}
\rput{0}(1,4){{\bf 2}}
\rput{0}(1,5){{\bf 4}}
\rput{0}(2,0){{\bf 2}}
\rput{0}(2,1){{\bf 3}}
\rput{0}(2,2){{\bf 1}}
\rput{0}(2,3){{\bf 4}}
\rput{0}(2,4){{\bf 1}}
\rput{0}(2,5){{\bf 3}}
\rput{0}(3,0){{\bf 1}}
\rput{0}(3,1){{\bf 4}}
\rput{0}(3,2){{\bf 2}}
\rput{0}(3,3){{\bf 3}}
\rput{0}(3,4){{\bf 2}}
\rput{0}(3,5){{\bf 4}}
\rput{0}(4,0){{\bf 2}}
\rput{0}(4,1){{\bf 3}}
\rput{0}(4,2){{\bf 1}}
\rput{0}(4,3){{\bf 4}}
\rput{0}(4,4){{\bf 1}}
\rput{0}(4,5){{\bf 3}}
\rput{0}(5,0){{\bf 1}}
\rput{0}(5,1){{\bf 4}}
\rput{0}(5,2){{\bf 2}}
\rput{0}(5,3){{\bf 3}}
\rput{0}(5,4){{\bf 2}}
\rput{0}(5,5){{\bf 4}}
\endpspicture
\\[2mm]
(c) & \phantom{(a)} (d) \\[5mm]
\end{tabular}
\caption{\label{figure_tri_L=6}
Four-colorings of the triangulation $T(6,6)$.
(a) Coloring $c_\text{ns}$ \protect\reff{def_coloring_ns} with 
$|\deg(c_\text{ns})|=18$. 
(b) Coloring $f_b$ obtained from $c_\text{ns}$ by swapping colors $1,2$ 
on the bottom row.
(c) Coloring $f_c$ obtained from $f_b$ by swapping colors $3,4$ on the 
second row from the bottom. 
(d) Coloring $f_d$ obtained from $f_b$ by swapping colors $3,4$ on the 
fourth row from the bottom. 
The coloring $c_\text{ns}$ in (a) has $|\deg(c_\text{ns})|=18$; 
the colorings $f_i$ in (b)--(d) have $|\deg(f_i)|=6$. 
}
\end{figure}

There is a single four-coloring $f$ with $|\deg(f)|=18$: 
the non-singular coloring $c_\text{ns}$ depicted in 
Figure~\ref{figure_tri_L=6}(a). 
Each row (horizontal, vertical or diagonal) contains exactly two colors, 
and for any choice of colors $a,b$, the induced subgraph $T_{ab}$ contains 
three parallel connected components, each of them being a cycle of length six. 
Then, the only non-trivial K-changes correspond to swapping colors on one 
of these cycles (as swapping colors simultaneously on two such cycles 
is equivalent to swapping colors on the third cycle and permute
colors $a,b$ globally). 
If we choose colors $1,2$ and swap colors on the bottom row, we get the 
four-coloring $f_b$ with degree $|\deg(f_b)|=6$ depicted in 
Figure~\ref{figure_tri_L=6}(b).
To obtain a new coloring we should choose the other pair of colors $3,4$,
as for any other choice $(a,b)\neq (1,2)$ or $(3,4)$, the induced subgraph 
$T_{ab}$ is connected, so we would not obtain a distinct coloring. 
Again, we only need to consider one of the three horizontal cycles of the 
induced subgraph $T_{34}$. Now we have two different choices: the second or 
the fourth rows from the bottom. The resulting colorings $f_c,f_d$ 
are depicted respectively in Figures~\ref{figure_tri_L=6}(c) and (d).
Both have $|\deg(f_i)|=6$, and all the induced subgraphs $T_{a,b}$ with 
$(a,b)\neq (1,2)$ or $(3,4)$ are again connected. Thus, all these
colorings form a closed class $\mathcal{C}_4^{(1)}$ under K-changes; but 
we still need to prove that there are no additional colorings $f$ with 
$|\deg f|=6$.  
 
To count the number of four-colorings $f$ with $|\deg(f)|=6$ belonging to the 
class  $\mathcal{C}_4^{(1)}$, we can fix the colors of the three vertices 
of a triangular face $t$. Then, all we can do is (for each of the three 
directions -- horizontal, vertical, and diagonal) to swap colors on any
non-empty subset of the four cycles in the chosen direction not intersecting 
$t$. Since there are 15 non-empty subsets, we have $15\times 3=45$ colorings 
$f$ with $|\deg(f)|=6$, and therefore, $|\mathcal{C}_4^{(1)}|=46$. 

Finally, we used a computer program (written in {\sc perl}) that enumerates 
all possible four-colorings on $T(6,6)$ and classify them according to 
$|\deg(f)|$. It finds $305192$ proper four-colorings with zero degree, 
$45$ colorings with $|\deg(f)|=6$, and a single coloring with $|\deg(f)|=18$.
Therefore, $\mathcal{C}^{(1)}_4$ contains all colorings with $|\deg(f)|=6,18$, 
$\mathcal{C}_4(T(6,6)) = \mathcal{C}_4^{(0)} \cup \mathcal{C}^{(1)}_4$, and 
$\Kc(T(6,6),4)=2$. Indeed, the number of all these colorings is equal to 
$P_{T(6,6)}(4)/4! = 305238$. \qed

%
%
%
%


\noindent
{\bf Remark.}  The class $\mathcal{C}^{(0)}_4$ is grossly larger 
than $\mathcal{C}^{(1)}_4$: to be more precise,  
$|\mathcal{C}^{(1)}_4|/|\mathcal{C}^{(0)}_4|\approx 1.5\times 10^{-4}$. 

\medskip

Let us now state a simple lemma which is the basic key in the proof of the
next theorems.

\begin{lemma} \label{lemma.tech}
{\rm (a)} If there is a four-coloring $f$ of the triangulation $T(r,s)$ 
with $\deg(f)\equiv 2 \pmod{4}$, then there exists a four-coloring $g$ of 
$T(3r,3s)$ with $\deg(g)\equiv 6 \pmod{12}$. 

{\rm (b)} If there is a four-coloring $f$ of $T(3r,s)$ or
$T(r,3s)$ with $\deg(f)\equiv 2 \pmod{4}$, then there exists a four-coloring
$g$ of $T(3r,3s)$ with $\deg(g)\equiv 6 \pmod{12}$. 

{\rm (c)} If there is a four-coloring $f$ of the triangulation $T(3r,3s)$ 
with $\deg(f)\equiv 6\pmod{12}$, then for any odd integers $p,q$, 
there exists a four-coloring $g$ of the triangulation $T(3rp,3sq)$ with 
$\deg(g)\equiv 6 \pmod{12}$.
\end{lemma}

\proof
(a) If $f$ is a four-coloring of $T(r,s)$, then we can
obtain a four-coloring $g$ of $T(3r,3s)$ by extending $f$
periodically three times in each direction. If $\deg(f)=2 + 4k$, 
with $k\in\Z$, then 
\be
\deg(g) \;=\; 9 \deg(f) \;=\; 18 + 36 k \;\equiv\; 6 \pmod{12} \,.\nonumber
\ee

(b) The same arguments as in (a) apply here; the only difference is that
the coloring of $T(3r,3s)$ is obtained from the coloring in $T(3r,s)$ 
(resp.\ $T(r,3s)$) by extending periodically the former three times in 
the vertical (resp.\ horizontal) direction. If $\deg(f)=2 + 4k$, then 
the degree of the periodically extended coloring $g$ is  
\be
\deg(g) \;=\; 3 \deg(f) \;=\; 6 + 12 k \;\equiv\; 6 \pmod{12} \,.\nonumber
\ee

(c) If $f$ is a four-coloring of $T(3r,3s)$ with $\deg(f)\equiv 6 \pmod{12}$,
then we can obtain a four-coloring $g$ of $T(3rp,3rq)$ by extending 
$f$ periodically $p$ times in the horizontal direction and $q$ times in the
vertical direction. If $\deg(f)=6+12k$ with $k\in\Z$, the degree of $g$ is
\be
\deg(g) \;=\; p q \deg(f) \;=\; 6pq + 12pqk \;\equiv\; 6 \pmod{12} \,\nonumber
\ee
if both $p$ and $q$ are odd integers. \qed

\subsection{Main results for $\bm{T(3L,3L)}$}

Our main results for triangulations of the type $T(3L,3L)$ can be 
summarized as follows: 

\begin{theorem} \label{theo.main}
For any triangulation $T(3L,3L)$ with $L\geq 2$ there exists a four-coloring
$f$ with $\deg(f)\equiv 6 \pmod{12}$. Hence, $\Kc(T(3L,3L),4)> 1$.
In other words, the WSK dynamics for four-colorings on $T(3L,3L)$ 
is non-ergodic.
\end{theorem}

\proof
The rest of this section is devoted to the proof of Theorem \ref{theo.main}.
We will show that $T(3L,3L)$ admits a four-coloring $f$ with 
$\deg(f)\equiv 6 \pmod{12}$. Then, Corollary~\ref{main.corollary} implies that 
$\Kc(T(3L,3L),4)> 1$ for any $L\geq 2$. The construction of $f$ will depend
on the value of $L$ modulo 4, and we will split the proof in four cases,
$L=4k-2, 4k-1, 4k$, or $L=4k+1$, with $k\in\N$. 

The basic strategy for all
these proofs is to explicitly construct the four-coloring with the 
desired degree. With this aim, it is useful to fix orientations
of both triangulations $T(3L,3L)$ and $\partial\Delta^3$ in order to
compute the degree of a given four-coloring (without ambiguity). 
We orient $T(3L,3L)$ and $\partial\Delta^3$ in such a way that
the boundaries of all triangular faces are always followed clockwise. 
The contribution of a triangular face $t$ of $T(3L,3L)$
to the degree is $+1$ (resp.\ $-1$) if the coloring is $123$ (resp.\ $132$)
if we move clockwise around the boundary of $t$. In our figures, 
those faces with orientation preserved (resp.\ reversed) by $f$ are 
depicted in light (resp.\ dark) gray. 

The easiest case is when $L=4k-2$. In this case, $T(3L,3L)$ admits the
non-singular 4-coloring, whose degree is congruent to 6 modulo 12 by
Proposition~\ref{prop.non-singular}.

\begin{figure}[htb]
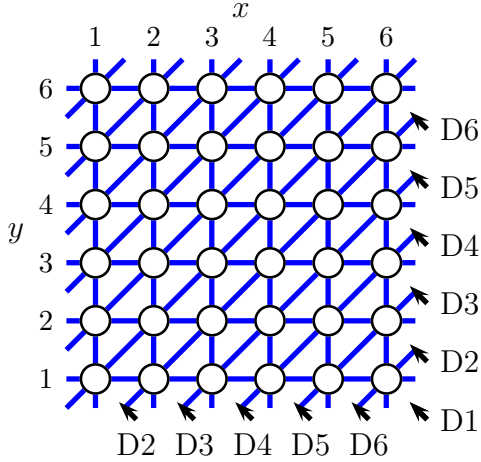

\centering
%
%
\psset{xunit=22pt}
\psset{yunit=22pt}
\pspicture(-1,-1)(6,6)
\psline[linewidth=2pt,linecolor=blue](-0.5,0)(5.5,0)
\psline[linewidth=2pt,linecolor=blue](-0.5,1)(5.5,1)
\psline[linewidth=2pt,linecolor=blue](-0.5,2)(5.5,2)
\psline[linewidth=2pt,linecolor=blue](-0.5,3)(5.5,3)
\psline[linewidth=2pt,linecolor=blue](-0.5,4)(5.5,4)
\psline[linewidth=2pt,linecolor=blue](-0.5,5)(5.5,5)
\psline[linewidth=2pt,linecolor=blue](0,-0.5)(0,5.5)
\psline[linewidth=2pt,linecolor=blue](1,-0.5)(1,5.5)
\psline[linewidth=2pt,linecolor=blue](2,-0.5)(2,5.5)
\psline[linewidth=2pt,linecolor=blue](3,-0.5)(3,5.5)
\psline[linewidth=2pt,linecolor=blue](4,-0.5)(4,5.5)
\psline[linewidth=2pt,linecolor=blue](5,-0.5)(5,5.5)
\psline[linewidth=2pt,linecolor=blue](-0.5,-0.5)(5.5,5.5)
\psline[linewidth=2pt,linecolor=blue](0.5,-0.5)(5.5,4.5)
\psline[linewidth=2pt,linecolor=blue](1.5,-0.5)(5.5,3.5)
\psline[linewidth=2pt,linecolor=blue](2.5,-0.5)(5.5,2.5)
\psline[linewidth=2pt,linecolor=blue](3.5,-0.5)(5.5,1.5)
\psline[linewidth=2pt,linecolor=blue](4.5,-0.5)(5.5,0.5)
\psline[linewidth=2pt,linecolor=blue](-0.5,0.5)(4.5,5.5)
\psline[linewidth=2pt,linecolor=blue](-0.5,1.5)(3.5,5.5)
\psline[linewidth=2pt,linecolor=blue](-0.5,2.5)(2.5,5.5)
\psline[linewidth=2pt,linecolor=blue](-0.5,3.5)(1.5,5.5)
\psline[linewidth=2pt,linecolor=blue](-0.5,4.5)(0.5,5.5)
\multirput{0}(0,0)(0,1){6}{%
  \multirput{0}(0,0)(1,0){6}{%
     \pscircle*[linecolor=white]{6pt}
     \pscircle[linewidth=1pt,linecolor=black] {6pt}
   }
}
%
%
\multirput{0}(5.5,-0.5)(0,1){6}{%
     \psline[linewidth=2pt,linecolor=black]{->}(0.2,-0.2)(-0.1,0.1)
}
\uput[0](5.7,-0.7){D1}
\uput[0](5.7, 0.3){D2}
\uput[0](5.7, 1.3){D3}
\uput[0](5.7, 2.3){D4}
\uput[0](5.7, 3.3){D5}
\uput[0](5.7, 4.3){D6}
\multirput{0}(5.5,-0.5)(-1,0){6}{%
     \psline[linewidth=2pt,linecolor=black]{->}(0.2,-0.2)(-0.1,0.1)
}
\uput[270](4.7,-0.7){D6}
\uput[270](3.7,-0.7){D5}
\uput[270](2.7,-0.7){D4}
\uput[270](1.7,-0.7){D3}
\uput[270](0.7,-0.7){D2}
\uput[90](2.5,6){$x$}
\uput[90](0,5.5){\small 1}
\uput[90](1,5.5){\small 2}
\uput[90](2,5.5){\small 3}
\uput[90](3,5.5){\small 4}
\uput[90](4,5.5){\small 5}
\uput[90](5,5.5){\small 6}
\uput[180](-0.5,0){\small 1}
\uput[180](-0.5,1){\small 2}
\uput[180](-0.5,2){\small 3}
\uput[180](-0.5,3){\small 4}
\uput[180](-0.5,4){\small 5}
\uput[180](-0.5,5){\small 6}
\uput[180](-1,2.5){$y$}
\endpspicture
\caption{\label{figure_tri_notation}
Notation used in the proof of Theorem~\protect\ref{theo.main}.
Given a triangulation $T(M,M)$ (here we depict the case $M=6$),
we label each vertex using Cartesian coordinates $(x,y)$ [$1\leq x,y\leq M$].
The arrows (pointing north-west) show the counter-diagonals D$j$ 
with $j=1,\ldots,M$.
}
\end{figure}

Other cases need a more elaborate
construction. The common strategy is to devise an algorithm to obtain 
the desired four-coloring; and the main ingredient is to use the
counter-diagonals of the triangulations: these counter-diagonals are
orthogonal to the inclined edges of the triangulation when embedded 
in a square grid. They will be denoted as D$j$ with $1\leq j \leq 3L$. 
In Figure~\ref{figure_tri_notation} we show the triangulation $T(6,6)$, and 
its six counter-diagonals D$j$. As we have embedded the triangulation into
a square grid, we will use Cartesian coordinates $(x,y)$, $1\leq x,y\leq 3L$, 
for labelling the vertices.

%
%
%
%

We will describe an algorithm that provides the
desired coloring $f$. It is useful to monitor the degree of the coloring as
we construct it. In particular, at a given step of the algorithm, the 
four-coloring $f$ will be defined on some region $R$ of $T=T(3L,3L)$ (i.e.,
the union of all properly colored triangular faces of $T$). What we mean
by the degree of $f$ at this stage, is the contribution to the degree of $f$
of the triangles belonging to $R$: $\deg(f|_R)$. 
Again, we will count only those triangular faces of $T$ colored $123$. 
Notice that at the end of the algorithm, when $R=T$, this partial degree 
will coincide with the standard one, $\deg(f)=\deg(f|_T)$.  

\medskip

\proofofcase{2}{$L=4k-1$}

Let us consider the triangulation $T=T(12k-3,12k-3)$ with $k\in\N$ (the
case $k=1$ will illustrate our ideas in Figures 
\ref{prop.12k-3.fig1}--\ref{prop.12k-3.fig3-4}). 
Our goal is to obtain a four-coloring 
$f$ of $T$ with degree $\deg(f)\equiv 6 \pmod{12}$. The algorithm to
obtain such a coloring consists of four steps: 

%
%
\begin{figure}[htb]
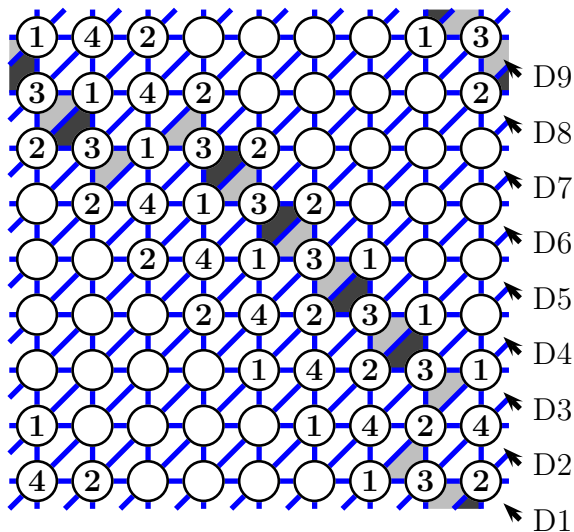

\centering
%
%
\psset{xunit=21pt} 
\psset{yunit=21pt} 
\psset{labelsep=5pt}
\pspicture(-1,-1)(9,9)
\psline*[linecolor=lightgray](1,5)(1,6)(3,6)(3,7)(1,5)
\psline*[linecolor=lightgray](6,0)(7,0)(7,2)(8,2)(6,0)
\psline*[linecolor=lightgray](0,6)(0,7)(1,7)(0,6)
\psline*[linecolor=lightgray](0,8)(-0.5,8)(-0.5,7.5)(0,8)
\psline*[linecolor=lightgray](8,8)(8.5,8)(8.5,7.5)(8,7)(8,8)
\psline*[linecolor=lightgray](3,5)(4,5)(4,6)(3,5)
\psline*[linecolor=lightgray](4,4)(6,4)(5,3)(5,5)(4,4)
\psline*[linecolor=lightgray](6,2)(6,3)(7,3)(6,2)
\psline*[linecolor=lightgray](7,8)(8,8)(8,8.5)(7.5,8.5)(7,8)
\psline*[linecolor=lightgray](7,0)(8,0)(7.5,-0.5)(7,-0.5)(7,0)
\psline*[linecolor=darkgray](8,7)(8.5,7)(8.5,7.5)(8,7)
\psline*[linecolor=darkgray](0,7)(0,8)(-0.5,7.5)(-0.5,7)(0,7)
\psline*[linecolor=darkgray](0,6)(1,6)(1,7)(0,6)
\psline*[linecolor=darkgray](3,5)(3,6)(4,6)(3,5)
\psline*[linecolor=darkgray](4,4)(4,5)(5,5)(4,4)
\psline*[linecolor=darkgray](5,3)(6,3)(6,4)(5,3)
\psline*[linecolor=darkgray](6,2)(7,2)(7,3)(6,2)
\psline*[linecolor=darkgray](7,8)(7,8.5)(7.5,8.5)(7,8)
\psline*[linecolor=darkgray](8,0)(8,-0.5)(7.5,-0.5)(8,0)
\psline[linewidth=2pt,linecolor=blue](-0.5,0)(8.5,0)
\psline[linewidth=2pt,linecolor=blue](-0.5,1)(8.5,1)
\psline[linewidth=2pt,linecolor=blue](-0.5,2)(8.5,2)
\psline[linewidth=2pt,linecolor=blue](-0.5,3)(8.5,3)
\psline[linewidth=2pt,linecolor=blue](-0.5,4)(8.5,4)
\psline[linewidth=2pt,linecolor=blue](-0.5,5)(8.5,5)
\psline[linewidth=2pt,linecolor=blue](-0.5,6)(8.5,6)
\psline[linewidth=2pt,linecolor=blue](-0.5,7)(8.5,7)
\psline[linewidth=2pt,linecolor=blue](-0.5,8)(8.5,8)
\psline[linewidth=2pt,linecolor=blue](0,-0.5)(0,8.5)
\psline[linewidth=2pt,linecolor=blue](1,-0.5)(1,8.5)
\psline[linewidth=2pt,linecolor=blue](2,-0.5)(2,8.5)
\psline[linewidth=2pt,linecolor=blue](3,-0.5)(3,8.5)
\psline[linewidth=2pt,linecolor=blue](4,-0.5)(4,8.5)
\psline[linewidth=2pt,linecolor=blue](5,-0.5)(5,8.5)
\psline[linewidth=2pt,linecolor=blue](6,-0.5)(6,8.5)
\psline[linewidth=2pt,linecolor=blue](7,-0.5)(7,8.5)
\psline[linewidth=2pt,linecolor=blue](8,-0.5)(8,8.5)
\psline[linewidth=2pt,linecolor=blue](-0.5,-0.5)(8.5,8.5)
\psline[linewidth=2pt,linecolor=blue](0.5,-0.5)(8.5,7.5)
\psline[linewidth=2pt,linecolor=blue](1.5,-0.5)(8.5,6.5)
\psline[linewidth=2pt,linecolor=blue](2.5,-0.5)(8.5,5.5)
\psline[linewidth=2pt,linecolor=blue](3.5,-0.5)(8.5,4.5)
\psline[linewidth=2pt,linecolor=blue](4.5,-0.5)(8.5,3.5)
\psline[linewidth=2pt,linecolor=blue](5.5,-0.5)(8.5,2.5)
\psline[linewidth=2pt,linecolor=blue](6.5,-0.5)(8.5,1.5)
\psline[linewidth=2pt,linecolor=blue](7.5,-0.5)(8.5,0.5)
\psline[linewidth=2pt,linecolor=blue](-0.5,0.5)(7.5,8.5)
\psline[linewidth=2pt,linecolor=blue](-0.5,1.5)(6.5,8.5)
\psline[linewidth=2pt,linecolor=blue](-0.5,2.5)(5.5,8.5)
\psline[linewidth=2pt,linecolor=blue](-0.5,3.5)(4.5,8.5)
\psline[linewidth=2pt,linecolor=blue](-0.5,4.5)(3.5,8.5)
\psline[linewidth=2pt,linecolor=blue](-0.5,5.5)(2.5,8.5)
\psline[linewidth=2pt,linecolor=blue](-0.5,6.5)(1.5,8.5)
\psline[linewidth=2pt,linecolor=blue](-0.5,7.5)(0.5,8.5)
\multirput{0}(0,0)(0,1){9}{%
  \multirput{0}(0,0)(1,0){9}{%
     \pscircle*[linecolor=white]{8pt}
     \pscircle[linewidth=1pt,linecolor=black] {8pt}
   }
}
\rput{0}(0,0){{\bf 4}}
\rput{0}(0,1){{\bf 1}}
\rput{0}(0,2){{\bf  }}
\rput{0}(0,3){{\bf  }}
\rput{0}(0,4){{\bf  }}
\rput{0}(0,5){{\bf  }}
\rput{0}(0,6){{\bf 2}}
\rput{0}(0,7){{\bf 3}}
\rput{0}(0,8){{\bf 1}}
\rput{0}(1,0){{\bf 2}}
\rput{0}(1,1){{\bf  }}
\rput{0}(1,2){{\bf  }}
\rput{0}(1,3){{\bf  }}
\rput{0}(1,4){{\bf  }}
\rput{0}(1,5){{\bf 2}}
\rput{0}(1,6){{\bf 3}}
\rput{0}(1,7){{\bf 1}}
\rput{0}(1,8){{\bf 4}}
\rput{0}(2,0){{\bf  }}
\rput{0}(2,1){{\bf  }}
\rput{0}(2,2){{\bf  }}
\rput{0}(2,3){{\bf  }}
\rput{0}(2,4){{\bf 2}}
\rput{0}(2,5){{\bf 4}}
\rput{0}(2,6){{\bf 1}}
\rput{0}(2,7){{\bf 4}}
\rput{0}(2,8){{\bf 2}}
\rput{0}(3,0){{\bf  }}
\rput{0}(3,1){{\bf  }}
\rput{0}(3,2){{\bf  }}
\rput{0}(3,3){{\bf 2}}
\rput{0}(3,4){{\bf 4}}
\rput{0}(3,5){{\bf 1}}
\rput{0}(3,6){{\bf 3}}
\rput{0}(3,7){{\bf 2}}
\rput{0}(3,8){{\bf  }}
\rput{0}(4,0){{\bf  }}
\rput{0}(4,1){{\bf  }}
\rput{0}(4,2){{\bf 1}}
\rput{0}(4,3){{\bf 4}}
\rput{0}(4,4){{\bf 1}}
\rput{0}(4,5){{\bf 3}}
\rput{0}(4,6){{\bf 2}}
\rput{0}(4,7){{\bf  }}
\rput{0}(4,8){{\bf  }}
\rput{0}(5,0){{\bf  }}
\rput{0}(5,1){{\bf 1}}
\rput{0}(5,2){{\bf 4}}
\rput{0}(5,3){{\bf 2}}
\rput{0}(5,4){{\bf 3}}
\rput{0}(5,5){{\bf 2}}
\rput{0}(5,6){{\bf  }}
\rput{0}(5,7){{\bf  }}
\rput{0}(5,8){{\bf  }}
\rput{0}(6,0){{\bf 1}}
\rput{0}(6,1){{\bf 4}}
\rput{0}(6,2){{\bf 2}}
\rput{0}(6,3){{\bf 3}}
\rput{0}(6,4){{\bf 1}}
\rput{0}(6,5){{\bf  }}
\rput{0}(6,6){{\bf  }}
\rput{0}(6,7){{\bf  }}
\rput{0}(6,8){{\bf  }}
\rput{0}(7,0){{\bf 3}}
\rput{0}(7,1){{\bf 2}}
\rput{0}(7,2){{\bf 3}}
\rput{0}(7,3){{\bf 1}}
\rput{0}(7,4){{\bf  }}
\rput{0}(7,5){{\bf  }}
\rput{0}(7,6){{\bf  }}
\rput{0}(7,7){{\bf  }}
\rput{0}(7,8){{\bf 1}}
\rput{0}(8,0){{\bf 2}}
\rput{0}(8,1){{\bf 4}}
\rput{0}(8,2){{\bf 1}}
\rput{0}(8,3){{\bf  }}
\rput{0}(8,4){{\bf  }}
\rput{0}(8,5){{\bf  }}
\rput{0}(8,6){{\bf  }}
\rput{0}(8,7){{\bf 2}}
\rput{0}(8,8){{\bf 3}}
%
%
\multirput{0}(8.5,-0.5)(0,1){9}{%
     \psline[linewidth=2pt,linecolor=black]{->}(0.2,-0.2)(-0.1,0.1) 
}
\uput[0](8.7,-0.7){D1}
\uput[0](8.7, 0.3){D2}
\uput[0](8.7, 1.3){D3}
\uput[0](8.7, 2.3){D4}
\uput[0](8.7, 3.3){D5}
\uput[0](8.7, 4.3){D6}
\uput[0](8.7, 5.3){D7}
\uput[0](8.7, 6.3){D8}
\uput[0](8.7, 7.3){D9}
\endpspicture
\caption{ \label{prop.12k-3.fig1}
The 4-coloring of $T(9,9)$ after Step~1 in the proof of the case $L=4k-1$.}
\end{figure}
%
%

\noindent 
{\bf Step 1.} 
We start by coloring the counter-diagonal D1: we color $1$ the vertices with 
$x$-coordinates $1\leq x \leq 6k-1$; the other $6k-2$ vertices are colored $2$.

On D2, we color $3$ those $6k-1$ vertices with $x$-coordinates 
$3k+1\leq x \leq 9k-1$. The other vertices on D2 are colored $4$. The 
vertices on D$(12k-3)$ are colored $3$ or $4$ in such a 
way that the resulting coloring is proper (for each vertex, there is a 
unique choice).

On D3 and D$(12k-4)$, we color all vertices $1$ or $2$ (there is a 
unique choice for each vertex). The resulting coloring is depicted on 
Figure~\ref{prop.12k-3.fig1}. The partial degree of $f$ is $\deg f|_R = 4$. 

\medskip
\noindent
{\bf Step 2.}
For $k>1$, we find that there are $12k-8$ counter-diagonals to be colored
and we need to sequentially color all of them but four. This can be 
achieved by performing the following procedure: suppose that we
have already colored counter-diagonals D$j$ and D$(12k-j-1)$ ($j\geq 3$) 
using colors $1$ and $2$. Then, we color D$(j+1)$ and D$(12k-j-2)$ using 
colors $3$ and $4$, and D$(j+2)$ and D$(12k-j-3)$ using colors $1$ and $2$. 
As in Step~1, for each vertex there is a unique choice. 

This procedure is repeated $3(k-1)$ times, so we add $12(k-1)$ 
counter-diagonals, and there are only four counter-diagonals not yet 
colored. Indeed, the last colored counter-diagonals D$(6k-3)$ and 
D$(6k+2)$ have colors $1$ and $2$, the same as it was at the end of Step~1. 

Each of these $3(k-1)$ steps adds $4$ to the degree of the coloring.
Thus, the partial degree of $f$ is $\deg f|_R = 4 + 12(k-1)$.  

\medskip
\noindent
{\bf Step 3.}
There remain only four counter-diagonals to be colored: D$(6k-2)$, D$(6k-1)$,
D$(6k)$, and D$(6k+1)$.
On D$(6k-2)$, the vertices $(3k-1,3k-1)$ and $(9k-2,9k-3)$ only admit  
a single color (which is $3$ for one of them, and $4$ for the other one).  
The rest of the vertices on D$(6k-2)$ are colored $1$ and $2$ (again, there 
is a unique choice for each vertex). 

We now color $3$ or $4$ all the vertices on D$(6k+1)$ (the choice is again 
unique for each vertex). The resulting coloring is depicted in 
Figure~\ref{prop.12k-3.fig3-4}(a). 
The contribution to the partial degree of the new triangles is zero; the 
partial degree of $f$ is given by $\deg f|_R = 4 + 12(k-1)$.

\medskip
\noindent
{\bf Step 4.}
On D$(6k-1)$, there are two pairs of nearby vertices which only admit a 
single color (which is $3$ for one pair, and $4$ for the other one). 
These vertices are located at $(3k-1,3k)$, $(3k,3k-1)$, $(9k-1,9k-3)$, and
$(9k-2,9k-2)$. The other vertices on D$(6k-1)$ can be colored $3$ or $4$
(with only one choice for each of them). The increment of the degree after 
coloring these vertices is $-2$, thus $\deg f|_R = 2 + 12(k-1)$. 

Finally, all vertices on D$(6k)$ are colored $1$ and $2$; and again
the choice is unique for each vertex. The final coloring is depicted on
Figure~\ref{prop.12k-3.fig3-4}(b). The increment in the degree is $4$, 
and therefore, the degree of the four-coloring $f$ is 
\begin{equation}
\deg f \;=\; 6 + 12(k-1) \;\equiv\; 6 \pmod{12} 
\end{equation}
This coloring $f$ of $T(12k-3,12k-3)$ satisfies the two needed properties:
it is a proper coloring and its degree is congruent to six modulo $12$.

\proofofcase{3}{$L=4k$}

Let us consider the triangulation $T=T(12k,12k)$ with $k\in\N$ (we will
illustrate the main steps with the case $k=1$). Our algorithm consists of
five steps: 

\medskip
\noindent
{\bf Step 1.}
On the counter-diagonal D1 we color $1$ the $6k$ consecutive vertices with 
$x$-coordinates $1\leq x \leq 6k$. The other $6k$ vertices on D1 are colored 
$2$.

\clearpage
%
%
\begin{figure}[htb]
\centering
\begin{tabular}{c}
%
%
\psset{xunit=21pt}  
\psset{yunit=21pt}  
\psset{labelsep=5pt}
\pspicture(-1,-1)(9,9)
\psline*[linecolor=lightgray](1,5)(1,6)(3,6)(3,7)(1,5)
\psline*[linecolor=lightgray](6,0)(7,0)(7,2)(8,2)(6,0)
\psline*[linecolor=lightgray](0,6)(0,7)(1,7)(0,6)
\psline*[linecolor=lightgray](0,8)(-0.5,8)(-0.5,7.5)(0,8)
\psline*[linecolor=lightgray](8,8)(8.5,8)(8.5,7.5)(8,7)(8,8)
\psline*[linecolor=lightgray](3,5)(4,5)(4,6)(3,5)
\psline*[linecolor=lightgray](4,4)(6,4)(5,3)(5,5)(4,4)
\psline*[linecolor=lightgray](6,2)(6,3)(7,3)(6,2)
\psline*[linecolor=lightgray](7,8)(8,8)(8,8.5)(7.5,8.5)(7,8)
\psline*[linecolor=lightgray](7,0)(8,0)(7.5,-0.5)(7,-0.5)(7,0)
\psline*[linecolor=lightgray](3,6)(4,6)(4,7)(3,6)
\psline*[linecolor=lightgray](4,5)(5,5)(5,6)(4,5)
\psline*[linecolor=lightgray](6,3)(6,4)(7,4)(6,3)
\psline*[linecolor=lightgray](7,2)(7,3)(8,3)(7,2)
\psline*[linecolor=darkgray](8,7)(8.5,7)(8.5,7.5)(8,7)
\psline*[linecolor=darkgray](0,7)(0,8)(-0.5,7.5)(-0.5,7)(0,7)
\psline*[linecolor=darkgray](0,6)(1,6)(1,7)(0,6)
\psline*[linecolor=darkgray](3,5)(3,6)(4,6)(3,5)
\psline*[linecolor=darkgray](4,4)(4,5)(5,5)(4,4)
\psline*[linecolor=darkgray](5,3)(6,3)(6,4)(5,3)
\psline*[linecolor=darkgray](6,2)(7,2)(7,3)(6,2)
\psline*[linecolor=darkgray](7,8)(7,8.5)(7.5,8.5)(7,8)
\psline*[linecolor=darkgray](8,0)(8,-0.5)(7.5,-0.5)(8,0)
\psline*[linecolor=darkgray](3,6)(3,7)(4,7)(3,6)
\psline*[linecolor=darkgray](4,5)(4,6)(5,6)(4,5)
\psline*[linecolor=darkgray](6,3)(7,3)(7,4)(6,3)
\psline*[linecolor=darkgray](7,2)(8,2)(8,3)(7,2)
\psline[linewidth=2pt,linecolor=blue](-0.5,0)(8.5,0)
\psline[linewidth=2pt,linecolor=blue](-0.5,1)(8.5,1)
\psline[linewidth=2pt,linecolor=blue](-0.5,2)(8.5,2)
\psline[linewidth=2pt,linecolor=blue](-0.5,3)(8.5,3)
\psline[linewidth=2pt,linecolor=blue](-0.5,4)(8.5,4)
\psline[linewidth=2pt,linecolor=blue](-0.5,5)(8.5,5)
\psline[linewidth=2pt,linecolor=blue](-0.5,6)(8.5,6)
\psline[linewidth=2pt,linecolor=blue](-0.5,7)(8.5,7)
\psline[linewidth=2pt,linecolor=blue](-0.5,8)(8.5,8)
\psline[linewidth=2pt,linecolor=blue](0,-0.5)(0,8.5)
\psline[linewidth=2pt,linecolor=blue](1,-0.5)(1,8.5)
\psline[linewidth=2pt,linecolor=blue](2,-0.5)(2,8.5)
\psline[linewidth=2pt,linecolor=blue](3,-0.5)(3,8.5)
\psline[linewidth=2pt,linecolor=blue](4,-0.5)(4,8.5)
\psline[linewidth=2pt,linecolor=blue](5,-0.5)(5,8.5)
\psline[linewidth=2pt,linecolor=blue](6,-0.5)(6,8.5)
\psline[linewidth=2pt,linecolor=blue](7,-0.5)(7,8.5)
\psline[linewidth=2pt,linecolor=blue](8,-0.5)(8,8.5)
\psline[linewidth=2pt,linecolor=blue](-0.5,-0.5)(8.5,8.5)
\psline[linewidth=2pt,linecolor=blue](0.5,-0.5)(8.5,7.5)
\psline[linewidth=2pt,linecolor=blue](1.5,-0.5)(8.5,6.5)
\psline[linewidth=2pt,linecolor=blue](2.5,-0.5)(8.5,5.5)
\psline[linewidth=2pt,linecolor=blue](3.5,-0.5)(8.5,4.5)
\psline[linewidth=2pt,linecolor=blue](4.5,-0.5)(8.5,3.5)
\psline[linewidth=2pt,linecolor=blue](5.5,-0.5)(8.5,2.5)
\psline[linewidth=2pt,linecolor=blue](6.5,-0.5)(8.5,1.5)
\psline[linewidth=2pt,linecolor=blue](7.5,-0.5)(8.5,0.5)
\psline[linewidth=2pt,linecolor=blue](-0.5,0.5)(7.5,8.5)
\psline[linewidth=2pt,linecolor=blue](-0.5,1.5)(6.5,8.5)
\psline[linewidth=2pt,linecolor=blue](-0.5,2.5)(5.5,8.5)
\psline[linewidth=2pt,linecolor=blue](-0.5,3.5)(4.5,8.5)
\psline[linewidth=2pt,linecolor=blue](-0.5,4.5)(3.5,8.5)
\psline[linewidth=2pt,linecolor=blue](-0.5,5.5)(2.5,8.5)
\psline[linewidth=2pt,linecolor=blue](-0.5,6.5)(1.5,8.5)
\psline[linewidth=2pt,linecolor=blue](-0.5,7.5)(0.5,8.5)
\multirput{0}(0,0)(0,1){9}{%
  \multirput{0}(0,0)(1,0){9}{%
     \pscircle*[linecolor=white]{8pt}
     \pscircle[linewidth=1pt,linecolor=black] {8pt}
   }
}
\rput{0}(0,0){{\bf 4}}
\rput{0}(0,1){{\bf 1}}
\rput{0}(0,2){{\bf 2}}
\rput{0}(0,3){{\bf  }}
\rput{0}(0,4){{\bf  }}
\rput{0}(0,5){{\bf 4}}
\rput{0}(0,6){{\bf 2}}
\rput{0}(0,7){{\bf 3}}
\rput{0}(0,8){{\bf 1}}
\rput{0}(1,0){{\bf 2}}
\rput{0}(1,1){{\bf 3}}
\rput{0}(1,2){{\bf  }}
\rput{0}(1,3){{\bf  }}
\rput{0}(1,4){{\bf 3}}
\rput{0}(1,5){{\bf 2}}
\rput{0}(1,6){{\bf 3}}
\rput{0}(1,7){{\bf 1}}
\rput{0}(1,8){{\bf 4}}
\rput{0}(2,0){{\bf 1}}
\rput{0}(2,1){{\bf  }}
\rput{0}(2,2){{\bf  }}
\rput{0}(2,3){{\bf 3}}
\rput{0}(2,4){{\bf 2}}
\rput{0}(2,5){{\bf 4}}
\rput{0}(2,6){{\bf 1}}
\rput{0}(2,7){{\bf 4}}
\rput{0}(2,8){{\bf 2}}
\rput{0}(3,0){{\bf  }}
\rput{0}(3,1){{\bf  }}
\rput{0}(3,2){{\bf 3}}
\rput{0}(3,3){{\bf 2}}
\rput{0}(3,4){{\bf 4}}
\rput{0}(3,5){{\bf 1}}
\rput{0}(3,6){{\bf 3}}
\rput{0}(3,7){{\bf 2}}
\rput{0}(3,8){{\bf 1}}
\rput{0}(4,0){{\bf  }}
\rput{0}(4,1){{\bf 3}}
\rput{0}(4,2){{\bf 1}}
\rput{0}(4,3){{\bf 4}}
\rput{0}(4,4){{\bf 1}}
\rput{0}(4,5){{\bf 3}}
\rput{0}(4,6){{\bf 2}}
\rput{0}(4,7){{\bf 1}}
\rput{0}(4,8){{\bf  }}
\rput{0}(5,0){{\bf 3}}
\rput{0}(5,1){{\bf 1}}
\rput{0}(5,2){{\bf 4}}
\rput{0}(5,3){{\bf 2}}
\rput{0}(5,4){{\bf 3}}
\rput{0}(5,5){{\bf 2}}
\rput{0}(5,6){{\bf 1}}
\rput{0}(5,7){{\bf  }}
\rput{0}(5,8){{\bf  }}
\rput{0}(6,0){{\bf 1}}
\rput{0}(6,1){{\bf 4}}
\rput{0}(6,2){{\bf 2}}
\rput{0}(6,3){{\bf 3}}
\rput{0}(6,4){{\bf 1}}
\rput{0}(6,5){{\bf 4}}
\rput{0}(6,6){{\bf  }}
\rput{0}(6,7){{\bf  }}
\rput{0}(6,8){{\bf 4}}
\rput{0}(7,0){{\bf 3}}
\rput{0}(7,1){{\bf 2}}
\rput{0}(7,2){{\bf 3}}
\rput{0}(7,3){{\bf 1}}
\rput{0}(7,4){{\bf 2}}
\rput{0}(7,5){{\bf  }}
\rput{0}(7,6){{\bf  }}
\rput{0}(7,7){{\bf 4}}
\rput{0}(7,8){{\bf 1}}
\rput{0}(8,0){{\bf 2}}
\rput{0}(8,1){{\bf 4}}
\rput{0}(8,2){{\bf 1}}
\rput{0}(8,3){{\bf 2}}
\rput{0}(8,4){{\bf  }}
\rput{0}(8,5){{\bf  }}
\rput{0}(8,6){{\bf 4}}
\rput{0}(8,7){{\bf 2}}
\rput{0}(8,8){{\bf 3}}
%
%
\multirput{0}(8.5,-0.5)(0,1){9}{%
     \psline[linewidth=2pt,linecolor=black]{->}(0.2,-0.2)(-0.1,0.1) 
}
\uput[0](8.7,-0.7){D1}
\uput[0](8.7, 0.3){D2}
\uput[0](8.7, 1.3){D3}
\uput[0](8.7, 2.3){D4}
\uput[0](8.7, 3.3){D5}
\uput[0](8.7, 4.3){D6}
\uput[0](8.7, 5.3){D7}
\uput[0](8.7, 6.3){D8}
\uput[0](8.7, 7.3){D9}
\endpspicture
\\[1mm]
(a) 
\\
%
%
\psset{xunit=21pt} 
\psset{yunit=21pt} 
\psset{labelsep=5pt}
\pspicture(-1,-1)(9,9)
\psline*[linecolor=lightgray](1,5)(1,6)(3,6)(3,7)(1,5)
\psline*[linecolor=lightgray](6,0)(7,0)(7,2)(8,2)(6,0)
\psline*[linecolor=lightgray](0,6)(0,7)(1,7)(0,6)
\psline*[linecolor=lightgray](0,8)(-0.5,8)(-0.5,7.5)(0,8)
\psline*[linecolor=lightgray](8,8)(8.5,8)(8.5,7.5)(8,7)(8,8)
\psline*[linecolor=lightgray](3,5)(4,5)(4,6)(3,5)
\psline*[linecolor=lightgray](4,4)(6,4)(5,3)(5,5)(4,4)
\psline*[linecolor=lightgray](6,2)(6,3)(7,3)(6,2)
\psline*[linecolor=lightgray](7,0)(8,0)(7.5,-0.5)(7,-0.5)(7,0)
\psline*[linecolor=lightgray](3,6)(4,6)(4,7)(3,6)
\psline*[linecolor=lightgray](4,5)(5,5)(5,6)(4,5)
\psline*[linecolor=lightgray](6,3)(6,4)(7,4)(6,3)
\psline*[linecolor=lightgray](7,2)(7,3)(8,3)(7,2)
\psline*[linecolor=lightgray](4,6)(5,6)(5,7)(4,6)
\psline*[linecolor=lightgray](7,3)(7,4)(8,4)(7,3)
\psline*[linecolor=lightgray](0,4)(1,4)(1,5)(0,4)
\psline*[linecolor=lightgray](1,3)(2,3)(2,4)(1,3)
\psline*[linecolor=lightgray](2,2)(3,3)(3,1)(4,2)(2,2)
\psline*[linecolor=lightgray](4,0)(4,1)(5,1)(4,0)
\psline*[linecolor=lightgray](5,0)(6,0)(5.5,-0.5)(5,-0.5)(5,0)
\psline*[linecolor=lightgray](5,8)(5,8.5)(5.5,8.5)(5,8)
\psline*[linecolor=lightgray](5,8)(5,7)(4,7)(5,8)
\psline*[linecolor=lightgray](6,7)(6,6)(5,6)(6,7)
\psline*[linecolor=lightgray](4,7)(5,7)(5,8.5)(5.5,8.5)(4,7)
\psline*[linecolor=lightgray](7,4)(7,5)(8,5)(7,4)
\psline*[linecolor=lightgray](8,3)(8,4)(8.5,4)(8.5,3.5)(8,3)
\psline*[linecolor=lightgray](0,4)(-0.5,4)(-0.5,3.5)(0,4)
\psline*[linecolor=lightgray](7,8)(7,8.5)(7.5,8.5)(7,8)
\psline*[linecolor=darkgray](8,7)(8.5,7)(8.5,7.5)(8,7)
\psline*[linecolor=darkgray](0,7)(0,8)(-0.5,7.5)(-0.5,7)(0,7)
\psline*[linecolor=darkgray](0,6)(1,6)(1,7)(0,6)
\psline*[linecolor=darkgray](3,5)(3,6)(4,6)(3,5)
\psline*[linecolor=darkgray](4,4)(4,5)(5,5)(4,4)
\psline*[linecolor=darkgray](5,3)(6,3)(6,4)(5,3)
\psline*[linecolor=darkgray](6,2)(7,2)(7,3)(6,2)
\psline*[linecolor=darkgray](8,0)(8,-0.5)(7.5,-0.5)(8,0)
\psline*[linecolor=darkgray](3,6)(3,7)(4,7)(3,6)
\psline*[linecolor=darkgray](4,5)(4,6)(5,6)(4,5)
\psline*[linecolor=darkgray](6,3)(7,3)(7,4)(6,3)
\psline*[linecolor=darkgray](7,2)(8,2)(8,3)(7,2)
\psline*[linecolor=darkgray](4,6)(4,7)(5,7)(4,6)
\psline*[linecolor=darkgray](5,5)(5,6)(6,6)(5,5)
\psline*[linecolor=darkgray](6,4)(7,4)(7,5)(6,4)
\psline*[linecolor=darkgray](7,3)(8,3)(8,4)(7,3)
\psline*[linecolor=darkgray](1,3)(1,4)(2,4)(1,3)
\psline*[linecolor=darkgray](2,2)(2,3)(3,3)(2,2)
\psline*[linecolor=darkgray](3,1)(4,1)(4,2)(3,1)
\psline*[linecolor=darkgray](4,0)(5,0)(5,1)(4,0)
\psline*[linecolor=darkgray](5,6)(5,7)(6,7)(5,6)
\psline*[linecolor=darkgray](7,4)(8,4)(8,5)(7,4)
\psline*[linecolor=darkgray](7,8)(8,8)(8,8.5)(7.5,8.5)(7,8)
\psline[linewidth=2pt,linecolor=blue](-0.5,0)(8.5,0)
\psline[linewidth=2pt,linecolor=blue](-0.5,1)(8.5,1)
\psline[linewidth=2pt,linecolor=blue](-0.5,2)(8.5,2)
\psline[linewidth=2pt,linecolor=blue](-0.5,3)(8.5,3)
\psline[linewidth=2pt,linecolor=blue](-0.5,4)(8.5,4)
\psline[linewidth=2pt,linecolor=blue](-0.5,5)(8.5,5)
\psline[linewidth=2pt,linecolor=blue](-0.5,6)(8.5,6)
\psline[linewidth=2pt,linecolor=blue](-0.5,7)(8.5,7)
\psline[linewidth=2pt,linecolor=blue](-0.5,8)(8.5,8)
\psline[linewidth=2pt,linecolor=blue](0,-0.5)(0,8.5)
\psline[linewidth=2pt,linecolor=blue](1,-0.5)(1,8.5)
\psline[linewidth=2pt,linecolor=blue](2,-0.5)(2,8.5)
\psline[linewidth=2pt,linecolor=blue](3,-0.5)(3,8.5)
\psline[linewidth=2pt,linecolor=blue](4,-0.5)(4,8.5)
\psline[linewidth=2pt,linecolor=blue](5,-0.5)(5,8.5)
\psline[linewidth=2pt,linecolor=blue](6,-0.5)(6,8.5)
\psline[linewidth=2pt,linecolor=blue](7,-0.5)(7,8.5)
\psline[linewidth=2pt,linecolor=blue](8,-0.5)(8,8.5)
\psline[linewidth=2pt,linecolor=blue](-0.5,-0.5)(8.5,8.5)
\psline[linewidth=2pt,linecolor=blue](0.5,-0.5)(8.5,7.5)
\psline[linewidth=2pt,linecolor=blue](1.5,-0.5)(8.5,6.5)
\psline[linewidth=2pt,linecolor=blue](2.5,-0.5)(8.5,5.5)
\psline[linewidth=2pt,linecolor=blue](3.5,-0.5)(8.5,4.5)
\psline[linewidth=2pt,linecolor=blue](4.5,-0.5)(8.5,3.5)
\psline[linewidth=2pt,linecolor=blue](5.5,-0.5)(8.5,2.5)
\psline[linewidth=2pt,linecolor=blue](6.5,-0.5)(8.5,1.5)
\psline[linewidth=2pt,linecolor=blue](7.5,-0.5)(8.5,0.5)
\psline[linewidth=2pt,linecolor=blue](-0.5,0.5)(7.5,8.5)
\psline[linewidth=2pt,linecolor=blue](-0.5,1.5)(6.5,8.5)
\psline[linewidth=2pt,linecolor=blue](-0.5,2.5)(5.5,8.5)
\psline[linewidth=2pt,linecolor=blue](-0.5,3.5)(4.5,8.5)
\psline[linewidth=2pt,linecolor=blue](-0.5,4.5)(3.5,8.5)
\psline[linewidth=2pt,linecolor=blue](-0.5,5.5)(2.5,8.5)
\psline[linewidth=2pt,linecolor=blue](-0.5,6.5)(1.5,8.5)
\psline[linewidth=2pt,linecolor=blue](-0.5,7.5)(0.5,8.5)
\multirput{0}(0,0)(0,1){9}{%
  \multirput{0}(0,0)(1,0){9}{%
     \pscircle*[linecolor=white]{8pt}
     \pscircle[linewidth=1pt,linecolor=black] {8pt}
   }
}
\rput{0}(0,0){{\bf 4}}
\rput{0}(0,1){{\bf 1}}
\rput{0}(0,2){{\bf 2}}
\rput{0}(0,3){{\bf 4}}
\rput{0}(0,4){{\bf 1}}
\rput{0}(0,5){{\bf 4}}
\rput{0}(0,6){{\bf 2}}
\rput{0}(0,7){{\bf 3}}
\rput{0}(0,8){{\bf 1}}
\rput{0}(1,0){{\bf 2}}
\rput{0}(1,1){{\bf 3}}
\rput{0}(1,2){{\bf 4}}
\rput{0}(1,3){{\bf 1}}
\rput{0}(1,4){{\bf 3}}
\rput{0}(1,5){{\bf 2}}
\rput{0}(1,6){{\bf 3}}
\rput{0}(1,7){{\bf 1}}
\rput{0}(1,8){{\bf 4}}
\rput{0}(2,0){{\bf 1}}
\rput{0}(2,1){{\bf 4}}
\rput{0}(2,2){{\bf 1}}
\rput{0}(2,3){{\bf 3}}
\rput{0}(2,4){{\bf 2}}
\rput{0}(2,5){{\bf 4}}
\rput{0}(2,6){{\bf 1}}
\rput{0}(2,7){{\bf 4}}
\rput{0}(2,8){{\bf 2}}
\rput{0}(3,0){{\bf 4}}
\rput{0}(3,1){{\bf 2}}
\rput{0}(3,2){{\bf 3}}
\rput{0}(3,3){{\bf 2}}
\rput{0}(3,4){{\bf 4}}
\rput{0}(3,5){{\bf 1}}
\rput{0}(3,6){{\bf 3}}
\rput{0}(3,7){{\bf 2}}
\rput{0}(3,8){{\bf 1}}
\rput{0}(4,0){{\bf 2}}
\rput{0}(4,1){{\bf 3}}
\rput{0}(4,2){{\bf 1}}
\rput{0}(4,3){{\bf 4}}
\rput{0}(4,4){{\bf 1}}
\rput{0}(4,5){{\bf 3}}
\rput{0}(4,6){{\bf 2}}
\rput{0}(4,7){{\bf 1}}
\rput{0}(4,8){{\bf 4}}
\rput{0}(5,0){{\bf 3}}
\rput{0}(5,1){{\bf 1}}
\rput{0}(5,2){{\bf 4}}
\rput{0}(5,3){{\bf 2}}
\rput{0}(5,4){{\bf 3}}
\rput{0}(5,5){{\bf 2}}
\rput{0}(5,6){{\bf 1}}
\rput{0}(5,7){{\bf 3}}
\rput{0}(5,8){{\bf 2}}
\rput{0}(6,0){{\bf 1}}
\rput{0}(6,1){{\bf 4}}
\rput{0}(6,2){{\bf 2}}
\rput{0}(6,3){{\bf 3}}
\rput{0}(6,4){{\bf 1}}
\rput{0}(6,5){{\bf 4}}
\rput{0}(6,6){{\bf 3}}
\rput{0}(6,7){{\bf 2}}
\rput{0}(6,8){{\bf 4}}
\rput{0}(7,0){{\bf 3}}
\rput{0}(7,1){{\bf 2}}
\rput{0}(7,2){{\bf 3}}
\rput{0}(7,3){{\bf 1}}
\rput{0}(7,4){{\bf 2}}
\rput{0}(7,5){{\bf 3}}
\rput{0}(7,6){{\bf 1}}
\rput{0}(7,7){{\bf 4}}
\rput{0}(7,8){{\bf 1}}
\rput{0}(8,0){{\bf 2}}
\rput{0}(8,1){{\bf 4}}
\rput{0}(8,2){{\bf 1}}
\rput{0}(8,3){{\bf 2}}
\rput{0}(8,4){{\bf 3}}
\rput{0}(8,5){{\bf 1}}
\rput{0}(8,6){{\bf 4}}
\rput{0}(8,7){{\bf 2}}
\rput{0}(8,8){{\bf 3}}
%
%
\multirput{0}(8.5,-0.5)(0,1){9}{%
     \psline[linewidth=2pt,linecolor=black]{->}(0.2,-0.2)(-0.1,0.1) 
}
\uput[0](8.7,-0.7){D1}
\uput[0](8.7, 0.3){D2}
\uput[0](8.7, 1.3){D3}
\uput[0](8.7, 2.3){D4}
\uput[0](8.7, 3.3){D5}
\uput[0](8.7, 4.3){D6}
\uput[0](8.7, 5.3){D7}
\uput[0](8.7, 6.3){D8}
\uput[0](8.7, 7.3){D9}
\endpspicture
\\[1mm]
(b) 
\end{tabular}
\caption{\label{prop.12k-3.fig3-4}
Four-colorings of the triangulation $T(9,9)$ after Steps~3 (a) and~4 (b) 
in the proof of the case $L=4k-1$.
}
\end{figure}
%
%

%
%
\begin{figure}[htb]
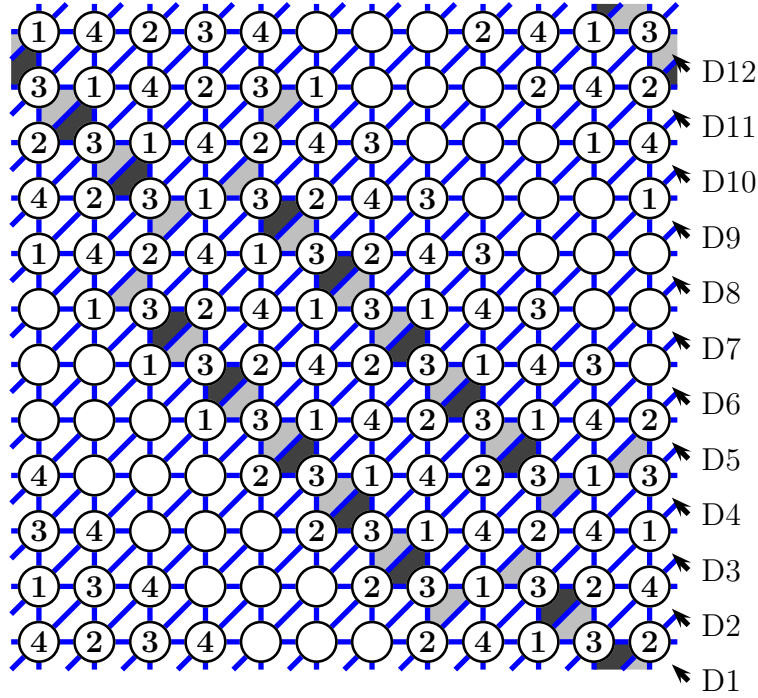

\centering
%
%
\psset{xunit=21pt}
\psset{yunit=21pt}
\psset{labelsep=5pt}
\pspicture(-1,-1)(13,12)
\psline*[linecolor=lightgray](0,9)(0,10)(1,10)(0,9)
\psline*[linecolor=lightgray](1,8)(1,9)(2,9)(1,8)
\psline*[linecolor=lightgray](2,7)(2,8)(4,8)(4,9)(2,7)
\psline*[linecolor=lightgray](4,7)(5,7)(5,8)(4,7)
\psline*[linecolor=lightgray](5,6)(7,6)(6,5)(6,7)(5,6)
\psline*[linecolor=lightgray](7,5)(8,5)(7,4)(7,5)
\psline*[linecolor=lightgray](8,4)(9,4)(8,3)(8,4)
\psline*[linecolor=lightgray](9,3)(10,3)(9,2)(9,3)
\psline*[linecolor=lightgray](8,1)(9,1)(9,2)(8,1)
\psline*[linecolor=lightgray](9,0)(10,0)(10,1)(9,0)
\psline*[linecolor=lightgray](11,0)(11,-0.5)(10.5,-0.5)(11,0)
\psline*[linecolor=lightgray](10,11)(10.5,11.5)(11,11.5)(11,11)(10,11)
\psline*[linecolor=lightgray](11,11)(11.5,11)(11.5,10.5)(11,10)(11,11)
\psline*[linecolor=lightgray](0,11)(-0.5,11)(-0.5,10.5)(0,11)
\psline*[linecolor=lightgray](1,6)(2,6)(2,7)(1,6)
\psline*[linecolor=lightgray](2,5)(3,5)(3,6)(2,5)
\psline*[linecolor=lightgray](3,4)(5,4)(4,3)(4,5)(3,4)
\psline*[linecolor=lightgray](5,2)(5,3)(6,3)(5,2)
\psline*[linecolor=lightgray](6,1)(6,2)(7,2)(6,1)
\psline*[linecolor=lightgray](7,0)(7,1)(8,1)(7,0)
\psline*[linecolor=lightgray](4,9)(4,10)(5,10)(4,9)
\psline*[linecolor=lightgray](10,3)(11,3)(11,4)(10,3)
 
\psline*[linecolor=darkgray](0,9)(1,9)(1,10)(0,9)
\psline*[linecolor=darkgray](1,8)(2,8)(2,9)(1,8)
\psline*[linecolor=darkgray](4,7)(4,8)(5,8)(4,7)
\psline*[linecolor=darkgray](5,6)(5,7)(6,7)(5,6)
\psline*[linecolor=darkgray](6,5)(7,5)(7,6)(6,5)
\psline*[linecolor=darkgray](7,4)(8,4)(8,5)(7,4)
\psline*[linecolor=darkgray](8,3)(9,3)(9,4)(8,3)
\psline*[linecolor=darkgray](9,0)(9,1)(10,1)(9,0)
\psline*[linecolor=darkgray](10,0)(10,-0.5)(10.5,-0.5)(11,0)(10,0)
\psline*[linecolor=darkgray](10,11)(10,11.5)(10.5,11.5)(10,11)
\psline*[linecolor=darkgray](0,10)(-0.5,10)(-0.5,10.5)(0,11)(0,10)
\psline*[linecolor=darkgray](11,10)(11.5,10)(11.5,10.5)(11,10)
\psline*[linecolor=darkgray](2,5)(2,6)(3,6)(2,5)
\psline*[linecolor=darkgray](3,4)(3,5)(4,5)(3,4)
\psline*[linecolor=darkgray](4,3)(5,3)(5,4)(4,3)
\psline*[linecolor=darkgray](5,2)(6,2)(6,3)(5,2)
\psline*[linecolor=darkgray](6,1)(7,1)(7,2)(6,1)
 
\psline[linewidth=2pt,linecolor=blue](-0.5,0)(11.5,0)
\psline[linewidth=2pt,linecolor=blue](-0.5,1)(11.5,1)
\psline[linewidth=2pt,linecolor=blue](-0.5,2)(11.5,2)
\psline[linewidth=2pt,linecolor=blue](-0.5,3)(11.5,3)
\psline[linewidth=2pt,linecolor=blue](-0.5,4)(11.5,4)
\psline[linewidth=2pt,linecolor=blue](-0.5,5)(11.5,5)
\psline[linewidth=2pt,linecolor=blue](-0.5,6)(11.5,6)
\psline[linewidth=2pt,linecolor=blue](-0.5,7)(11.5,7)
\psline[linewidth=2pt,linecolor=blue](-0.5,8)(11.5,8)
\psline[linewidth=2pt,linecolor=blue](-0.5,9)(11.5,9)
\psline[linewidth=2pt,linecolor=blue](-0.5,10)(11.5,10)
\psline[linewidth=2pt,linecolor=blue](-0.5,11)(11.5,11)
\psline[linewidth=2pt,linecolor=blue](0,-0.5)(0,11.5)
\psline[linewidth=2pt,linecolor=blue](1,-0.5)(1,11.5)
\psline[linewidth=2pt,linecolor=blue](2,-0.5)(2,11.5)
\psline[linewidth=2pt,linecolor=blue](3,-0.5)(3,11.5)
\psline[linewidth=2pt,linecolor=blue](4,-0.5)(4,11.5)
\psline[linewidth=2pt,linecolor=blue](5,-0.5)(5,11.5)
\psline[linewidth=2pt,linecolor=blue](6,-0.5)(6,11.5)
\psline[linewidth=2pt,linecolor=blue](7,-0.5)(7,11.5)
\psline[linewidth=2pt,linecolor=blue](8,-0.5)(8,11.5)
\psline[linewidth=2pt,linecolor=blue](9,-0.5)(9,11.5)
\psline[linewidth=2pt,linecolor=blue](10,-0.5)(10,11.5)
\psline[linewidth=2pt,linecolor=blue](11,-0.5)(11,11.5)
\psline[linewidth=2pt,linecolor=blue](-0.5,-0.5)(11.5,11.5)
\psline[linewidth=2pt,linecolor=blue](0.5,-0.5)(11.5,10.5)
\psline[linewidth=2pt,linecolor=blue](1.5,-0.5)(11.5,9.5)
\psline[linewidth=2pt,linecolor=blue](2.5,-0.5)(11.5,8.5)
\psline[linewidth=2pt,linecolor=blue](3.5,-0.5)(11.5,7.5)
\psline[linewidth=2pt,linecolor=blue](4.5,-0.5)(11.5,6.5)
\psline[linewidth=2pt,linecolor=blue](5.5,-0.5)(11.5,5.5)
\psline[linewidth=2pt,linecolor=blue](6.5,-0.5)(11.5,4.5)
\psline[linewidth=2pt,linecolor=blue](7.5,-0.5)(11.5,3.5)
\psline[linewidth=2pt,linecolor=blue](8.5,-0.5)(11.5,2.5)
\psline[linewidth=2pt,linecolor=blue](9.5,-0.5)(11.5,1.5)
\psline[linewidth=2pt,linecolor=blue](10.5,-0.5)(11.5,0.5)
\psline[linewidth=2pt,linecolor=blue](-0.5,0.5)(10.5,11.5)
\psline[linewidth=2pt,linecolor=blue](-0.5,1.5)(9.5,11.5)
\psline[linewidth=2pt,linecolor=blue](-0.5,2.5)(8.5,11.5)
\psline[linewidth=2pt,linecolor=blue](-0.5,3.5)(7.5,11.5)
\psline[linewidth=2pt,linecolor=blue](-0.5,4.5)(6.5,11.5)
\psline[linewidth=2pt,linecolor=blue](-0.5,5.5)(5.5,11.5)
\psline[linewidth=2pt,linecolor=blue](-0.5,6.5)(4.5,11.5)
\psline[linewidth=2pt,linecolor=blue](-0.5,7.5)(3.5,11.5)
\psline[linewidth=2pt,linecolor=blue](-0.5,8.5)(2.5,11.5)
\psline[linewidth=2pt,linecolor=blue](-0.5,9.5)(1.5,11.5)
\psline[linewidth=2pt,linecolor=blue](-0.5,10.5)(0.5,11.5)
\multirput{0}(0,0)(0,1){12}{%
  \multirput{0}(0,0)(1,0){12}{%
     \pscircle*[linecolor=white]{8pt}
     \pscircle[linewidth=1pt,linecolor=black]{8pt}
   }
}
\rput{0}(0,0){{\bf 4}}
\rput{0}(0,1){{\bf 1}}
\rput{0}(0,2){{\bf 3}}
\rput{0}(0,3){{\bf 4}}
\rput{0}(0,4){{\bf  }}
\rput{0}(0,5){{\bf  }}
\rput{0}(0,6){{\bf  }}
\rput{0}(0,7){{\bf 1}}
\rput{0}(0,8){{\bf 4}}
\rput{0}(0,9){{\bf 2}}
\rput{0}(0,10){{\bf 3}}
\rput{0}(0,11){{\bf 1}}

\rput{0}(1,0){{\bf 2}}
\rput{0}(1,1){{\bf 3}}
\rput{0}(1,2){{\bf 4}}
\rput{0}(1,3){{\bf  }}
\rput{0}(1,4){{\bf  }}
\rput{0}(1,5){{\bf  }}
\rput{0}(1,6){{\bf 1}}
\rput{0}(1,7){{\bf 4}}
\rput{0}(1,8){{\bf 2}}
\rput{0}(1,9){{\bf 3}}
\rput{0}(1,10){{\bf 1}}
\rput{0}(1,11){{\bf 4}}

\rput{0}(2,0){{\bf 3}}
\rput{0}(2,1){{\bf 4}}
\rput{0}(2,2){{\bf  }}
\rput{0}(2,3){{\bf  }}
\rput{0}(2,4){{\bf  }}
\rput{0}(2,5){{\bf 1}}
\rput{0}(2,6){{\bf 3}}
\rput{0}(2,7){{\bf 2}}
\rput{0}(2,8){{\bf 3}}
\rput{0}(2,9){{\bf 1}}
\rput{0}(2,10){{\bf 4}}
\rput{0}(2,11){{\bf 2}}

\rput{0}(3,0){{\bf 4}}
\rput{0}(3,1){{\bf  }}
\rput{0}(3,2){{\bf  }}
\rput{0}(3,3){{\bf  }}
\rput{0}(3,4){{\bf 1}}
\rput{0}(3,5){{\bf 3}}
\rput{0}(3,6){{\bf 2}}
\rput{0}(3,7){{\bf 4}}
\rput{0}(3,8){{\bf 1}}
\rput{0}(3,9){{\bf 4}}
\rput{0}(3,10){{\bf 2}}
\rput{0}(3,11){{\bf 3}}

\rput{0}(4,0){{\bf  }}
\rput{0}(4,1){{\bf  }}
\rput{0}(4,2){{\bf  }}
\rput{0}(4,3){{\bf 2}}
\rput{0}(4,4){{\bf 3}}
\rput{0}(4,5){{\bf 2}}
\rput{0}(4,6){{\bf 4}}
\rput{0}(4,7){{\bf 1}}
\rput{0}(4,8){{\bf 3}}
\rput{0}(4,9){{\bf 2}}
\rput{0}(4,10){{\bf 3}}
\rput{0}(4,11){{\bf 4}}

\rput{0}(5,0){{\bf  }}
\rput{0}(5,1){{\bf  }}
\rput{0}(5,2){{\bf 2}}
\rput{0}(5,3){{\bf 3}}
\rput{0}(5,4){{\bf 1}}
\rput{0}(5,5){{\bf 4}}
\rput{0}(5,6){{\bf 1}}
\rput{0}(5,7){{\bf 3}}
\rput{0}(5,8){{\bf 2}}
\rput{0}(5,9){{\bf 4}}
\rput{0}(5,10){{\bf 1}}
\rput{0}(5,11){{\bf  }}

\rput{0}(6,0){{\bf  }}
\rput{0}(6,1){{\bf 2}}
\rput{0}(6,2){{\bf 3}}
\rput{0}(6,3){{\bf 1}}
\rput{0}(6,4){{\bf 4}}
\rput{0}(6,5){{\bf 2}}
\rput{0}(6,6){{\bf 3}}
\rput{0}(6,7){{\bf 2}}
\rput{0}(6,8){{\bf 4}}
\rput{0}(6,9){{\bf 3}}
\rput{0}(6,10){{\bf  }}
\rput{0}(6,11){{\bf  }}

\rput{0}(7,0){{\bf 2}}
\rput{0}(7,1){{\bf 3}}
\rput{0}(7,2){{\bf 1}}
\rput{0}(7,3){{\bf 4}}
\rput{0}(7,4){{\bf 2}}
\rput{0}(7,5){{\bf 3}}
\rput{0}(7,6){{\bf 1}}
\rput{0}(7,7){{\bf 4}}
\rput{0}(7,8){{\bf 3}}
\rput{0}(7,9){{\bf  }}
\rput{0}(7,10){{\bf  }}
\rput{0}(7,11){{\bf  }}

\rput{0}(8,0){{\bf 4}}
\rput{0}(8,1){{\bf 1}}
\rput{0}(8,2){{\bf 4}}
\rput{0}(8,3){{\bf 2}}
\rput{0}(8,4){{\bf 3}}
\rput{0}(8,5){{\bf 1}}
\rput{0}(8,6){{\bf 4}}
\rput{0}(8,7){{\bf 3}}
\rput{0}(8,8){{\bf  }}
\rput{0}(8,9){{\bf  }}
\rput{0}(8,10){{\bf  }}
\rput{0}(8,11){{\bf 2}}

\rput{0}(9,0){{\bf 1}}
\rput{0}(9,1){{\bf 3}}
\rput{0}(9,2){{\bf 2}}
\rput{0}(9,3){{\bf 3}}
\rput{0}(9,4){{\bf 1}}
\rput{0}(9,5){{\bf 4}}
\rput{0}(9,6){{\bf 3}}
\rput{0}(9,7){{\bf  }}
\rput{0}(9,8){{\bf  }}
\rput{0}(9,9){{\bf  }}
\rput{0}(9,10){{\bf 2}}
\rput{0}(9,11){{\bf 4}}

\rput{0}(10,0){{\bf  3}}
\rput{0}(10,1){{\bf  2}}
\rput{0}(10,2){{\bf  4}}
\rput{0}(10,3){{\bf  1}}
\rput{0}(10,4){{\bf  4}}
\rput{0}(10,5){{\bf  3}}
\rput{0}(10,6){{\bf   }}
\rput{0}(10,7){{\bf   }}
\rput{0}(10,8){{\bf   }}
\rput{0}(10,9){{\bf  1}}
\rput{0}(10,10){{\bf 4}}
\rput{0}(10,11){{\bf 1}}

\rput{0}(11,0){{\bf 2}}
\rput{0}(11,1){{\bf 4}}
\rput{0}(11,2){{\bf 1}}
\rput{0}(11,3){{\bf 3}}
\rput{0}(11,4){{\bf 2}}
\rput{0}(11,5){{\bf  }}
\rput{0}(11,6){{\bf  }}
\rput{0}(11,7){{\bf  }}
\rput{0}(11,8){{\bf 1}}
\rput{0}(11,9){{\bf 4}}
\rput{0}(11,10){{\bf 2}}
\rput{0}(11,11){{\bf 3}}
%
%
\multirput{0}(11.5,-0.5)(0,1){12}{%
     \psline[linewidth=2pt,linecolor=black]{->}(0.2,-0.2)(-0.1,0.1) 
}
\uput[0](11.7,-0.7){D1} 
\uput[0](11.7, 0.3){D2} 
\uput[0](11.7, 1.3){D3} 
\uput[0](11.7, 2.3){D4}
\uput[0](11.7, 3.3){D5}
\uput[0](11.7, 4.3){D6}
\uput[0](11.7, 5.3){D7}
\uput[0](11.7, 6.3){D8}
\uput[0](11.7, 7.3){D9} 
\uput[0](11.7, 8.3){D10} 
\uput[0](11.7, 9.3){D11} 
\uput[0](11.7,10.3){D12} 
\endpspicture
\caption{ \label{prop.12k.fig2}
The 4-coloring of $T(12,12)$ after Step~3 in the case $L=4k$.}
\end{figure}
%
%

On D2, we color $3$ the $6k$ consecutive vertices with $x$-coordinates
$3k+2\leq x \leq 9k+1$. The other vertices on D2 are colored $4$. 
The vertices on D$(12k)$ are colored $3$ or $4$ in such a way that the 
resulting coloring is proper (for each vertex, the choice is unique).

We color all vertices on D3 and D$(12k-1)$ using colors $1$ and $2$. We then
color D4 and D$(12k-2)$ using colors $3$ and $4$. Again the condition
that $f$ is proper implies that for each vertex the choice is unique.  
The partial degree of $f$ is $\deg f|_R = 4$. 

\medskip
\noindent
{\bf Step 2.}
For $k>1$, we find that there are $12k-7$ counter-diagonals to be colored,
and we need to sequentially color all of them but five. This can be 
achieved by performing the following procedure: suppose that we
have already colored counter-diagonals D$j$ and D$(12k-j-2)$ ($j\geq 4$) 
using colors $3$ and $4$. Then, we color D$(j+1)$ and D$(12k-j+1)$ using 
$1$ and $2$, and then, we color D$(j+2)$ and D$(12k-j)$ using $3$ and $4$. 
Again, for each vertex we have only one choice. This step is repeated 
$3(k-1)$ times: we add $12(k-1)$ counter-diagonals, and there are only
five counter-diagonals not yet colored. Indeed, the last colored 
counter-diagonals use colors $3$ and $4$, as it was at the end of Step~1. 

Each of these $3(k-1)$ steps adds $4$ to the degree of the coloring.
Thus, the partial degree of the coloring is  $\deg f|_R = 4 + 12(k-1)$.  

\medskip
\noindent
{\bf Step 3.}
The last colored counter-diagonals are D$(6k-2)$ and D$(6k+4)$. 

On D$(6k-1)$, the vertices at $(6k,12k-1)$ and $(12k,6k-1)$ only admit one
color: one of them should have color $1$ and the other one $2$.
The rest of the vertices on D$(6k-1)$ are colored $3$ or $4$  
(again, there is a unique choice for each vertex). 

We color $1$ or $2$ all vertices on D$(6k+3)$; again there is a unique 
choice for each vertex.  
As shown in Figure~\ref{prop.12k.fig2}, the contribution to the degree of 
these new triangles is $4$; thus, the partial degree of $f$ is  
$\deg f|_R = 8 + 12(k-1)$.

%
%
\begin{figure}[htb]
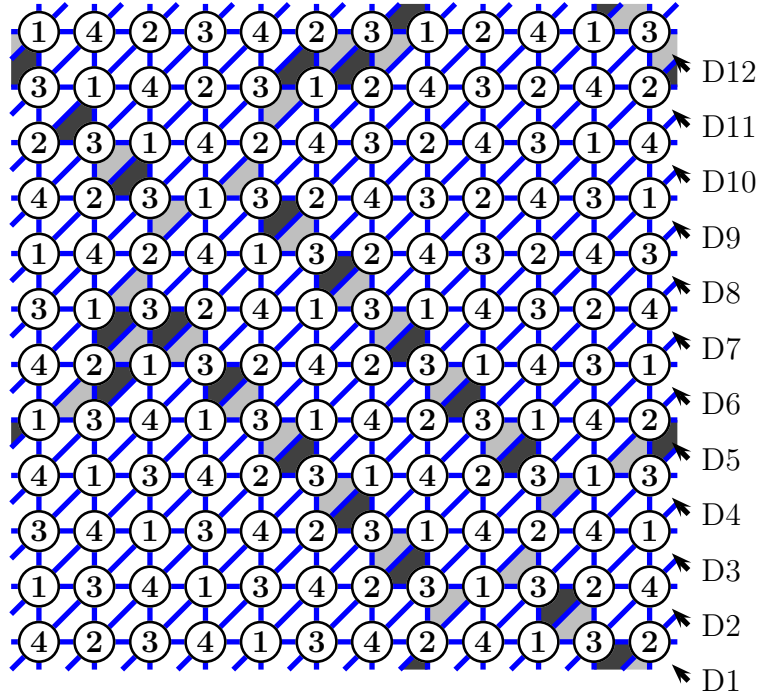

\centering
%
%
\psset{xunit=21pt}
\psset{yunit=21pt}
\psset{labelsep=5pt}
\pspicture(-1,-1)(13,12)
\psline*[linecolor=lightgray](1,8)(1,9)(2,9)(1,8)
\psline*[linecolor=lightgray](2,7)(2,8)(4,8)(4,9)(2,7)
\psline*[linecolor=lightgray](4,7)(5,7)(5,8)(4,7)
\psline*[linecolor=lightgray](5,6)(7,6)(6,5)(6,7)(5,6)
\psline*[linecolor=lightgray](7,5)(8,5)(7,4)(7,5)
\psline*[linecolor=lightgray](8,4)(9,4)(8,3)(8,4)
\psline*[linecolor=lightgray](9,3)(10,3)(9,2)(9,3)
\psline*[linecolor=lightgray](8,1)(9,1)(9,2)(8,1)
\psline*[linecolor=lightgray](9,0)(10,0)(10,1)(9,0)
\psline*[linecolor=lightgray](11,0)(11,-0.5)(10.5,-0.5)(11,0)
\psline*[linecolor=lightgray](10,11)(10.5,11.5)(11,11.5)(11,11)(10,11)
\psline*[linecolor=lightgray](11,11)(11.5,11)(11.5,10.5)(11,10)(11,11)
\psline*[linecolor=lightgray](0,11)(-0.5,11)(-0.5,10.5)(0,11)
\psline*[linecolor=lightgray](1,6)(2,6)(2,7)(1,6)
\psline*[linecolor=lightgray](2,5)(3,5)(3,6)(2,5)
\psline*[linecolor=lightgray](3,4)(5,4)(4,3)(4,5)(3,4)
\psline*[linecolor=lightgray](5,2)(5,3)(6,3)(5,2)
\psline*[linecolor=lightgray](6,1)(6,2)(7,2)(6,1)
\psline*[linecolor=lightgray](7,0)(7,1)(8,1)(7,0)
\psline*[linecolor=lightgray](4,9)(4,10)(5,10)(4,9)
\psline*[linecolor=lightgray](10,3)(11,3)(11,4)(10,3)
\psline*[linecolor=lightgray](5,10)(5,11)(6,11)(5,10)
\psline*[linecolor=lightgray](1,5)(2,5)(2,6)(1,5)
\psline*[linecolor=lightgray](0,4)(1,4)(1,5)(0,4)
\psline*[linecolor=lightgray](6,10)(6,11)(7,11)(6,10)
 
\psline*[linecolor=darkgray](0,9)(1,9)(1,10)(0,9)
\psline*[linecolor=darkgray](1,8)(2,8)(2,9)(1,8)
\psline*[linecolor=darkgray](4,7)(4,8)(5,8)(4,7)
\psline*[linecolor=darkgray](5,6)(5,7)(6,7)(5,6)
\psline*[linecolor=darkgray](6,5)(7,5)(7,6)(6,5)
\psline*[linecolor=darkgray](7,4)(8,4)(8,5)(7,4)
\psline*[linecolor=darkgray](8,3)(9,3)(9,4)(8,3)
\psline*[linecolor=darkgray](9,0)(9,1)(10,1)(9,0)
\psline*[linecolor=darkgray](10,0)(10,-0.5)(10.5,-0.5)(11,0)(10,0)
\psline*[linecolor=darkgray](10,11)(10,11.5)(10.5,11.5)(10,11)
\psline*[linecolor=darkgray](0,10)(-0.5,10)(-0.5,10.5)(0,11)(0,10)
\psline*[linecolor=darkgray](11,10)(11.5,10)(11.5,10.5)(11,10)
\psline*[linecolor=darkgray](2,5)(2,6)(3,6)(2,5)
\psline*[linecolor=darkgray](3,4)(3,5)(4,5)(3,4)
\psline*[linecolor=darkgray](4,3)(5,3)(5,4)(4,3)
\psline*[linecolor=darkgray](5,2)(6,2)(6,3)(5,2)
\psline*[linecolor=darkgray](6,1)(7,1)(7,2)(6,1)
\psline*[linecolor=darkgray](5,10)(6,10)(6,11)(5,10)
\psline*[linecolor=darkgray](4,10)(5,10)(5,11)(4,10)
\psline*[linecolor=darkgray](1,5)(1,6)(2,6)(1,5)
\psline*[linecolor=darkgray](11,4)(11.5,4)(11.5,3.5)(11,3)
\psline*[linecolor=darkgray](0,4)(-0.5,4)(-0.5,3.5)(0,4)
\psline*[linecolor=darkgray](1,4)(1,5)(2,5)(1,4)
\psline*[linecolor=darkgray](6,11)(6.5,11.5)(7,11.5)(7,11)(6,11)
\psline*[linecolor=darkgray](7,0)(7,-0.5)(6.5,-0.5)(7,0)
 
\psline[linewidth=2pt,linecolor=blue](-0.5,0)(11.5,0)
\psline[linewidth=2pt,linecolor=blue](-0.5,1)(11.5,1)
\psline[linewidth=2pt,linecolor=blue](-0.5,2)(11.5,2)
\psline[linewidth=2pt,linecolor=blue](-0.5,3)(11.5,3)
\psline[linewidth=2pt,linecolor=blue](-0.5,4)(11.5,4)
\psline[linewidth=2pt,linecolor=blue](-0.5,5)(11.5,5)
\psline[linewidth=2pt,linecolor=blue](-0.5,6)(11.5,6)
\psline[linewidth=2pt,linecolor=blue](-0.5,7)(11.5,7)
\psline[linewidth=2pt,linecolor=blue](-0.5,8)(11.5,8)
\psline[linewidth=2pt,linecolor=blue](-0.5,9)(11.5,9)
\psline[linewidth=2pt,linecolor=blue](-0.5,10)(11.5,10)
\psline[linewidth=2pt,linecolor=blue](-0.5,11)(11.5,11)
\psline[linewidth=2pt,linecolor=blue](0,-0.5)(0,11.5)
\psline[linewidth=2pt,linecolor=blue](1,-0.5)(1,11.5)
\psline[linewidth=2pt,linecolor=blue](2,-0.5)(2,11.5)
\psline[linewidth=2pt,linecolor=blue](3,-0.5)(3,11.5)
\psline[linewidth=2pt,linecolor=blue](4,-0.5)(4,11.5)
\psline[linewidth=2pt,linecolor=blue](5,-0.5)(5,11.5)
\psline[linewidth=2pt,linecolor=blue](6,-0.5)(6,11.5)
\psline[linewidth=2pt,linecolor=blue](7,-0.5)(7,11.5)
\psline[linewidth=2pt,linecolor=blue](8,-0.5)(8,11.5)
\psline[linewidth=2pt,linecolor=blue](9,-0.5)(9,11.5)
\psline[linewidth=2pt,linecolor=blue](10,-0.5)(10,11.5)
\psline[linewidth=2pt,linecolor=blue](11,-0.5)(11,11.5)
\psline[linewidth=2pt,linecolor=blue](-0.5,-0.5)(11.5,11.5)
\psline[linewidth=2pt,linecolor=blue](0.5,-0.5)(11.5,10.5)
\psline[linewidth=2pt,linecolor=blue](1.5,-0.5)(11.5,9.5)
\psline[linewidth=2pt,linecolor=blue](2.5,-0.5)(11.5,8.5)
\psline[linewidth=2pt,linecolor=blue](3.5,-0.5)(11.5,7.5)
\psline[linewidth=2pt,linecolor=blue](4.5,-0.5)(11.5,6.5)
\psline[linewidth=2pt,linecolor=blue](5.5,-0.5)(11.5,5.5)
\psline[linewidth=2pt,linecolor=blue](6.5,-0.5)(11.5,4.5)
\psline[linewidth=2pt,linecolor=blue](7.5,-0.5)(11.5,3.5)
\psline[linewidth=2pt,linecolor=blue](8.5,-0.5)(11.5,2.5)
\psline[linewidth=2pt,linecolor=blue](9.5,-0.5)(11.5,1.5)
\psline[linewidth=2pt,linecolor=blue](10.5,-0.5)(11.5,0.5)
\psline[linewidth=2pt,linecolor=blue](-0.5,0.5)(10.5,11.5)
\psline[linewidth=2pt,linecolor=blue](-0.5,1.5)(9.5,11.5)
\psline[linewidth=2pt,linecolor=blue](-0.5,2.5)(8.5,11.5)
\psline[linewidth=2pt,linecolor=blue](-0.5,3.5)(7.5,11.5)
\psline[linewidth=2pt,linecolor=blue](-0.5,4.5)(6.5,11.5)
\psline[linewidth=2pt,linecolor=blue](-0.5,5.5)(5.5,11.5)
\psline[linewidth=2pt,linecolor=blue](-0.5,6.5)(4.5,11.5)
\psline[linewidth=2pt,linecolor=blue](-0.5,7.5)(3.5,11.5)
\psline[linewidth=2pt,linecolor=blue](-0.5,8.5)(2.5,11.5)
\psline[linewidth=2pt,linecolor=blue](-0.5,9.5)(1.5,11.5)
\psline[linewidth=2pt,linecolor=blue](-0.5,10.5)(0.5,11.5)
\multirput{0}(0,0)(0,1){12}{%
  \multirput{0}(0,0)(1,0){12}{%
     \pscircle*[linecolor=white]{8pt}
     \pscircle[linewidth=1pt,linecolor=black] {8pt}
   }
}
\rput{0}(0,0){{\bf 4}}
\rput{0}(0,1){{\bf 1}}
\rput{0}(0,2){{\bf 3}}
\rput{0}(0,3){{\bf 4}}
\rput{0}(0,4){{\bf 1}}
\rput{0}(0,5){{\bf 4}}
\rput{0}(0,6){{\bf 3}}
\rput{0}(0,7){{\bf 1}}
\rput{0}(0,8){{\bf 4}}
\rput{0}(0,9){{\bf 2}}
\rput{0}(0,10){{\bf 3}}
\rput{0}(0,11){{\bf 1}}

\rput{0}(1,0){{\bf 2}}
\rput{0}(1,1){{\bf 3}}
\rput{0}(1,2){{\bf 4}}
\rput{0}(1,3){{\bf 1}}
\rput{0}(1,4){{\bf 3}}
\rput{0}(1,5){{\bf 2}}
\rput{0}(1,6){{\bf 1}}
\rput{0}(1,7){{\bf 4}}
\rput{0}(1,8){{\bf 2}}
\rput{0}(1,9){{\bf 3}}
\rput{0}(1,10){{\bf 1}}
\rput{0}(1,11){{\bf 4}}

\rput{0}(2,0){{\bf 3}}
\rput{0}(2,1){{\bf 4}}
\rput{0}(2,2){{\bf 1}}
\rput{0}(2,3){{\bf 3}}
\rput{0}(2,4){{\bf 4}}
\rput{0}(2,5){{\bf 1}}
\rput{0}(2,6){{\bf 3}}
\rput{0}(2,7){{\bf 2}}
\rput{0}(2,8){{\bf 3}}
\rput{0}(2,9){{\bf 1}}
\rput{0}(2,10){{\bf 4}}
\rput{0}(2,11){{\bf 2}}

\rput{0}(3,0){{\bf 4}}
\rput{0}(3,1){{\bf 1}}
\rput{0}(3,2){{\bf 3}}
\rput{0}(3,3){{\bf 4}}
\rput{0}(3,4){{\bf 1}}
\rput{0}(3,5){{\bf 3}}
\rput{0}(3,6){{\bf 2}}
\rput{0}(3,7){{\bf 4}}
\rput{0}(3,8){{\bf 1}}
\rput{0}(3,9){{\bf 4}}
\rput{0}(3,10){{\bf 2}}
\rput{0}(3,11){{\bf 3}}

\rput{0}(4,0){{\bf 1}}
\rput{0}(4,1){{\bf 3}}
\rput{0}(4,2){{\bf 4}}
\rput{0}(4,3){{\bf 2}}
\rput{0}(4,4){{\bf 3}}
\rput{0}(4,5){{\bf 2}}
\rput{0}(4,6){{\bf 4}}
\rput{0}(4,7){{\bf 1}}
\rput{0}(4,8){{\bf 3}}
\rput{0}(4,9){{\bf 2}}
\rput{0}(4,10){{\bf 3}}
\rput{0}(4,11){{\bf 4}}

\rput{0}(5,0){{\bf 3}}
\rput{0}(5,1){{\bf 4}}
\rput{0}(5,2){{\bf 2}}
\rput{0}(5,3){{\bf 3}}
\rput{0}(5,4){{\bf 1}}
\rput{0}(5,5){{\bf 4}}
\rput{0}(5,6){{\bf 1}}
\rput{0}(5,7){{\bf 3}}
\rput{0}(5,8){{\bf 2}}
\rput{0}(5,9){{\bf 4}}
\rput{0}(5,10){{\bf 1}}
\rput{0}(5,11){{\bf 2}}

\rput{0}(6,0){{\bf 4}}
\rput{0}(6,1){{\bf 2}}
\rput{0}(6,2){{\bf 3}}
\rput{0}(6,3){{\bf 1}}
\rput{0}(6,4){{\bf 4}}
\rput{0}(6,5){{\bf 2}}
\rput{0}(6,6){{\bf 3}}
\rput{0}(6,7){{\bf 2}}
\rput{0}(6,8){{\bf 4}}
\rput{0}(6,9){{\bf 3}}
\rput{0}(6,10){{\bf 2}}
\rput{0}(6,11){{\bf 3}}

\rput{0}(7,0){{\bf 2}}
\rput{0}(7,1){{\bf 3}}
\rput{0}(7,2){{\bf 1}}
\rput{0}(7,3){{\bf 4}}
\rput{0}(7,4){{\bf 2}}
\rput{0}(7,5){{\bf 3}}
\rput{0}(7,6){{\bf 1}}
\rput{0}(7,7){{\bf 4}}
\rput{0}(7,8){{\bf 3}}
\rput{0}(7,9){{\bf 2}}
\rput{0}(7,10){{\bf 4}}
\rput{0}(7,11){{\bf 1}}

\rput{0}(8,0){{\bf 4}}
\rput{0}(8,1){{\bf 1}}
\rput{0}(8,2){{\bf 4}}
\rput{0}(8,3){{\bf 2}}
\rput{0}(8,4){{\bf 3}}
\rput{0}(8,5){{\bf 1}}
\rput{0}(8,6){{\bf 4}}
\rput{0}(8,7){{\bf 3}}
\rput{0}(8,8){{\bf 2}}
\rput{0}(8,9){{\bf 4}}
\rput{0}(8,10){{\bf 3}}
\rput{0}(8,11){{\bf 2}}

\rput{0}(9,0){{\bf 1}}
\rput{0}(9,1){{\bf 3}}
\rput{0}(9,2){{\bf 2}}
\rput{0}(9,3){{\bf 3}}
\rput{0}(9,4){{\bf 1}}
\rput{0}(9,5){{\bf 4}}
\rput{0}(9,6){{\bf 3}}
\rput{0}(9,7){{\bf 2}}
\rput{0}(9,8){{\bf 4}}
\rput{0}(9,9){{\bf 3}}
\rput{0}(9,10){{\bf 2}}
\rput{0}(9,11){{\bf 4}}

\rput{0}(10,0){{\bf  3}}
\rput{0}(10,1){{\bf  2}}
\rput{0}(10,2){{\bf  4}}
\rput{0}(10,3){{\bf  1}}
\rput{0}(10,4){{\bf  4}}
\rput{0}(10,5){{\bf  3}}
\rput{0}(10,6){{\bf  2}}
\rput{0}(10,7){{\bf  4}}
\rput{0}(10,8){{\bf  3}}
\rput{0}(10,9){{\bf  1}}
\rput{0}(10,10){{\bf 4}}
\rput{0}(10,11){{\bf 1}}

\rput{0}(11,0){{\bf 2}}
\rput{0}(11,1){{\bf 4}}
\rput{0}(11,2){{\bf 1}}
\rput{0}(11,3){{\bf 3}}
\rput{0}(11,4){{\bf 2}}
\rput{0}(11,5){{\bf 1}}
\rput{0}(11,6){{\bf 4}}
\rput{0}(11,7){{\bf 3}}
\rput{0}(11,8){{\bf 1}}
\rput{0}(11,9){{\bf 4}}
\rput{0}(11,10){{\bf 2}}
\rput{0}(11,11){{\bf 3}}
%
%
\multirput{0}(11.5,-0.5)(0,1){12}{%
     \psline[linewidth=2pt,linecolor=black]{->}(0.2,-0.2)(-0.1,0.1) 
}
\uput[0](11.7,-0.7){D1} 
\uput[0](11.7, 0.3){D2} 
\uput[0](11.7, 1.3){D3} 
\uput[0](11.7, 2.3){D4}
\uput[0](11.7, 3.3){D5}
\uput[0](11.7, 4.3){D6}
\uput[0](11.7, 5.3){D7}
\uput[0](11.7, 6.3){D8}
\uput[0](11.7, 7.3){D9} 
\uput[0](11.7, 8.3){D10} 
\uput[0](11.7, 9.3){D11} 
\uput[0](11.7,10.3){D12} 
\endpspicture
\caption{ \label{prop.12k.fig4}
The 4-coloring of $T(12,12)$ after Step~5 in the case $L=4k$.
}
\end{figure}
%
%

\noindent
{\bf Step 4.}
On D$(6k)$ the vertices at $(1,6k-1)$, $(12k,6k)$, $(6k+1,12k-1)$, and
$(6k,12k)$ only admit a unique color choice: either $1$ or $2$. The first
two vertices should be colored alike, while the last two vertices take the
other color. We color the other vertices on D$(6k)$ 
with 1 and 2 in such a way that those
vertices with $x$-coordinate satisfying $1\leq x < 6k$ take the same
color as the vertex at $(1,6k-1)$; the rest are colored the same as
the vertex at $(6k,12k)$.

All vertices on D$(6k+1)$ are colored $3$ or $4$. For all of them, except for
those at $(1,6k)$ and $(6k+1,12k)$, there is unique possibility to do so.
We color $4$ the vertex at $(1,6k)$, and color $3$ the vertex at $(6k+1,12k)$. 
The increment of the partial degree is $-2$, thus $\deg f|_R = 6 + 12(k-1)$. 

\medskip
\noindent
{\bf Step 5.}
Finally, on D$(6k+2)$, there are two vertices which only admit a single color 
chosen among $1$ and $2$. For odd $k$ these vertices are $(2,6k)$ and 
$(6k+2,12k)$; while for even $k$, these vertices are $(1,6k+1)$ and $(6k+1,1)$.
The other vertices on D$(6k+2)$ can be colored $3$ and $4$ (uniquely). 
The resulting coloring is depicted in Figure~\ref{prop.12k.fig4}. 
In this step, the increment in the degree is zero. Therefore, the degree 
of the obtained four-coloring is 
\begin{equation}
\deg f \;=\; 6 + 12(k-1) \;\equiv\; 6 \pmod{12} \nonumber
\end{equation}
This coloring $f$ of $T(12k,12k)$ is proper and its degree is congruent 
to six modulo $12$, as claimed.

\medskip

\proofofcase{4}{$L=4k+1$}

Let us consider the triangulation $T=T(12k+3,12k+3)$ with $k\in\N$ (we will
illustrate the main steps with the case $k=1$).

\medskip
\noindent
{\bf Step 1.} 
On D1 we color $1$ the $6k+2$ consecutive vertices with 
$x$-coordinate $1\leq x \leq 6k+2$. The other $6k+1$ vertices on D1 
are colored $2$.

On D2 we color $3$ the $6k+1$ consecutive vertices with               
$x$-coordinate $3k+3\leq x \leq 9k+3$. The other vertices on D2 are  
colored $4$. We color $3$ or $4$ all vertices on D$(12k+3)$; the choice
is unique for each vertex.  

We color $1$ or $2$ all vertices on D3, D5, D$(12k+2)$, and D$(12k)$. 
And we color $3$ or $4$ all vertices on D4 and D$(12k+1)$. In all cases,
the choice is unique for each vertex. 

The resulting (partial) coloring is depicted in Figure~\ref{prop.12k+3.fig1}. 
The partial degree of this coloring is $\deg f|_R = 8$. 

\medskip
\noindent
{\bf Step 2.}
For $k>1$, we find that there are $12k-6$ counter-diagonals to be colored
and in this step we will sequentially color all of them but six. This can be 
achieved by performing the following procedure: suppose that we
have already colored D$j$ and D$(12k-j+5)$ ($j\geq 5$) 
using colors $1$ and $2$. Then, we color D$(j+1)$ and D$(12k-j+4)$ using
colors $3$ and $4$, and D$(j+2)$ and D$(12k-j+3)$ using colors $1$ and $2$.
Again, for each vertex the choice is unique.
 
This step is repeated $3(k-1)$ times; thus, we add $12(k-1)$ counter-diagonals,
and there are only six counter-diagonals not yet colored. Indeed, the last 
colored counter-diagonals use colors $1$ and $2$, as it was at the end of 
Step~1. 

Each of these $3(k-1)$ steps adds $4$ to the degree of the coloring.
Thus, the partial degree is $\deg f|_R = 8 + 12(k-1)$.  

%
%
\begin{figure}[htb]
\centering
%
%
\psset{xunit=21pt}
\psset{yunit=21pt}
\psset{labelsep=5pt}
\pspicture(-1,-1)(16,15)
\psline*[linecolor=lightgray](0,2)(0,3)(-0.5,2.5)(-0.5,2)(0,2)
\psline*[linecolor=lightgray](0,1)(2,1)(1,0)(1,2)(0,1)
\psline*[linecolor=lightgray](2,0)(3,0)(2.5,-0.5)(2,-0.5)(2,0)
\psline*[linecolor=lightgray](0,14)(-0.5,14)(-0.5,13.5)(0,14)
\psline*[linecolor=lightgray](0,13)(1,13)(0,12)(0,13)
\psline*[linecolor=lightgray](1,12)(2,12)(1,11)(1,12)
\psline*[linecolor=lightgray](2,11)(3,11)(2,10)(2,11)
\psline*[linecolor=lightgray](3,10)(4,10)(3,9)(3,10)
\psline*[linecolor=lightgray](3,9)(3,8)(2,8)(3,9)
\psline*[linecolor=lightgray](4,8)(4,7)(3,7)(4,8)
\psline*[linecolor=lightgray](5,7)(5,6)(4,6)(5,7)
\psline*[linecolor=lightgray](5,5)(7,5)(6,4)(6,6)(5,5)
\psline*[linecolor=lightgray](7,3)(7,4)(8,4)(7,3)
\psline*[linecolor=lightgray](8,2)(8,3)(9,3)(8,2)
\psline*[linecolor=lightgray](9,1)(9,2)(10,2)(9,1)
\psline*[linecolor=lightgray](10,2)(11,2)(11,4)(12,4)(10,2)
\psline*[linecolor=lightgray](2,14)(2,14.5)(2.5,14.5)(2,14)
\psline*[linecolor=lightgray](3,14)(3,13)(4,14)(3,14)
\psline*[linecolor=lightgray](4,13)(4,12)(5,13)(4,13)
\psline*[linecolor=lightgray](5,12)(5,11)(6,12)(5,12)
\psline*[linecolor=lightgray](4,10)(5,10)(5,11)(4,10)
\psline*[linecolor=lightgray](5,9)(6,9)(6,10)(5,9)
\psline*[linecolor=lightgray](6,8)(7,8)(7,9)(6,8)
\psline*[linecolor=lightgray](7,7)(9,7)(8,6)(8,8)(7,7)
\psline*[linecolor=lightgray](9,5)(9,6)(10,6)(9,5)
\psline*[linecolor=lightgray](10,4)(10,5)(11,5)(10,4)
\psline*[linecolor=lightgray](12,4)(13,4)(13,5)(12,4)
\psline*[linecolor=lightgray](13,3)(14,3)(14,4)(13,3)
\psline*[linecolor=lightgray](14,2)(14.5,2)(14.5,2.5)(14,2)
\psline*[linecolor=lightgray](13,14)(13.5,14.5)(14,14.5)(14,14)(13,14)
\psline*[linecolor=lightgray](14,13)(14,14)(14.5,14)(14.5,13.5)(14,13)
\psline*[linecolor=lightgray](11,1)(12,1)(12,2)(11,1)
\psline*[linecolor=lightgray](12,0)(13,0)(13,1)(12,0)
\psline*[linecolor=lightgray](14,0)(14,-0.5)(13.5,-0.5)(14,0)

\psline*[linecolor=darkgray](0,3)(-0.5,3)(-0.5,2.5)(0,3)
\psline*[linecolor=darkgray](0,2)(1,2)(0,1)(0,2)
\psline*[linecolor=darkgray](1,0)(2,0)(2,1)(1,0)
\psline*[linecolor=darkgray](3,0)(3,-0.5)(2.5,-0.5)(3,0)
\psline*[linecolor=darkgray](0,13)(-0.5,13)(-0.5,13.5)(0,14)(0,13)
\psline*[linecolor=darkgray](0,12)(1,12)(1,13)(0,12)
\psline*[linecolor=darkgray](1,11)(2,11)(2,12)(1,11)
\psline*[linecolor=darkgray](2,10)(3,10)(3,11)(2,10)
\psline*[linecolor=darkgray](3,7)(3,8)(4,8)(3,7)
\psline*[linecolor=darkgray](4,6)(4,7)(5,7)(4,6)
\psline*[linecolor=darkgray](5,5)(5,6)(6,6)(5,5)
\psline*[linecolor=darkgray](6,4)(7,4)(7,5)(6,4)
\psline*[linecolor=darkgray](7,3)(8,3)(8,4)(7,3)
\psline*[linecolor=darkgray](8,2)(9,2)(9,3)(8,2)
\psline*[linecolor=darkgray](2,14)(2.5,14.5)(3,14.5)(3,14)(2,14)
\psline*[linecolor=darkgray](3,13)(4,13)(4,14)(3,13)
\psline*[linecolor=darkgray](4,12)(5,12)(5,13)(4,12)
\psline*[linecolor=darkgray](5,9)(5,10)(6,10)(5,9)
\psline*[linecolor=darkgray](6,8)(6,9)(7,9)(6,8)
\psline*[linecolor=darkgray](7,7)(7,8)(8,8)(7,7)
\psline*[linecolor=darkgray](8,6)(9,6)(9,7)(8,6)
\psline*[linecolor=darkgray](9,5)(10,5)(10,6)(9,5)
\psline*[linecolor=darkgray](10,4)(11,4)(11,5)(10,4)
\psline*[linecolor=darkgray](13,3)(13,4)(14,4)(13,3)
\psline*[linecolor=darkgray](14,2)(14,3)(14.5,3)(14.5,2.5)(14,2)
\psline*[linecolor=darkgray](13,14)(13,14.5)(13.5,14.5)(13,14)
\psline*[linecolor=darkgray](14,13)(14.5,13)(14.5,13.5)(14,13)
\psline*[linecolor=darkgray](11,1)(11,2)(12,2)(11,1)
\psline*[linecolor=darkgray](12,0)(12,1)(13,1)(12,0)
\psline*[linecolor=darkgray](13,0)(14,0)(13.5,-0.5)(13,-0.5)(13,0)

\psline[linewidth=2pt,linecolor=blue](-0.5,0)(14.5,0)
\psline[linewidth=2pt,linecolor=blue](-0.5,1)(14.5,1)
\psline[linewidth=2pt,linecolor=blue](-0.5,2)(14.5,2)
\psline[linewidth=2pt,linecolor=blue](-0.5,3)(14.5,3)
\psline[linewidth=2pt,linecolor=blue](-0.5,4)(14.5,4)
\psline[linewidth=2pt,linecolor=blue](-0.5,5)(14.5,5)
\psline[linewidth=2pt,linecolor=blue](-0.5,6)(14.5,6)
\psline[linewidth=2pt,linecolor=blue](-0.5,7)(14.5,7)
\psline[linewidth=2pt,linecolor=blue](-0.5,8)(14.5,8)
\psline[linewidth=2pt,linecolor=blue](-0.5,9)(14.5,9)
\psline[linewidth=2pt,linecolor=blue](-0.5,10)(14.5,10)
\psline[linewidth=2pt,linecolor=blue](-0.5,11)(14.5,11)
\psline[linewidth=2pt,linecolor=blue](-0.5,12)(14.5,12)
\psline[linewidth=2pt,linecolor=blue](-0.5,13)(14.5,13)
\psline[linewidth=2pt,linecolor=blue](-0.5,14)(14.5,14)
\psline[linewidth=2pt,linecolor=blue](0,-0.5)(0,14.5)
\psline[linewidth=2pt,linecolor=blue](1,-0.5)(1,14.5)
\psline[linewidth=2pt,linecolor=blue](2,-0.5)(2,14.5)
\psline[linewidth=2pt,linecolor=blue](3,-0.5)(3,14.5)
\psline[linewidth=2pt,linecolor=blue](4,-0.5)(4,14.5)
\psline[linewidth=2pt,linecolor=blue](5,-0.5)(5,14.5)
\psline[linewidth=2pt,linecolor=blue](6,-0.5)(6,14.5)
\psline[linewidth=2pt,linecolor=blue](7,-0.5)(7,14.5)
\psline[linewidth=2pt,linecolor=blue](8,-0.5)(8,14.5)
\psline[linewidth=2pt,linecolor=blue](9,-0.5)(9,14.5)
\psline[linewidth=2pt,linecolor=blue](10,-0.5)(10,14.5)
\psline[linewidth=2pt,linecolor=blue](11,-0.5)(11,14.5)
\psline[linewidth=2pt,linecolor=blue](12,-0.5)(12,14.5)
\psline[linewidth=2pt,linecolor=blue](13,-0.5)(13,14.5)
\psline[linewidth=2pt,linecolor=blue](14,-0.5)(14,14.5)
\psline[linewidth=2pt,linecolor=blue](-0.5,-0.5)(14.5,14.5)
\psline[linewidth=2pt,linecolor=blue](0.5,-0.5)(14.5,13.5)
\psline[linewidth=2pt,linecolor=blue](1.5,-0.5)(14.5,12.5)
\psline[linewidth=2pt,linecolor=blue](2.5,-0.5)(14.5,11.5)
\psline[linewidth=2pt,linecolor=blue](3.5,-0.5)(14.5,10.5)
\psline[linewidth=2pt,linecolor=blue](4.5,-0.5)(14.5,9.5)
\psline[linewidth=2pt,linecolor=blue](5.5,-0.5)(14.5,8.5)
\psline[linewidth=2pt,linecolor=blue](6.5,-0.5)(14.5,7.5)
\psline[linewidth=2pt,linecolor=blue](7.5,-0.5)(14.5,6.5)
\psline[linewidth=2pt,linecolor=blue](8.5,-0.5)(14.5,5.5)
\psline[linewidth=2pt,linecolor=blue](9.5,-0.5)(14.5,4.5)
\psline[linewidth=2pt,linecolor=blue](10.5,-0.5)(14.5,3.5)
\psline[linewidth=2pt,linecolor=blue](11.5,-0.5)(14.5,2.5)
\psline[linewidth=2pt,linecolor=blue](12.5,-0.5)(14.5,1.5)
\psline[linewidth=2pt,linecolor=blue](13.5,-0.5)(14.5,0.5)
\psline[linewidth=2pt,linecolor=blue](-0.5,0.5)(13.5,14.5)
\psline[linewidth=2pt,linecolor=blue](-0.5,1.5)(12.5,14.5)
\psline[linewidth=2pt,linecolor=blue](-0.5,2.5)(11.5,14.5)
\psline[linewidth=2pt,linecolor=blue](-0.5,3.5)(10.5,14.5)
\psline[linewidth=2pt,linecolor=blue](-0.5,4.5)(9.5,14.5)
\psline[linewidth=2pt,linecolor=blue](-0.5,5.5)(8.5,14.5)
\psline[linewidth=2pt,linecolor=blue](-0.5,6.5)(7.5,14.5)
\psline[linewidth=2pt,linecolor=blue](-0.5,7.5)(6.5,14.5)
\psline[linewidth=2pt,linecolor=blue](-0.5,8.5)(5.5,14.5)
\psline[linewidth=2pt,linecolor=blue](-0.5,9.5)(4.5,14.5)
\psline[linewidth=2pt,linecolor=blue](-0.5,10.5)(3.5,14.5)
\psline[linewidth=2pt,linecolor=blue](-0.5,11.5)(2.5,14.5)
\psline[linewidth=2pt,linecolor=blue](-0.5,12.5)(1.5,14.5)
\psline[linewidth=2pt,linecolor=blue](-0.5,13.5)(0.5,14.5)
\multirput{0}(0,0)(0,1){15}{%
  \multirput{0}(0,0)(1,0){15}{%
     \pscircle*[linecolor=white]{8pt}
     \pscircle[linewidth=1pt,linecolor=black] {8pt}
   }
}
\rput{0}(0,0)  {{\bf 4}}
\rput{0}(0,1)  {{\bf 1}}
\rput{0}(0,2)  {{\bf 3}}
\rput{0}(0,3)  {{\bf 2}}
\rput{0}(0,4)  {{\bf  }}
\rput{0}(0,5)  {{\bf  }}
\rput{0}(0,6)  {{\bf  }}
\rput{0}(0,7)  {{\bf  }}
\rput{0}(0,8)  {{\bf  }}
\rput{0}(0,9)  {{\bf  }}
\rput{0}(0,10) {{\bf 1}}
\rput{0}(0,11) {{\bf 4}}
\rput{0}(0,12) {{\bf 2}}
\rput{0}(0,13) {{\bf 3}}
\rput{0}(0,14) {{\bf 1}}

\rput{0}(1,0)  {{\bf 2}}
\rput{0}(1,1)  {{\bf 3}}
\rput{0}(1,2)  {{\bf 2}}
\rput{0}(1,3)  {{\bf  }}
\rput{0}(1,4)  {{\bf  }}
\rput{0}(1,5)  {{\bf  }}
\rput{0}(1,6)  {{\bf  }}
\rput{0}(1,7)  {{\bf  }}
\rput{0}(1,8)  {{\bf  }}
\rput{0}(1,9)  {{\bf 1}}
\rput{0}(1,10) {{\bf 4}}
\rput{0}(1,11) {{\bf 2}}
\rput{0}(1,12) {{\bf 3}}
\rput{0}(1,13) {{\bf 1}}
\rput{0}(1,14) {{\bf 4}}

\rput{0}(2,0)  {{\bf 3}}
\rput{0}(2,1)  {{\bf 1}}
\rput{0}(2,2)  {{\bf  }}
\rput{0}(2,3)  {{\bf  }}
\rput{0}(2,4)  {{\bf  }}
\rput{0}(2,5)  {{\bf  }}
\rput{0}(2,6)  {{\bf  }}
\rput{0}(2,7)  {{\bf  }}
\rput{0}(2,8)  {{\bf 1}}
\rput{0}(2,9)  {{\bf 4}}
\rput{0}(2,10) {{\bf 2}}
\rput{0}(2,11) {{\bf 3}}
\rput{0}(2,12) {{\bf 1}}
\rput{0}(2,13) {{\bf 4}}
\rput{0}(2,14) {{\bf 2}}

\rput{0}(3,0)  {{\bf 1}}
\rput{0}(3,1)  {{\bf  }}
\rput{0}(3,2)  {{\bf  }}
\rput{0}(3,3)  {{\bf  }}
\rput{0}(3,4)  {{\bf  }}
\rput{0}(3,5)  {{\bf  }}
\rput{0}(3,6)  {{\bf  }}
\rput{0}(3,7)  {{\bf 1}}
\rput{0}(3,8)  {{\bf 3}}
\rput{0}(3,9)  {{\bf 2}}
\rput{0}(3,10) {{\bf 3}}
\rput{0}(3,11) {{\bf 1}}
\rput{0}(3,12) {{\bf 4}}
\rput{0}(3,13) {{\bf 2}}
\rput{0}(3,14) {{\bf 3}}

\rput{0}(4,0)  {{\bf  }}
\rput{0}(4,1)  {{\bf  }}
\rput{0}(4,2)  {{\bf  }}
\rput{0}(4,3)  {{\bf  }}
\rput{0}(4,4)  {{\bf  }}
\rput{0}(4,5)  {{\bf  }}
\rput{0}(4,6)  {{\bf 1}}
\rput{0}(4,7)  {{\bf 3}}
\rput{0}(4,8)  {{\bf 2}}
\rput{0}(4,9)  {{\bf 4}}
\rput{0}(4,10) {{\bf 1}}
\rput{0}(4,11) {{\bf 4}}
\rput{0}(4,12) {{\bf 2}}
\rput{0}(4,13) {{\bf 3}}
\rput{0}(4,14) {{\bf 1}}

\rput{0}(5,0)  {{\bf  }}
\rput{0}(5,1)  {{\bf  }}
\rput{0}(5,2)  {{\bf  }}
\rput{0}(5,3)  {{\bf  }}
\rput{0}(5,4)  {{\bf  }}
\rput{0}(5,5)  {{\bf 1}}
\rput{0}(5,6)  {{\bf 3}}
\rput{0}(5,7)  {{\bf 2}}
\rput{0}(5,8)  {{\bf 4}}
\rput{0}(5,9)  {{\bf 1}}
\rput{0}(5,10) {{\bf 3}}
\rput{0}(5,11) {{\bf 2}}
\rput{0}(5,12) {{\bf 4}}
\rput{0}(5,13) {{\bf 1}}
\rput{0}(5,14) {{\bf  }}

\rput{0}(6,0)  {{\bf  }}
\rput{0}(6,1)  {{\bf  }}
\rput{0}(6,2)  {{\bf  }}
\rput{0}(6,3)  {{\bf  }}
\rput{0}(6,4)  {{\bf 2}}
\rput{0}(6,5)  {{\bf 3}}
\rput{0}(6,6)  {{\bf 2}}
\rput{0}(6,7)  {{\bf 4}}
\rput{0}(6,8)  {{\bf 1}}
\rput{0}(6,9)  {{\bf 3}}
\rput{0}(6,10) {{\bf 2}}
\rput{0}(6,11) {{\bf 4}}
\rput{0}(6,12) {{\bf 1}}
\rput{0}(6,13) {{\bf  }}
\rput{0}(6,14) {{\bf  }}

\rput{0}(7,0)  {{\bf  }}
\rput{0}(7,1)  {{\bf  }}
\rput{0}(7,2)  {{\bf  }}
\rput{0}(7,3)  {{\bf 2}}
\rput{0}(7,4)  {{\bf 3}}
\rput{0}(7,5)  {{\bf 1}}
\rput{0}(7,6)  {{\bf 4}}
\rput{0}(7,7)  {{\bf 1}}
\rput{0}(7,8)  {{\bf 3}}
\rput{0}(7,9)  {{\bf 2}}
\rput{0}(7,10) {{\bf 4}}
\rput{0}(7,11) {{\bf 1}}
\rput{0}(7,12) {{\bf  }}
\rput{0}(7,13) {{\bf  }}
\rput{0}(7,14) {{\bf  }}

\rput{0}(8,0)  {{\bf  }}
\rput{0}(8,1)  {{\bf  }}
\rput{0}(8,2)  {{\bf 2}}
\rput{0}(8,3)  {{\bf 3}}
\rput{0}(8,4)  {{\bf 1}}
\rput{0}(8,5)  {{\bf 4}}
\rput{0}(8,6)  {{\bf 2}}
\rput{0}(8,7)  {{\bf 3}}
\rput{0}(8,8)  {{\bf 2}}
\rput{0}(8,9)  {{\bf 4}}
\rput{0}(8,10) {{\bf 1}}
\rput{0}(8,11) {{\bf  }}
\rput{0}(8,12) {{\bf  }}
\rput{0}(8,13) {{\bf  }}
\rput{0}(8,14) {{\bf  }}

\rput{0}(9,0)  {{\bf  }}
\rput{0}(9,1)  {{\bf 2}}
\rput{0}(9,2)  {{\bf 3}}
\rput{0}(9,3)  {{\bf 1}}
\rput{0}(9,4)  {{\bf 4}}
\rput{0}(9,5)  {{\bf 2}}
\rput{0}(9,6)  {{\bf 3}}
\rput{0}(9,7)  {{\bf 1}}
\rput{0}(9,8)  {{\bf 4}}
\rput{0}(9,9)  {{\bf 1}}
\rput{0}(9,10) {{\bf  }}
\rput{0}(9,11) {{\bf  }}
\rput{0}(9,12) {{\bf  }}
\rput{0}(9,13) {{\bf  }}
\rput{0}(9,14) {{\bf  }}

\rput{0}(10,0) {{\bf 2}}
\rput{0}(10,1) {{\bf 4}}
\rput{0}(10,2) {{\bf 1}}
\rput{0}(10,3) {{\bf 4}}
\rput{0}(10,4) {{\bf 2}}
\rput{0}(10,5) {{\bf 3}}
\rput{0}(10,6) {{\bf 1}}
\rput{0}(10,7) {{\bf 4}}
\rput{0}(10,8) {{\bf 2}}
\rput{0}(10,9) {{\bf  }}
\rput{0}(10,10){{\bf  }}
\rput{0}(10,11){{\bf  }}
\rput{0}(10,12){{\bf  }}
\rput{0}(10,13){{\bf  }}
\rput{0}(10,14){{\bf  }}

\rput{0}(11,0) {{\bf 4}}
\rput{0}(11,1) {{\bf 1}}
\rput{0}(11,2) {{\bf 3}}
\rput{0}(11,3) {{\bf 2}}
\rput{0}(11,4) {{\bf 3}}
\rput{0}(11,5) {{\bf 1}}
\rput{0}(11,6) {{\bf 4}}
\rput{0}(11,7) {{\bf 2}}
\rput{0}(11,8) {{\bf  }}
\rput{0}(11,9) {{\bf  }}
\rput{0}(11,10){{\bf  }}
\rput{0}(11,11){{\bf  }}
\rput{0}(11,12){{\bf  }}
\rput{0}(11,13){{\bf  }}
\rput{0}(11,14){{\bf 2}}

\rput{0}(12,0) {{\bf 1}}
\rput{0}(12,1) {{\bf 3}}
\rput{0}(12,2) {{\bf 2}}
\rput{0}(12,3) {{\bf 4}}
\rput{0}(12,4) {{\bf 1}}
\rput{0}(12,5) {{\bf 4}}
\rput{0}(12,6) {{\bf 2}}
\rput{0}(12,7) {{\bf  }}
\rput{0}(12,8) {{\bf  }}
\rput{0}(12,9) {{\bf  }}
\rput{0}(12,10){{\bf  }}
\rput{0}(12,11){{\bf  }}
\rput{0}(12,12){{\bf  }}
\rput{0}(12,13){{\bf 2}}
\rput{0}(12,14){{\bf 4}}

\rput{0}(13,0) {{\bf 3}}
\rput{0}(13,1) {{\bf 2}}
\rput{0}(13,2) {{\bf 4}}
\rput{0}(13,3) {{\bf 1}}
\rput{0}(13,4) {{\bf 3}}
\rput{0}(13,5) {{\bf 2}}
\rput{0}(13,6) {{\bf  }}
\rput{0}(13,7) {{\bf  }}
\rput{0}(13,8) {{\bf  }}
\rput{0}(13,9) {{\bf  }}
\rput{0}(13,10){{\bf  }}
\rput{0}(13,11){{\bf  }}
\rput{0}(13,12){{\bf 1}}
\rput{0}(13,13){{\bf 4}}
\rput{0}(13,14){{\bf 1}}

\rput{0}(14,0) {{\bf 2}}
\rput{0}(14,1) {{\bf 4}}
\rput{0}(14,2) {{\bf 1}}
\rput{0}(14,3) {{\bf 3}}
\rput{0}(14,4) {{\bf 2}}
\rput{0}(14,5) {{\bf  }}
\rput{0}(14,6) {{\bf  }}
\rput{0}(14,7) {{\bf  }}
\rput{0}(14,8) {{\bf  }}
\rput{0}(14,9) {{\bf  }}
\rput{0}(14,10){{\bf  }}
\rput{0}(14,11){{\bf 1}}
\rput{0}(14,12){{\bf 4}}
\rput{0}(14,13){{\bf 2}}
\rput{0}(14,14){{\bf 3}}
%
%
\multirput{0}(14.5,-0.5)(0,1){15}{%
     \psline[linewidth=2pt,linecolor=black]{->}(0.2,-0.2)(-0.1,0.1)
}
\uput[0](14.7,-0.7){D1}
\uput[0](14.7, 0.3){D2}
\uput[0](14.7, 1.3){D3}
\uput[0](14.7, 2.3){D4}
\uput[0](14.7, 3.3){D5}
\uput[0](14.7, 4.3){D6}
\uput[0](14.7, 5.3){D7}
\uput[0](14.7, 6.3){D8}
\uput[0](14.7, 7.3){D9}
\uput[0](14.7, 8.3){D10}
\uput[0](14.7, 9.3){D11}
\uput[0](14.7,10.3){D12}
\uput[0](14.7,11.3){D13}
\uput[0](14.7,12.3){D14}
\uput[0](14.7,13.3){D15}
\endpspicture
\caption{ \label{prop.12k+3.fig1}
The 4-coloring of $T(15,15)$ after Step~1 in the case $L=4k+1$.
}
\end{figure}
%
%

\noindent
{\bf Step 3.}
The last colored counter-diagonals are D$(6k-1)$ and D$(6k+6)$. 
On D$(6k)$ the vertices at $(3k,3k)$ and $(9k+2,9k+1)$ only admit a single
color: either $3$ or $4$. We color the rest of the vertices of D$(6k)$ 
with colors $1$ and $2$ (again, uniquely). 
On D$(6k+5)$ we perform the same procedure; here the vertices with only one
color choice are located at $(3k+3,3k)$ and $(9k+4,9k+4)$.  
The contribution to the degree 
of the newly colored triangles is zero: the partial degree is still  
$\deg f|_R = 8 + 12(k-1)$.

%
%
\begin{figure}[htb]
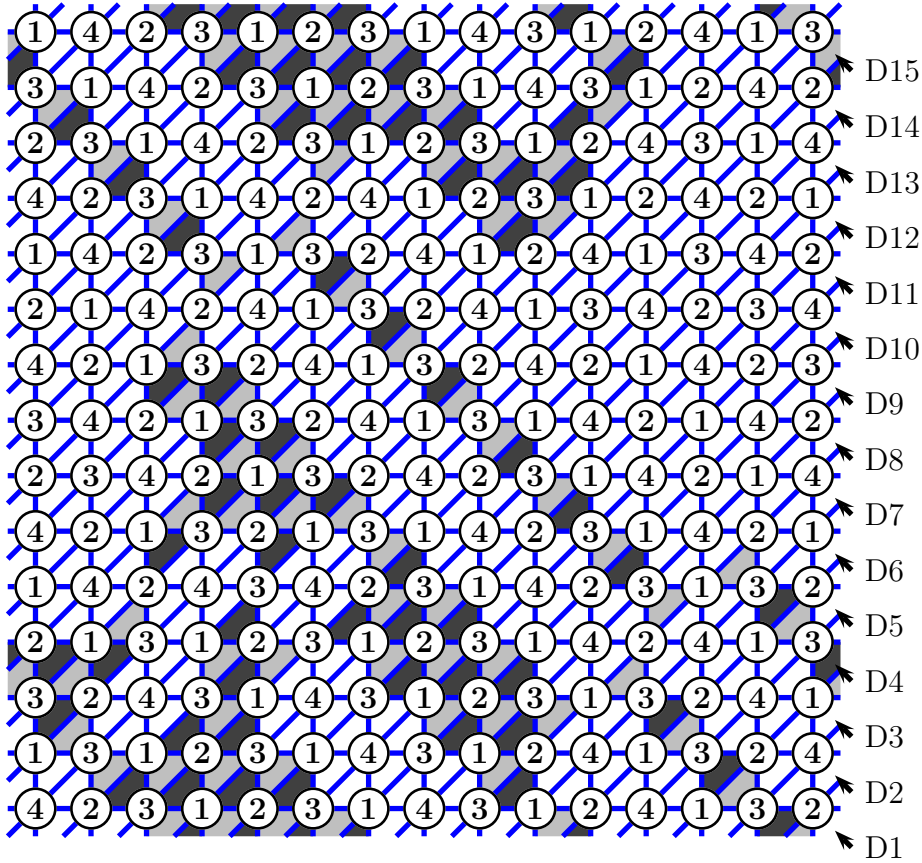

\centering
%
%
\psset{xunit=21pt}
\psset{yunit=21pt}
\psset{labelsep=5pt}
\pspicture(-1,-1)(16,15)
\psline*[linecolor=lightgray](0,2)(0,3)(-0.5,2.5)(-0.5,2)(0,2)
\psline*[linecolor=lightgray](0,1)(2,1)(1,0)(1,2)(0,1)
\psline*[linecolor=lightgray](2,0)(3,0)(2.5,-0.5)(2,-0.5)(2,0)
\psline*[linecolor=lightgray](0,14)(-0.5,14)(-0.5,13.5)(0,14)
\psline*[linecolor=lightgray](0,13)(1,13)(0,12)(0,13)
\psline*[linecolor=lightgray](1,12)(2,12)(1,11)(1,12)
\psline*[linecolor=lightgray](2,11)(3,11)(2,10)(2,11)
\psline*[linecolor=lightgray](3,10)(4,10)(3,9)(3,10)
\psline*[linecolor=lightgray](3,9)(3,8)(2,8)(3,9)
\psline*[linecolor=lightgray](4,8)(4,7)(3,7)(4,8)
\psline*[linecolor=lightgray](5,7)(5,6)(4,6)(5,7)
\psline*[linecolor=lightgray](5,5)(7,5)(6,4)(6,6)(5,5)
\psline*[linecolor=lightgray](7,3)(7,4)(8,4)(7,3)
\psline*[linecolor=lightgray](8,2)(8,3)(9,3)(8,2)
\psline*[linecolor=lightgray](9,1)(9,2)(10,2)(9,1)
\psline*[linecolor=lightgray](10,2)(11,2)(11,4)(12,4)(10,2)
\psline*[linecolor=lightgray](2,14)(2,14.5)(2.5,14.5)(2,14)
\psline*[linecolor=lightgray](3,14)(3,13)(4,14)(3,14)
\psline*[linecolor=lightgray](4,13)(4,12)(5,13)(4,13)
\psline*[linecolor=lightgray](5,12)(5,11)(6,12)(5,12)
\psline*[linecolor=lightgray](4,10)(5,10)(5,11)(4,10)
\psline*[linecolor=lightgray](5,9)(6,9)(6,10)(5,9)
\psline*[linecolor=lightgray](6,8)(7,8)(7,9)(6,8)
\psline*[linecolor=lightgray](7,7)(9,7)(8,6)(8,8)(7,7)
\psline*[linecolor=lightgray](9,5)(9,6)(10,6)(9,5)
\psline*[linecolor=lightgray](10,4)(10,5)(11,5)(10,4)
\psline*[linecolor=lightgray](12,4)(13,4)(13,5)(12,4)
\psline*[linecolor=lightgray](13,3)(14,3)(14,4)(13,3)
\psline*[linecolor=lightgray](14,2)(14.5,2)(14.5,2.5)(14,2)
\psline*[linecolor=lightgray](13,14)(13.5,14.5)(14,14.5)(14,14)(13,14)
\psline*[linecolor=lightgray](14,13)(14,14)(14.5,14)(14.5,13.5)(14,13)
\psline*[linecolor=lightgray](0,2)(1,2)(1,3)(0,2)
\psline*[linecolor=lightgray](2,0)(2,1)(3,1)(2,0)
\psline*[linecolor=lightgray](3,0)(3,1)(4,1)(3,0)
\psline*[linecolor=lightgray](8,0)(8,1)(9,1)(8,0)
\psline*[linecolor=lightgray](8,1)(8,2)(9,2)(8,1)
\psline*[linecolor=lightgray](7,1)(7,2)(8,2)(7,1)
\psline*[linecolor=lightgray](7,2)(7,3)(8,3)(7,2)
\psline*[linecolor=lightgray](6,2)(6,3)(7,3)(6,2)
\psline*[linecolor=lightgray](6,3)(6,4)(7,4)(6,3)
\psline*[linecolor=lightgray](2,7)(3,7)(3,8)(2,7)
\psline*[linecolor=lightgray](3,6)(4,6)(4,7)(3,6)
\psline*[linecolor=lightgray](3,5)(4,5)(4,6)(3,5)
\psline*[linecolor=lightgray](4,5)(5,5)(5,6)(4,5)
\psline*[linecolor=lightgray](4,14)(4,14.5)(4.5,14.5)(4,14)
\psline*[linecolor=lightgray](4,0)(4,-0.5)(4.5,-0.5)(5,0)(4,0)
\psline*[linecolor=lightgray](3,14)(3,14.5)(3.5,14.5)(3,14)
\psline*[linecolor=lightgray](3,0)(3,-0.5)(3.5,-0.5)(4,0)(3,0)
\psline*[linecolor=lightgray](4,14)(5,14)(4,13)(4,14)
\psline*[linecolor=lightgray](5,14)(6,14)(5,13)(5,14)
\psline*[linecolor=lightgray](5,13)(6,13)(5,12)(5,13)
\psline*[linecolor=lightgray](6,13)(7,13)(6,12)(6,13)
\psline*[linecolor=lightgray](7,12)(8,12)(7,11)(7,12)
\psline*[linecolor=lightgray](8,11)(9,11)(8,10)(8,11)
\psline*[linecolor=lightgray](9,14)(9,14.5)(9.5,14.5)(9,14)
\psline*[linecolor=lightgray](10,13)(10,14)(11,14)(10,13)
\psline*[linecolor=lightgray](9,0)(10,0)(9.5,-0.5)(9,-0.5)(9,0)
\psline*[linecolor=lightgray](1,3)(2,3)(2,4)(1,3)
\psline*[linecolor=lightgray](2,5)(3,5)(3,6)(2,5)
\psline*[linecolor=lightgray](5,0)(6,0)(5.5,-0.5)(5,-0.5)(5,0)
\psline*[linecolor=lightgray](4,0)(4,1)(5,1)(4,0)
\psline*[linecolor=lightgray](3,1)(3,2)(4,2)(3,1)
\psline*[linecolor=lightgray](3,2)(3,3)(4,3)(3,2)
\psline*[linecolor=lightgray](4,2)(4,3)(5,3)(4,2)
\psline*[linecolor=lightgray](5,14)(5,14.5)(5.5,14.5)(5,14)
\psline*[linecolor=lightgray](6,13)(6,14)(7,14)(6,13)
\psline*[linecolor=lightgray](7,12)(7,13)(8,13)(7,12)
\psline*[linecolor=lightgray](8,11)(8,12)(9,12)(8,11)
\psline*[linecolor=lightgray](9,10)(9,11)(10,11)(9,10)
\psline*[linecolor=lightgray](9,11)(9,12)(10,12)(9,11)
\psline*[linecolor=lightgray](10,12)(10,13)(11,13)(10,12)
\psline*[linecolor=lightgray](11,1)(12,1)(12,2)(11,1)
\psline*[linecolor=lightgray](12,0)(13,0)(13,1)(12,0)
\psline*[linecolor=lightgray](14,0)(14,-0.5)(13.5,-0.5)(14,0)

\psline*[linecolor=darkgray](0,3)(-0.5,3)(-0.5,2.5)(0,3)
\psline*[linecolor=darkgray](0,2)(1,2)(0,1)(0,2)
\psline*[linecolor=darkgray](1,0)(2,0)(2,1)(1,0)
\psline*[linecolor=darkgray](3,0)(3,-0.5)(2.5,-0.5)(3,0)
\psline*[linecolor=darkgray](0,13)(-0.5,13)(-0.5,13.5)(0,14)(0,13)
\psline*[linecolor=darkgray](0,12)(1,12)(1,13)(0,12)
\psline*[linecolor=darkgray](1,11)(2,11)(2,12)(1,11)
\psline*[linecolor=darkgray](2,10)(3,10)(3,11)(2,10)
\psline*[linecolor=darkgray](3,7)(3,8)(4,8)(3,7)
\psline*[linecolor=darkgray](4,6)(4,7)(5,7)(4,6)
\psline*[linecolor=darkgray](5,5)(5,6)(6,6)(5,5)
\psline*[linecolor=darkgray](6,4)(7,4)(7,5)(6,4)
\psline*[linecolor=darkgray](7,3)(8,3)(8,4)(7,3)
\psline*[linecolor=darkgray](8,2)(9,2)(9,3)(8,2)
\psline*[linecolor=darkgray](2,14)(2.5,14.5)(3,14.5)(3,14)(2,14)
\psline*[linecolor=darkgray](3,13)(4,13)(4,14)(3,13)
\psline*[linecolor=darkgray](4,12)(5,12)(5,13)(4,12)
\psline*[linecolor=darkgray](5,9)(5,10)(6,10)(5,9)
\psline*[linecolor=darkgray](6,8)(6,9)(7,9)(6,8)
\psline*[linecolor=darkgray](7,7)(7,8)(8,8)(7,7)
\psline*[linecolor=darkgray](8,6)(9,6)(9,7)(8,6)
\psline*[linecolor=darkgray](9,5)(10,5)(10,6)(9,5)
\psline*[linecolor=darkgray](10,4)(11,4)(11,5)(10,4)
\psline*[linecolor=darkgray](13,3)(13,4)(14,4)(13,3)
\psline*[linecolor=darkgray](14,2)(14,3)(14.5,3)(14.5,2.5)(14,2)
\psline*[linecolor=darkgray](13,14)(13,14.5)(13.5,14.5)(13,14)
\psline*[linecolor=darkgray](14,13)(14.5,13)(14.5,13.5)(14,13)
\psline*[linecolor=darkgray](0,2)(0,3)(1,3)(0,2)
\psline*[linecolor=darkgray](1,2)(1,3)(2,3)(1,2)
\psline*[linecolor=darkgray](2,0)(3,0)(3,1)(2,0)
\psline*[linecolor=darkgray](2,1)(3,1)(3,2)(2,1)
\psline*[linecolor=darkgray](3,0)(4,0)(4,1)(3,0)
\psline*[linecolor=darkgray](8,0)(9,0)(9,1)(8,0)
\psline*[linecolor=darkgray](8,1)(9,1)(9,2)(8,1)
\psline*[linecolor=darkgray](7,1)(8,1)(8,2)(7,1)
\psline*[linecolor=darkgray](7,2)(8,2)(8,3)(7,2)
\psline*[linecolor=darkgray](6,2)(7,2)(7,3)(6,2)
\psline*[linecolor=darkgray](6,3)(7,3)(7,4)(6,3)
\psline*[linecolor=darkgray](5,3)(6,3)(6,4)(5,3)
\psline*[linecolor=darkgray](2,7)(2,8)(3,8)(2,7)
\psline*[linecolor=darkgray](3,6)(3,7)(4,7)(3,6)
\psline*[linecolor=darkgray](3,5)(3,6)(4,6)(3,5)
\psline*[linecolor=darkgray](4,5)(4,6)(5,6)(4,5)
\psline*[linecolor=darkgray](4,4)(4,5)(5,5)(4,4)
\psline*[linecolor=darkgray](4,0)(4,-0.5)(3.5,-0.5)(4,0)
\psline*[linecolor=darkgray](5,0)(5,-0.5)(4.5,-0.5)(5,0)
\psline*[linecolor=darkgray](3,14)(3.5,14.5)(4,14.5)(4,14)(3,14)
\psline*[linecolor=darkgray](4,14)(4.5,14.5)(5,14.5)(5,14)(4,14)
\psline*[linecolor=darkgray](4,13)(5,13)(5,14)(4,13)
\psline*[linecolor=darkgray](5,13)(6,13)(6,14)(5,13)
\psline*[linecolor=darkgray](5,12)(6,12)(6,13)(5,12)
\psline*[linecolor=darkgray](6,12)(7,12)(7,13)(6,12)
\psline*[linecolor=darkgray](7,11)(8,11)(8,12)(7,11)
\psline*[linecolor=darkgray](8,10)(9,10)(9,11)(8,10)
\psline*[linecolor=darkgray](10,13)(11,13)(11,14)(10,13)
\psline*[linecolor=darkgray](9,14)(10,14)(10,14.5)(9.5,14.5)(9,14)
\psline*[linecolor=darkgray](10,0)(10,-0.5)(9.5,-0.5)(10,0)
\psline*[linecolor=darkgray](2,4)(2,5)(3,5)(2,4)
\psline*[linecolor=darkgray](3,3)(4,3)(4,4)(3,3)
\psline*[linecolor=darkgray](3,2)(4,2)(4,3)(3,2)
\psline*[linecolor=darkgray](3,1)(4,1)(4,2)(3,1)
\psline*[linecolor=darkgray](4,0)(5,0)(5,1)(4,0)
\psline*[linecolor=darkgray](6,0)(6,-0.5)(5.5,-0.5)(6,0)
\psline*[linecolor=darkgray](5,14)(5.5,14.5)(6,14.5)(6,14)(5,14)
\psline*[linecolor=darkgray](6,13)(7,13)(7,14)(6,13)
\psline*[linecolor=darkgray](7,12)(8,12)(8,13)(7,12)
\psline*[linecolor=darkgray](8,11)(9,11)(9,12)(8,11)
\psline*[linecolor=darkgray](9,11)(10,11)(10,12)(9,11)
\psline*[linecolor=darkgray](9,12)(10,12)(10,13)(9,12)
\psline*[linecolor=darkgray](11,1)(11,2)(12,2)(11,1)
\psline*[linecolor=darkgray](12,0)(12,1)(13,1)(12,0)
\psline*[linecolor=darkgray](13,0)(14,0)(13.5,-0.5)(13,-0.5)(13,0)

\psline[linewidth=2pt,linecolor=blue](-0.5,0)(14.5,0)
\psline[linewidth=2pt,linecolor=blue](-0.5,1)(14.5,1)
\psline[linewidth=2pt,linecolor=blue](-0.5,2)(14.5,2)
\psline[linewidth=2pt,linecolor=blue](-0.5,3)(14.5,3)
\psline[linewidth=2pt,linecolor=blue](-0.5,4)(14.5,4)
\psline[linewidth=2pt,linecolor=blue](-0.5,5)(14.5,5)
\psline[linewidth=2pt,linecolor=blue](-0.5,6)(14.5,6)
\psline[linewidth=2pt,linecolor=blue](-0.5,7)(14.5,7)
\psline[linewidth=2pt,linecolor=blue](-0.5,8)(14.5,8)
\psline[linewidth=2pt,linecolor=blue](-0.5,9)(14.5,9)
\psline[linewidth=2pt,linecolor=blue](-0.5,10)(14.5,10)
\psline[linewidth=2pt,linecolor=blue](-0.5,11)(14.5,11)
\psline[linewidth=2pt,linecolor=blue](-0.5,12)(14.5,12)
\psline[linewidth=2pt,linecolor=blue](-0.5,13)(14.5,13)
\psline[linewidth=2pt,linecolor=blue](-0.5,14)(14.5,14)
\psline[linewidth=2pt,linecolor=blue](0,-0.5)(0,14.5)
\psline[linewidth=2pt,linecolor=blue](1,-0.5)(1,14.5)
\psline[linewidth=2pt,linecolor=blue](2,-0.5)(2,14.5)
\psline[linewidth=2pt,linecolor=blue](3,-0.5)(3,14.5)
\psline[linewidth=2pt,linecolor=blue](4,-0.5)(4,14.5)
\psline[linewidth=2pt,linecolor=blue](5,-0.5)(5,14.5)
\psline[linewidth=2pt,linecolor=blue](6,-0.5)(6,14.5)
\psline[linewidth=2pt,linecolor=blue](7,-0.5)(7,14.5)
\psline[linewidth=2pt,linecolor=blue](8,-0.5)(8,14.5)
\psline[linewidth=2pt,linecolor=blue](9,-0.5)(9,14.5)
\psline[linewidth=2pt,linecolor=blue](10,-0.5)(10,14.5)
\psline[linewidth=2pt,linecolor=blue](11,-0.5)(11,14.5)
\psline[linewidth=2pt,linecolor=blue](12,-0.5)(12,14.5)
\psline[linewidth=2pt,linecolor=blue](13,-0.5)(13,14.5)
\psline[linewidth=2pt,linecolor=blue](14,-0.5)(14,14.5)
\psline[linewidth=2pt,linecolor=blue](-0.5,-0.5)(14.5,14.5)
\psline[linewidth=2pt,linecolor=blue](0.5,-0.5)(14.5,13.5)
\psline[linewidth=2pt,linecolor=blue](1.5,-0.5)(14.5,12.5)
\psline[linewidth=2pt,linecolor=blue](2.5,-0.5)(14.5,11.5)
\psline[linewidth=2pt,linecolor=blue](3.5,-0.5)(14.5,10.5)
\psline[linewidth=2pt,linecolor=blue](4.5,-0.5)(14.5,9.5)
\psline[linewidth=2pt,linecolor=blue](5.5,-0.5)(14.5,8.5)
\psline[linewidth=2pt,linecolor=blue](6.5,-0.5)(14.5,7.5)
\psline[linewidth=2pt,linecolor=blue](7.5,-0.5)(14.5,6.5)
\psline[linewidth=2pt,linecolor=blue](8.5,-0.5)(14.5,5.5)
\psline[linewidth=2pt,linecolor=blue](9.5,-0.5)(14.5,4.5)
\psline[linewidth=2pt,linecolor=blue](10.5,-0.5)(14.5,3.5)
\psline[linewidth=2pt,linecolor=blue](11.5,-0.5)(14.5,2.5)
\psline[linewidth=2pt,linecolor=blue](12.5,-0.5)(14.5,1.5)
\psline[linewidth=2pt,linecolor=blue](13.5,-0.5)(14.5,0.5)
\psline[linewidth=2pt,linecolor=blue](-0.5,0.5)(13.5,14.5)
\psline[linewidth=2pt,linecolor=blue](-0.5,1.5)(12.5,14.5)
\psline[linewidth=2pt,linecolor=blue](-0.5,2.5)(11.5,14.5)
\psline[linewidth=2pt,linecolor=blue](-0.5,3.5)(10.5,14.5)
\psline[linewidth=2pt,linecolor=blue](-0.5,4.5)(9.5,14.5)
\psline[linewidth=2pt,linecolor=blue](-0.5,5.5)(8.5,14.5)
\psline[linewidth=2pt,linecolor=blue](-0.5,6.5)(7.5,14.5)
\psline[linewidth=2pt,linecolor=blue](-0.5,7.5)(6.5,14.5)
\psline[linewidth=2pt,linecolor=blue](-0.5,8.5)(5.5,14.5)
\psline[linewidth=2pt,linecolor=blue](-0.5,9.5)(4.5,14.5)
\psline[linewidth=2pt,linecolor=blue](-0.5,10.5)(3.5,14.5)
\psline[linewidth=2pt,linecolor=blue](-0.5,11.5)(2.5,14.5)
\psline[linewidth=2pt,linecolor=blue](-0.5,12.5)(1.5,14.5)
\psline[linewidth=2pt,linecolor=blue](-0.5,13.5)(0.5,14.5)
\multirput{0}(0,0)(0,1){15}{%
  \multirput{0}(0,0)(1,0){15}{%
     \pscircle*[linecolor=white]{8pt}
     \pscircle[linewidth=1pt,linecolor=black] {8pt}
   }
}
\rput{0}(0,0)  {{\bf 4}}
\rput{0}(0,1)  {{\bf 1}}
\rput{0}(0,2)  {{\bf 3}}
\rput{0}(0,3)  {{\bf 2}}
\rput{0}(0,4)  {{\bf 1}}
\rput{0}(0,5)  {{\bf 4}}
\rput{0}(0,6)  {{\bf 2}}
\rput{0}(0,7)  {{\bf 3}}
\rput{0}(0,8)  {{\bf 4}}
\rput{0}(0,9)  {{\bf 2}}
\rput{0}(0,10) {{\bf 1}}
\rput{0}(0,11) {{\bf 4}}
\rput{0}(0,12) {{\bf 2}}
\rput{0}(0,13) {{\bf 3}}
\rput{0}(0,14) {{\bf 1}}

\rput{0}(1,0)  {{\bf 2}}
\rput{0}(1,1)  {{\bf 3}}
\rput{0}(1,2)  {{\bf 2}}
\rput{0}(1,3)  {{\bf 1}}
\rput{0}(1,4)  {{\bf 4}}
\rput{0}(1,5)  {{\bf 2}}
\rput{0}(1,6)  {{\bf 3}}
\rput{0}(1,7)  {{\bf 4}}
\rput{0}(1,8)  {{\bf 2}}
\rput{0}(1,9)  {{\bf 1}}
\rput{0}(1,10) {{\bf 4}}
\rput{0}(1,11) {{\bf 2}}
\rput{0}(1,12) {{\bf 3}}
\rput{0}(1,13) {{\bf 1}}
\rput{0}(1,14) {{\bf 4}}

\rput{0}(2,0)  {{\bf 3}}
\rput{0}(2,1)  {{\bf 1}}
\rput{0}(2,2)  {{\bf 4}}
\rput{0}(2,3)  {{\bf 3}}
\rput{0}(2,4)  {{\bf 2}}
\rput{0}(2,5)  {{\bf 1}}
\rput{0}(2,6)  {{\bf 4}}
\rput{0}(2,7)  {{\bf 2}}
\rput{0}(2,8)  {{\bf 1}}
\rput{0}(2,9)  {{\bf 4}}
\rput{0}(2,10) {{\bf 2}}
\rput{0}(2,11) {{\bf 3}}
\rput{0}(2,12) {{\bf 1}}
\rput{0}(2,13) {{\bf 4}}
\rput{0}(2,14) {{\bf 2}}

\rput{0}(3,0)  {{\bf 1}}
\rput{0}(3,1)  {{\bf 2}}
\rput{0}(3,2)  {{\bf 3}}
\rput{0}(3,3)  {{\bf 1}}
\rput{0}(3,4)  {{\bf 4}}
\rput{0}(3,5)  {{\bf 3}}
\rput{0}(3,6)  {{\bf 2}}
\rput{0}(3,7)  {{\bf 1}}
\rput{0}(3,8)  {{\bf 3}}
\rput{0}(3,9)  {{\bf 2}}
\rput{0}(3,10) {{\bf 3}}
\rput{0}(3,11) {{\bf 1}}
\rput{0}(3,12) {{\bf 4}}
\rput{0}(3,13) {{\bf 2}}
\rput{0}(3,14) {{\bf 3}}

\rput{0}(4,0)  {{\bf 2}}
\rput{0}(4,1)  {{\bf 3}}
\rput{0}(4,2)  {{\bf 1}}
\rput{0}(4,3)  {{\bf 2}}
\rput{0}(4,4)  {{\bf 3}}
\rput{0}(4,5)  {{\bf 2}}
\rput{0}(4,6)  {{\bf 1}}
\rput{0}(4,7)  {{\bf 3}}
\rput{0}(4,8)  {{\bf 2}}
\rput{0}(4,9)  {{\bf 4}}
\rput{0}(4,10) {{\bf 1}}
\rput{0}(4,11) {{\bf 4}}
\rput{0}(4,12) {{\bf 2}}
\rput{0}(4,13) {{\bf 3}}
\rput{0}(4,14) {{\bf 1}}

\rput{0}(5,0)  {{\bf 3}}
\rput{0}(5,1)  {{\bf 1}}
\rput{0}(5,2)  {{\bf 4}}
\rput{0}(5,3)  {{\bf 3}}
\rput{0}(5,4)  {{\bf 4}}
\rput{0}(5,5)  {{\bf 1}}
\rput{0}(5,6)  {{\bf 3}}
\rput{0}(5,7)  {{\bf 2}}
\rput{0}(5,8)  {{\bf 4}}
\rput{0}(5,9)  {{\bf 1}}
\rput{0}(5,10) {{\bf 3}}
\rput{0}(5,11) {{\bf 2}}
\rput{0}(5,12) {{\bf 3}}
\rput{0}(5,13) {{\bf 1}}
\rput{0}(5,14) {{\bf 2}}

\rput{0}(6,0)  {{\bf 1}}
\rput{0}(6,1)  {{\bf 4}}
\rput{0}(6,2)  {{\bf 3}}
\rput{0}(6,3)  {{\bf 1}}
\rput{0}(6,4)  {{\bf 2}}
\rput{0}(6,5)  {{\bf 3}}
\rput{0}(6,6)  {{\bf 2}}
\rput{0}(6,7)  {{\bf 4}}
\rput{0}(6,8)  {{\bf 1}}
\rput{0}(6,9)  {{\bf 3}}
\rput{0}(6,10) {{\bf 2}}
\rput{0}(6,11) {{\bf 4}}
\rput{0}(6,12) {{\bf 1}}
\rput{0}(6,13) {{\bf 2}}
\rput{0}(6,14) {{\bf 3}}

\rput{0}(7,0)  {{\bf 4}}
\rput{0}(7,1)  {{\bf 3}}
\rput{0}(7,2)  {{\bf 1}}
\rput{0}(7,3)  {{\bf 2}}
\rput{0}(7,4)  {{\bf 3}}
\rput{0}(7,5)  {{\bf 1}}
\rput{0}(7,6)  {{\bf 4}}
\rput{0}(7,7)  {{\bf 1}}
\rput{0}(7,8)  {{\bf 3}}
\rput{0}(7,9)  {{\bf 2}}
\rput{0}(7,10) {{\bf 4}}
\rput{0}(7,11) {{\bf 1}}
\rput{0}(7,12) {{\bf 2}}
\rput{0}(7,13) {{\bf 3}}
\rput{0}(7,14) {{\bf 1}}

\rput{0}(8,0)  {{\bf 3}}
\rput{0}(8,1)  {{\bf 1}}
\rput{0}(8,2)  {{\bf 2}}
\rput{0}(8,3)  {{\bf 3}}
\rput{0}(8,4)  {{\bf 1}}
\rput{0}(8,5)  {{\bf 4}}
\rput{0}(8,6)  {{\bf 2}}
\rput{0}(8,7)  {{\bf 3}}
\rput{0}(8,8)  {{\bf 2}}
\rput{0}(8,9)  {{\bf 4}}
\rput{0}(8,10) {{\bf 1}}
\rput{0}(8,11) {{\bf 2}}
\rput{0}(8,12) {{\bf 3}}
\rput{0}(8,13) {{\bf 1}}
\rput{0}(8,14) {{\bf 4}}

\rput{0}(9,0)  {{\bf 1}}
\rput{0}(9,1)  {{\bf 2}}
\rput{0}(9,2)  {{\bf 3}}
\rput{0}(9,3)  {{\bf 1}}
\rput{0}(9,4)  {{\bf 4}}
\rput{0}(9,5)  {{\bf 2}}
\rput{0}(9,6)  {{\bf 3}}
\rput{0}(9,7)  {{\bf 1}}
\rput{0}(9,8)  {{\bf 4}}
\rput{0}(9,9)  {{\bf 1}}
\rput{0}(9,10) {{\bf 2}}
\rput{0}(9,11) {{\bf 3}}
\rput{0}(9,12) {{\bf 1}}
\rput{0}(9,13) {{\bf 4}}
\rput{0}(9,14) {{\bf 3}}

\rput{0}(10,0) {{\bf 2}}
\rput{0}(10,1) {{\bf 4}}
\rput{0}(10,2) {{\bf 1}}
\rput{0}(10,3) {{\bf 4}}
\rput{0}(10,4) {{\bf 2}}
\rput{0}(10,5) {{\bf 3}}
\rput{0}(10,6) {{\bf 1}}
\rput{0}(10,7) {{\bf 4}}
\rput{0}(10,8) {{\bf 2}}
\rput{0}(10,9) {{\bf 3}}
\rput{0}(10,10){{\bf 4}}
\rput{0}(10,11){{\bf 1}}
\rput{0}(10,12){{\bf 2}}
\rput{0}(10,13){{\bf 3}}
\rput{0}(10,14){{\bf 1}}

\rput{0}(11,0) {{\bf 4}}
\rput{0}(11,1) {{\bf 1}}
\rput{0}(11,2) {{\bf 3}}
\rput{0}(11,3) {{\bf 2}}
\rput{0}(11,4) {{\bf 3}}
\rput{0}(11,5) {{\bf 1}}
\rput{0}(11,6) {{\bf 4}}
\rput{0}(11,7) {{\bf 2}}
\rput{0}(11,8) {{\bf 1}}
\rput{0}(11,9) {{\bf 4}}
\rput{0}(11,10){{\bf 1}}
\rput{0}(11,11){{\bf 2}}
\rput{0}(11,12){{\bf 4}}
\rput{0}(11,13){{\bf 1}}
\rput{0}(11,14){{\bf 2}}

\rput{0}(12,0) {{\bf 1}}
\rput{0}(12,1) {{\bf 3}}
\rput{0}(12,2) {{\bf 2}}
\rput{0}(12,3) {{\bf 4}}
\rput{0}(12,4) {{\bf 1}}
\rput{0}(12,5) {{\bf 4}}
\rput{0}(12,6) {{\bf 2}}
\rput{0}(12,7) {{\bf 1}}
\rput{0}(12,8) {{\bf 4}}
\rput{0}(12,9) {{\bf 2}}
\rput{0}(12,10){{\bf 3}}
\rput{0}(12,11){{\bf 4}}
\rput{0}(12,12){{\bf 3}}
\rput{0}(12,13){{\bf 2}}
\rput{0}(12,14){{\bf 4}}

\rput{0}(13,0) {{\bf 3}}
\rput{0}(13,1) {{\bf 2}}
\rput{0}(13,2) {{\bf 4}}
\rput{0}(13,3) {{\bf 1}}
\rput{0}(13,4) {{\bf 3}}
\rput{0}(13,5) {{\bf 2}}
\rput{0}(13,6) {{\bf 1}}
\rput{0}(13,7) {{\bf 4}}
\rput{0}(13,8) {{\bf 2}}
\rput{0}(13,9) {{\bf 3}}
\rput{0}(13,10){{\bf 4}}
\rput{0}(13,11){{\bf 2}}
\rput{0}(13,12){{\bf 1}}
\rput{0}(13,13){{\bf 4}}
\rput{0}(13,14){{\bf 1}}

\rput{0}(14,0) {{\bf 2}}
\rput{0}(14,1) {{\bf 4}}
\rput{0}(14,2) {{\bf 1}}
\rput{0}(14,3) {{\bf 3}}
\rput{0}(14,4) {{\bf 2}}
\rput{0}(14,5) {{\bf 1}}
\rput{0}(14,6) {{\bf 4}}
\rput{0}(14,7) {{\bf 2}}
\rput{0}(14,8) {{\bf 3}}
\rput{0}(14,9) {{\bf 4}}
\rput{0}(14,10){{\bf 2}}
\rput{0}(14,11){{\bf 1}}
\rput{0}(14,12){{\bf 4}}
\rput{0}(14,13){{\bf 2}}
\rput{0}(14,14){{\bf 3}}
%
%
\multirput{0}(14.5,-0.5)(0,1){15}{%
     \psline[linewidth=2pt,linecolor=black]{->}(0.2,-0.2)(-0.1,0.1)
}
\uput[0](14.7,-0.7){D1}
\uput[0](14.7, 0.3){D2}
\uput[0](14.7, 1.3){D3}
\uput[0](14.7, 2.3){D4}
\uput[0](14.7, 3.3){D5}
\uput[0](14.7, 4.3){D6}
\uput[0](14.7, 5.3){D7}
\uput[0](14.7, 6.3){D8}
\uput[0](14.7, 7.3){D9}
\uput[0](14.7, 8.3){D10}
\uput[0](14.7, 9.3){D11}
\uput[0](14.7,10.3){D12}
\uput[0](14.7,11.3){D13}
\uput[0](14.7,12.3){D14}
\uput[0](14.7,13.3){D15}
\endpspicture
\caption{ \label{prop.12k+3.fig3}
The 4-coloring of $T(15,15)$ after Step~4 in the case $L=4k+1$.
}
\end{figure}
%
%

On D$(6k+1)$ there are two pairs of nearby vertices which only admit
one color among $3$ and $4$. One pair is $(3k+1,3k)$ and 
$(3k,3k+1)$; the other one is $(9k+3,9k+1)$ and $(9k+2,9k+2)$. 
We color the other vertices on D$(6k+1)$ by colors $3$ and $4$ 
while using the following rule: those with $x$-coordinate satisfying
$3k+1<x<9k+2$ are colored $3$ (resp.\ $4$) if $k$ is odd (resp.\ even). 
At the end, there are $6k+2$ and $6k+1$ vertices colored alike on D$(6k+1)$. 

On D$(6k+4)$ we also find two pairs of vertices which only admit one color
among $3$ and $4$: one pair is $(3k+3,3k+1)$ and $(3k+2,3k+2)$; the other 
one is $(9k+4,9k+3)$ and $(9k+3,9k+4)$. The other vertices on D$(6k+4)$ are
then colored $3$ and $4$ with the help of the following rules: 1) those
with $x$-coordinate satisfying $3k+3<x<9k+3$ are colored $3$ (resp.\ $4$) 
if $k$ is odd (resp.\  even); 2) the number of vertices colored $3$ is
the same as on D$(6k+1)$. This second rule is used to determine the color
of the vertex at $(3k+1,3k+3)$.  
 
The contribution to the partial degree of these
new triangles is $-4$; thus, the partial degree of $f$ is  
$\deg f|_R = 4 + 12(k-1)$.

\medskip
\noindent
{\bf Step 4.}
On D$(6k+2)$ there are two vertices located at $(3k,3k+2)$ and $(9k+2,9k+3)$
whose colors are fixed to either $1$ or $2$. Color with the same color as
$(9k+2,9k+3)$ the two vertices $(3k+1,3k+1)$ and $(9k+3,9k+2)$. 
At the end, there are $6k+4$ vertices having one color, and $6k+1$ having 
the other one. 

On D$(6k+3)$ there are two vertices whose colors are fixed to either $3$ or 
$4$. There are also four additional vertices whose colors are fixed to either
$1$ or  $2$. These six vertices are located at 
$(3k+2,3k+1)$, $(3k+1,3k+2)$, $(3k,3k+3)$, $(9k+4,9k+2)$, $(9k+3,9k+3)$,
and $(9k+2,9k+4)$. The other vertices on D$(6k+3)$ are colored $3$ or $4$
(the choice for each vertex is unique).  

In Figure~\ref{prop.12k+3.fig3} the final coloring $f$ is depicted.
The increment in the partial degree is $2$. Therefore,
\begin{equation}
\deg f \;=\; 6 + 12(k-1) \;\equiv\; 6 \pmod{12}. \nonumber
\end{equation}
The coloring $f$ of $T(12k+3,12k+3)$ is proper and its degree is congruent
to 6 modulo $12$, as claimed. This completes the proof. \qed

%
%
\section{Further results for $\bm{T(3L,3M)}$} \label{sec.asym} 

In the previous section we have proven that $T(3L,3L)$ has 
at least one coloring with degree $\equiv 6 \pmod{12}$ for any $L\geq 2$,
and hence $\Kc(T(3L,3L),4)>1$. This result can be used for some other
triangulations with aspect ratio different from $1$:

\begin{theorem} \label{theo.cases}
The number of Kempe equivalence classes $\Kc(T,4)$ is at least two 
for any triangulation $T(3Lp,3Lq)$ for $L\geq 2$ and any odd integers 
$p,q$. 
\end{theorem}

\proof
Theorem~\ref{theo.main} shows that there is a coloring $f$ of 
$T(3L,3L)$ for $L\geq 2$ with $\deg(f)\equiv 6 \pmod{12}$. Then,
Lemma~\ref{lemma.tech}(c) proves the claimed result. \qed 

\medskip

In order to obtain more general results, it is convenient to prove the 
following simple proposition.

\begin{proposition} \label{prop.T_Lx3}
The degree of any four-coloring of any triangulation $T(L,3)$ or $T(3,L)$ 
with $L\geq 1$ is zero.
\end{proposition}

\begin{figure}[htb]
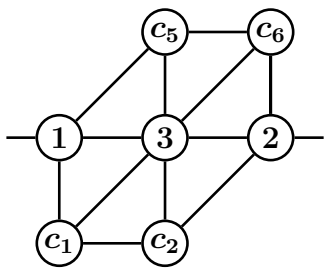

\centering
\psset{xunit=40pt}
\psset{yunit=40pt}
\pspicture(-0.5,-0.5)(2.5,2.5)
\psline[linewidth=1pt](0,0)(0,1)(1,1)(0,0)(1,0)(1,1)(1,2)(0,1)
\psline[linewidth=1pt](1,0)(2,1)(1,1)
\psline[linewidth=1pt](1,1)(2,2)(2,1)(2,2)(1,2)
\psline[linewidth=1pt](0,1)(-0.5,1)
\psline[linewidth=1pt](2,1)(2.5,1)
\multirput{0}(0,0)(1,0){2}{%
\rput{0}(0,0){%
 \pscircle*[linewidth=1pt,linecolor=white]{9pt}%
 \pscircle[linewidth=1pt]{9pt}%
}}
\multirput{0}(0,1)(1,0){3}{%
\rput{0}(0,0){%
 \pscircle*[linewidth=1pt,linecolor=white]{9pt}%
 \pscircle[linewidth=1pt]{9pt}%
}}
\multirput{0}(1,2)(1,0){2}{%
\rput{0}(0,0){%
 \pscircle*[linewidth=1pt,linecolor=white]{9pt}%
 \pscircle[linewidth=1pt]{9pt}%
}}
\rput[c]{0}(0,0){$\bm{c_1}$}
\rput[c]{0}(1,0){$\bm{c_2}$}
\rput[c]{0}(0,1){$\bm{1}$}
\rput[c]{0}(1,1){$\bm{3}$}
\rput[c]{0}(1,2){$\bm{c_5}$}
\rput[c]{0}(2,1){$\bm{2}$}
\rput[c]{0}(2,2){$\bm{c_6}$}
\endpspicture
\caption{ \label{prop.T_Lx3.fig}
Subset of the triangulation $T(3,L)$ used in the proof of 
Proposition~\ref{prop.T_Lx3}.
}
\end{figure}

\proof
Suppose we compute the degree of a given 4-coloring $c$ of the triangulation
$T(3,L)$ by counting those triangular faces colored $123$. We can focus
on those sites colored $3$. Let us suppose the vertex $x$ is colored $3$.
Because the 4-coloring $c$ is proper, none of the neighbors of $x$ 
can be colored $3$.
And because the triangulation has width $3$, the two neighbors along the 
horizontal axis are also adjacent to each other, so they have different colors,
say 1 and 2. This situation is depicted in Figure~\ref{prop.T_Lx3.fig}.
There are only $9$ different four-colorings of the above graph, and all of
them contribute zero to the degree.
Therefore, the contribution of all vertices colored $3$ to the degree 
is zero, and the claimed result is proven. \qed 

%
%


The following lemma shows how to build a four-coloring of the triangulation
$T(L,M+3)$ by ``gluing'' four-colorings of the triangulations $T(L,M)$ and 
$T(L,3)$ that have the same coloring on the top row. One key point is that
the degree is an invariant under this operation. 

\begin{lemma} \label{lemma_plusLx3}
Let us suppose that $c$ is a four-coloring of a triangulation $T(L,M)$
with degree $d$, and that the coloring on the top row is $c_{\rm top}$.  
Let us further suppose there exists a four-coloring $c'$ of the triangulation
$T(L,3)$ with the same coloring on the top row $c'_{\rm top} = c_{\rm top}$. 
Then, there exists a four-coloring of the triangulation $T(L,M+3)$ 
with degree $d$. 
\end{lemma}

\proof
Because both $T(L,M)$ and $T(L,3)$ are triangulations of a torus with the same
width $L$, and the corresponding colorings $c$ and $c'$ 
both have the same top-row coloring $c_{\rm top}$, we can
obtain a four-coloring $c''$ of the triangulation $T(L,M+3)$ by
``gluing'' together these two colorings.
This is indeed a proper coloring of $T(L,M+3)$ and its 
degree can be computed as $\deg(c'')=\deg(c)+\deg(c')=\deg(c)=d$,
since $\deg(c')=0$ by Proposition~\ref{prop.T_Lx3}.
\qed

\medskip

This lemma gives us the opportunity to devise an inductive proof that 
there is a four-coloring with degree $6\pmod{12}$ for any triangulation
$T(3L,3M)$ with $M\geq L$. The base case $L=M$ is already verified by 
Theorem~\ref{theo.main}. If we can find a proper four-coloring of the
triangulation $T(3L,3)$ with a top-row coloring equal to the top-row
coloring of the coloring obtained in the proof of Theorem~\ref{theo.main},
then the above lemma can be used to prove the inductive step. The main 
issue is therefore, to prove the existence of such coloring for
$T(3L,3)$.  

\begin{theorem} \label{theo.asym}
For any triangulation $T(3L,3M)$ with any $L\geq 3$ and $M\geq L$,
there exists a four-coloring $f$ with $\deg(f)\equiv 6\pmod{12}$. 
Consequently, the WSK dynamics for four-colorings of $T(3L,3M)$ is 
non-ergodic.
\end{theorem}

\proof
The proof is by induction on $M$. The base case $M=L\geq 3$ is proven by 
Theorem~\ref{theo.main}. Now suppose that there exist such colorings for
all triangulations $T(3L,3M')$ with $L\le M'\leq M$, and we wish to prove that
such configuration exists also for $M$. 
The main idea is to prove the existence of a proper four-coloring of 
the triangulation $T(3L,3)$ such that its top row coloring coincides with the 
one obtained in the proof of the corresponding case in Theorem~\ref{theo.main}. 

To simplify the notation we will denote by $c_i$ the sequence of colors
in the row $i$ of $T(3L,3)$ and by $c_0$ the coloring of the top row of
$T(3L,3L)$ obtained in the proof of Theorem~\ref{theo.main}. Of course, our
goal is to have $c_0=c_3$.

To describe a sequence of colors, we will use the following notation: 
$[a_1 a_2 \cdots a_s]^t$ will be the sequence of length $st$ in
which $a_1 a_2 \cdots a_s$ is repeated $t$ times. For example,
$12[34]^32 = 123434342$.

Our basic strategy is, as in Theorem~\ref{theo.main}, to explicitly construct
four-colorings of $T(3L,3)$ with $L\geq 3$. The construction of such
a coloring will depend on the value of $L$ modulo 4, and we will split 
the proof in four cases, $L=4k-2, 4k-1, 4k$, or $L=4k+1$, with $k\in\N$.

The case $L=4k-2$ was the easiest one in the proof of Theorem~\ref{theo.main};
however, in this case it is the most elaborate. Thus, we will start the proof
by considering the easiest cases, and delay the most complex one to the end. 

\proofofcase{1}{$L=4k-1$}

Let $t = \lfloor \tfrac{3k-2}{2} \rfloor$.
The top-row coloring obtained from the proof of Case~2 in 
Theorem~\ref{theo.main} can be written as
$$
   c_0 = c_3 = [1423]^t 1231 [3241]^t 3
$$
when $k$ is even. Then we define $c_1$ and $c_2$ as:
\begin{eqnarray*}
 c_2 & = & 3[1423]^t 142 [1324]^t 2 \\
 c_1 & = & 23[1423]^t 14 [2413]^t 4.
\end{eqnarray*}
If $k$ is odd, then we have:
\begin{eqnarray*}
 c_0 = c_3 & = & [1423]^t 14214241 [3241]^t 3 \\
       c_2 & = & 3[1423]^t 1423124 [1324]^t 2 \;=\; 
                 3[1423]^{t+1} 124 [1324]^t 2\\
       c_1 & = & 23[1423]^t 14231 [3241]^t 34 \;=\; 
                 23[1423]^{t+1} 1 [3241]^t 34.
\end{eqnarray*}
It is easy to verify that this gives a proper 4-coloring of $T(3L,3)$.
By Proposition~\ref{prop.T_Lx3}, it has zero degree. This completes the proof 
of this case.

\proofofcase{2}{$L=4k$}

As for the previous case, let $t = \lfloor \tfrac{3k-2}{2} \rfloor$.
The top-row coloring $c_3=c_0$ is obtained from the proof of Case~3 
in Theorem~\ref{theo.main}. When $k$ is even, the sought 4-coloring is
defined as follows: 
\begin{eqnarray*}
c_0 = c_3 &=&  [1423]^t 1431341 [3241]^t 3 \\
      c_2 &=& 3[1423]^t 124132  [4132]^t 4 \;=\; 3[1423]^t 12 [4132]^{t+1} 4\\
      c_1 &=& 4[2314]^t 312413  [2413]^t 2 \;=\; 4[2314]^t 31 [2413]^{t+1} 2
\,.
\end{eqnarray*}
If $k$ is odd, then we have:
\begin{eqnarray*}
 c_0 = c_3& =&   [1423]^t 14234231241 [3241]^t 3 \;=\; 
                 [1423]^{t+1} 4231241 [3241]^t 3  \\
       c_2& =& 3 [1423]^t 1423423132  [4132]^t 4 \;=\; 
               3 [1423]^{t+1} 423132  [4132]^t 4 \\ 
       c_1& =& 4 [2314]^t 2342312413  [2413]^t 2 \;=\; 
               4 [2314]^t 234231 [2413]^{t+1} 2\,.  
\end{eqnarray*}
Again, it is easy to verify that this gives a proper 4-coloring of $T(3L,3)$,
and by Proposition~\ref{prop.T_Lx3}, it has zero degree. 
This completes the proof of this case.

\clearpage 
\proofofcase{3}{$L=4k+1$}

Let $t = \lfloor \tfrac{3k-2}{2} \rfloor$.
The top-row coloring $c_3=c_0$ is obtained from the proof of Case~4 
in Theorem~\ref{theo.main}. When $k$ is even, the sought 4-coloring is
defined as follows: 
\begin{eqnarray*}
 c_0 = c_3 & = &  [1423]^t 1421423421 [3241]^t 3 \\
       c_2 & = & 3[1423]^t 14214213   [2413]^t 42 \\ 
       c_1 & = & 2[3142]^t 314214213  [2413]^t 4 \;=\; 
                 2[3142]^{t+1} 14213  [2413]^t 4\,.
\end{eqnarray*}
If $k$ is odd, then we have:
\begin{eqnarray*}
 c_0 = c_3 & = & [1423]^{t+1} 1231431241 [3241]^t 3 \\
      c_2 & = & [1423]^{t+1} 312312413  [2413]^t 42 \;=\; 
                [1423]^{t+1} 31231      [2413]^{t+1} 42 \\
      c_1 & = & [2314]^{t+1} 2312312413 [2413]^t 2 \;=\; 
                [2314]^{t+1} 2312312413 [2413]^{t+1} 2 \,. 
\end{eqnarray*}
Again, it is easy to verify that this gives a proper 4-coloring of $T(3L,3)$,
and by Proposition~\ref{prop.T_Lx3}, it has zero degree. 
This completes the proof of this case.

\proofofcase{4}{$L=4k-2$}

We cannot use the results of the proof of Theorem~\ref{theo.main},
as the resulting four-coloring for $T(3L,3L)$ is characterized by 
the fact that any row (horizontal, vertical or inclined) is bi-colored.
Thus, we cannot obtain a four-coloring of $T(12k-6,3)$ with a bi-colored
horizontal row.

We first need to obtain a proper four-coloring $f$ of $T(12k-6,12k-6)$ with
$\deg(f)\equiv 6\pmod{12}$, and such as there is a proper four-coloring
of $T(12k-6,3)$ compatible with the coloring of one of the horizontal rows
of $f$. We obtain such coloring $f$ by a constructive proof
similar to those explained in the proof of Theorem~\ref{theo.main}. 
The notation we use is the same as in Theorem~\ref{theo.main}.

Let us consider the triangulation $T=T(12k-6,12k-6)$ with integer $k\ge 2$ 
(the case $k=2$ will illustrate our ideas). Our goal is to obtain 
a four-coloring $f$ of $T$ with degree $\deg(f)\equiv 6 \pmod{12}$. 
The algorithm to obtain such a coloring consists of four steps: 

\medskip
\noindent
{\bf Step 1.}
We start by coloring counter-diagonal D1: we color $1$ the vertices with
$x$-coordinates $1\leq x \leq 6k-3$; the other $6k-3$ vertices are colored $2$.

On D2, we color $3$ those $6k-3$ vertices with $x$-coordinates
$3k\leq x \leq 9k-4$. The other vertices on D2 are colored $4$. The
vertices on D$(12k-6)$ are colored $3$ or $4$ in such a 
way that the resulting coloring is proper (for each vertex, there is a
unique choice).

\begin{figure}[htb]
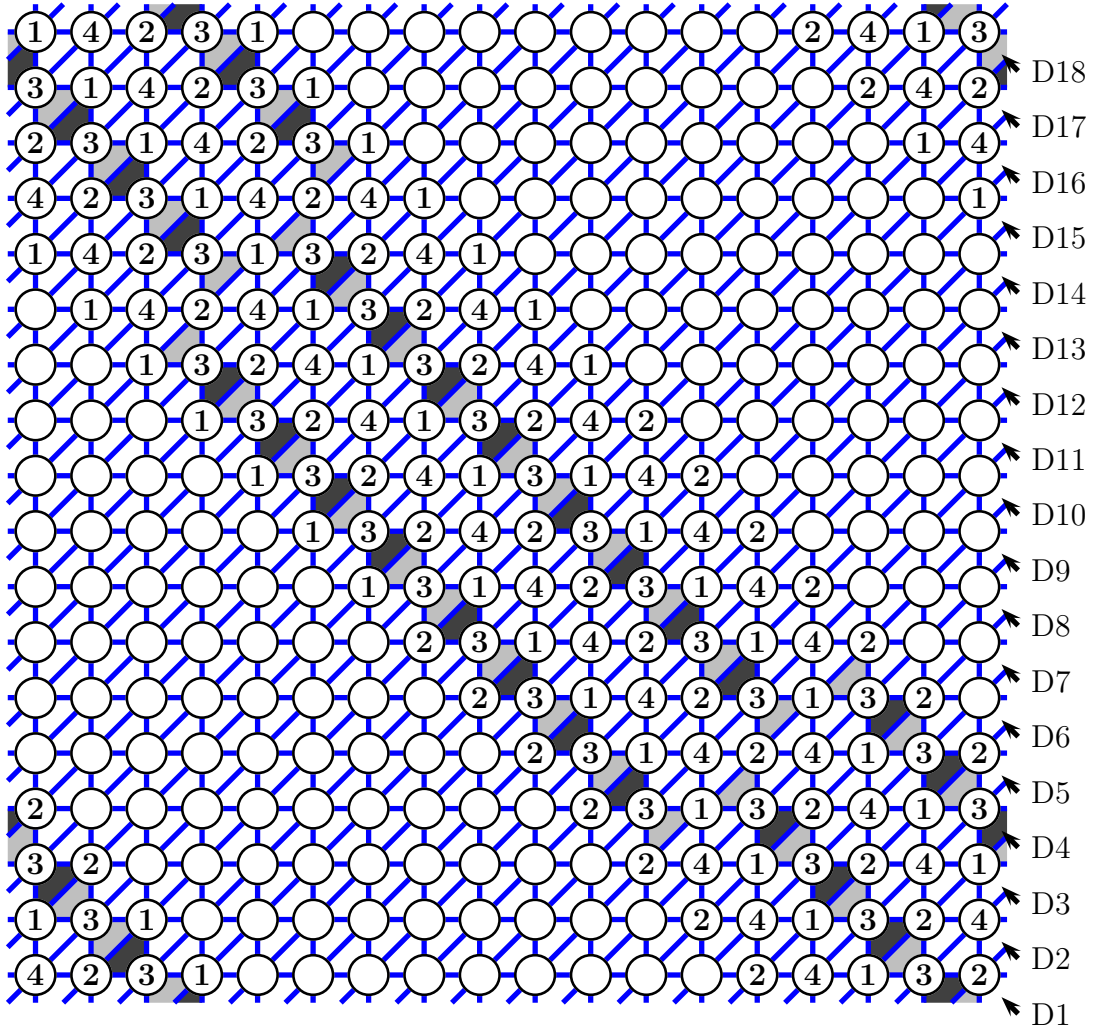

\centering
%
%
\psset{xunit=21pt}
\psset{yunit=21pt}
\psset{labelsep=5pt}
\pspicture(-0.5,-0.9)(18.6,17.5)
\psline*[linewidth=2pt,linecolor=lightgray](0,1)(1,1)(1,2)(0,1)
\psline*[linewidth=2pt,linecolor=darkgray](0,1)(0,2)(1,2)(0,1)
\psline*[linewidth=2pt,linecolor=darkgray](0,15)(1,15)(1,16)(0,15)
\psline*[linewidth=2pt,linecolor=lightgray](0,15)(0,16)(1,16)(0,15)
\psline*[linewidth=2pt,linecolor=darkgray](1,0)(2,0)(2,1)(1,0)
\psline*[linewidth=2pt,linecolor=lightgray](1,0)(1,1)(2,1)(1,0)
\psline*[linewidth=2pt,linecolor=darkgray](1,14)(2,14)(2,15)(1,14)
\psline*[linewidth=2pt,linecolor=lightgray](1,14)(1,15)(2,15)(1,14)
\psline*[linewidth=2pt,linecolor=lightgray](2,11)(3,11)(3,12)(2,11)
\psline*[linewidth=2pt,linecolor=darkgray](2,13)(3,13)(3,14)(2,13)
\psline*[linewidth=2pt,linecolor=lightgray](2,13)(2,14)(3,14)(2,13)
\psline*[linewidth=2pt,linecolor=darkgray](2,17)(3,17)(3,17.5)(2.5,17.5)(2,17)
\psline*[linewidth=2pt,linecolor=darkgray](3,0)(3,-0.5)(2.5,-0.5)(3,0)
\psline*[linewidth=2pt,linecolor=lightgray](2,17)(2,17.5)(2.5,17.5)(2,17)
\psline*[linewidth=2pt,linecolor=lightgray](2,0)(2,-0.5)(2.5,-0.5)(3,0)(2,0)
\psline*[linewidth=2pt,linecolor=lightgray](3,10)(4,10)(4,11)(3,10)
\psline*[linewidth=2pt,linecolor=darkgray](3,10)(3,11)(4,11)(3,10)
\psline*[linewidth=2pt,linecolor=lightgray](3,12)(3,13)(4,13)(3,12)
\psline*[linewidth=2pt,linecolor=darkgray](3,16)(4,16)(4,17)(3,16)
\psline*[linewidth=2pt,linecolor=lightgray](3,16)(3,17)(4,17)(3,16)
\psline*[linewidth=2pt,linecolor=lightgray](4,9)(5,9)(5,10)(4,9)
\psline*[linewidth=2pt,linecolor=darkgray](4,9)(4,10)(5,10)(4,9)
\psline*[linewidth=2pt,linecolor=lightgray](4,13)(5,13)(5,14)(4,13)
\psline*[linewidth=2pt,linecolor=darkgray](4,15)(5,15)(5,16)(4,15)
\psline*[linewidth=2pt,linecolor=lightgray](4,15)(4,16)(5,16)(4,15)
\psline*[linewidth=2pt,linecolor=lightgray](5,8)(6,8)(6,9)(5,8)
\psline*[linewidth=2pt,linecolor=darkgray](5,8)(5,9)(6,9)(5,8)
\psline*[linewidth=2pt,linecolor=lightgray](5,12)(6,12)(6,13)(5,12)
\psline*[linewidth=2pt,linecolor=darkgray](5,12)(5,13)(6,13)(5,12)
\psline*[linewidth=2pt,linecolor=lightgray](5,14)(5,15)(6,15)(5,14)
\psline*[linewidth=2pt,linecolor=lightgray](6,7)(7,7)(7,8)(6,7)
\psline*[linewidth=2pt,linecolor=darkgray](6,7)(6,8)(7,8)(6,7)
\psline*[linewidth=2pt,linecolor=lightgray](6,11)(7,11)(7,12)(6,11)
\psline*[linewidth=2pt,linecolor=darkgray](6,11)(6,12)(7,12)(6,11)
\psline*[linewidth=2pt,linecolor=darkgray](7,6)(8,6)(8,7)(7,6)
\psline*[linewidth=2pt,linecolor=lightgray](7,6)(7,7)(8,7)(7,6)
\psline*[linewidth=2pt,linecolor=lightgray](7,10)(8,10)(8,11)(7,10)
\psline*[linewidth=2pt,linecolor=darkgray](7,10)(7,11)(8,11)(7,10)
\psline*[linewidth=2pt,linecolor=darkgray](8,5)(9,5)(9,6)(8,5)
\psline*[linewidth=2pt,linecolor=lightgray](8,5)(8,6)(9,6)(8,5)
\psline*[linewidth=2pt,linecolor=lightgray](8,9)(9,9)(9,10)(8,9)
\psline*[linewidth=2pt,linecolor=darkgray](8,9)(8,10)(9,10)(8,9)
\psline*[linewidth=2pt,linecolor=darkgray](9,4)(10,4)(10,5)(9,4)
\psline*[linewidth=2pt,linecolor=lightgray](9,4)(9,5)(10,5)(9,4)
\psline*[linewidth=2pt,linecolor=darkgray](9,8)(10,8)(10,9)(9,8)
\psline*[linewidth=2pt,linecolor=lightgray](9,8)(9,9)(10,9)(9,8)
\psline*[linewidth=2pt,linecolor=darkgray](10,3)(11,3)(11,4)(10,3)
\psline*[linewidth=2pt,linecolor=lightgray](10,3)(10,4)(11,4)(10,3)
\psline*[linewidth=2pt,linecolor=darkgray](10,7)(11,7)(11,8)(10,7)
\psline*[linewidth=2pt,linecolor=lightgray](10,7)(10,8)(11,8)(10,7)
\psline*[linewidth=2pt,linecolor=lightgray](11,2)(11,3)(12,3)(11,2)
\psline*[linewidth=2pt,linecolor=darkgray](11,6)(12,6)(12,7)(11,6)
\psline*[linewidth=2pt,linecolor=lightgray](11,6)(11,7)(12,7)(11,6)
\psline*[linewidth=2pt,linecolor=lightgray](12,3)(13,3)(13,4)(12,3)
\psline*[linewidth=2pt,linecolor=darkgray](12,5)(13,5)(13,6)(12,5)
\psline*[linewidth=2pt,linecolor=lightgray](12,5)(12,6)(13,6)(12,5)
\psline*[linewidth=2pt,linecolor=lightgray](13,2)(14,2)(14,3)(13,2)
\psline*[linewidth=2pt,linecolor=darkgray](13,2)(13,3)(14,3)(13,2)
\psline*[linewidth=2pt,linecolor=lightgray](13,4)(13,5)(14,5)(13,4)
\psline*[linewidth=2pt,linecolor=lightgray](14,1)(15,1)(15,2)(14,1)
\psline*[linewidth=2pt,linecolor=darkgray](14,1)(14,2)(15,2)(14,1)
\psline*[linewidth=2pt,linecolor=lightgray](14,5)(15,5)(15,6)(14,5)
\psline*[linewidth=2pt,linecolor=lightgray](15,0)(16,0)(16,1)(15,0)
\psline*[linewidth=2pt,linecolor=darkgray](15,0)(15,1)(16,1)(15,0)
\psline*[linewidth=2pt,linecolor=lightgray](15,4)(16,4)(16,5)(15,4)
\psline*[linewidth=2pt,linecolor=darkgray](15,4)(15,5)(16,5)(15,4)
\psline*[linewidth=2pt,linecolor=lightgray](16,3)(17,3)(17,4)(16,3)
\psline*[linewidth=2pt,linecolor=darkgray](16,3)(16,4)(17,4)(16,3)
\psline*[linewidth=2pt,linecolor=lightgray](16,17)(17,17)(17,17.5)(16.5,17.5)(16,17)
\psline*[linewidth=2pt,linecolor=lightgray](17,0)(17,-0.5)(16.5,-0.5)(17,0)
\psline*[linewidth=2pt,linecolor=darkgray](16,17)(16,17.5)(16.5,17.5)(16,17)
\psline*[linewidth=2pt,linecolor=darkgray](16,0)(16,-0.5)(16.5,-0.5)(17,0)(16,0)
\psline*[linewidth=2pt,linecolor=lightgray](17,2)(17.5,2)(17.5,2.5)(17,2)
\psline*[linewidth=2pt,linecolor=lightgray](0,2)(-0.5,2)(-0.5,2.5)(0,3,)(0,2)
\psline*[linewidth=2pt,linecolor=darkgray](17,2)(17,3)(17.5,3)(17.5,2.5)(17,2)
\psline*[linewidth=2pt,linecolor=darkgray](0,3)(-0.5,3)(-0.5,2.5)(0,3)
\psline*[linewidth=2pt,linecolor=darkgray](17,16)(17.5,16)(17.5,16.5)(17,16)
\psline*[linewidth=2pt,linecolor=darkgray](0,16)(-0.5,16)(-0.5,16.5)(0,17,)(0,16)
\psline*[linewidth=2pt,linecolor=lightgray](17,16)(17,17)(17.5,17)(17.5,16.5)(17,16)
\psline*[linewidth=2pt,linecolor=lightgray](0,17)(-0.5,17)(-0.5,16.5)(0,17)
\psline[linewidth=2pt,linecolor=blue](-0.5,0)(17.5,0)
\psline[linewidth=2pt,linecolor=blue](-0.5,1)(17.5,1)
\psline[linewidth=2pt,linecolor=blue](-0.5,2)(17.5,2)
\psline[linewidth=2pt,linecolor=blue](-0.5,3)(17.5,3)
\psline[linewidth=2pt,linecolor=blue](-0.5,4)(17.5,4)
\psline[linewidth=2pt,linecolor=blue](-0.5,5)(17.5,5)
\psline[linewidth=2pt,linecolor=blue](-0.5,6)(17.5,6)
\psline[linewidth=2pt,linecolor=blue](-0.5,7)(17.5,7)
\psline[linewidth=2pt,linecolor=blue](-0.5,8)(17.5,8)
\psline[linewidth=2pt,linecolor=blue](-0.5,9)(17.5,9)
\psline[linewidth=2pt,linecolor=blue](-0.5,10)(17.5,10)
\psline[linewidth=2pt,linecolor=blue](-0.5,11)(17.5,11)
\psline[linewidth=2pt,linecolor=blue](-0.5,12)(17.5,12)
\psline[linewidth=2pt,linecolor=blue](-0.5,13)(17.5,13)
\psline[linewidth=2pt,linecolor=blue](-0.5,14)(17.5,14)
\psline[linewidth=2pt,linecolor=blue](-0.5,15)(17.5,15)
\psline[linewidth=2pt,linecolor=blue](-0.5,16)(17.5,16)
\psline[linewidth=2pt,linecolor=blue](-0.5,17)(17.5,17)
\psline[linewidth=2pt,linecolor=blue](0,-0.5)(0,17.5)
\psline[linewidth=2pt,linecolor=blue](1,-0.5)(1,17.5)
\psline[linewidth=2pt,linecolor=blue](2,-0.5)(2,17.5)
\psline[linewidth=2pt,linecolor=blue](3,-0.5)(3,17.5)
\psline[linewidth=2pt,linecolor=blue](4,-0.5)(4,17.5)
\psline[linewidth=2pt,linecolor=blue](5,-0.5)(5,17.5)
\psline[linewidth=2pt,linecolor=blue](6,-0.5)(6,17.5)
\psline[linewidth=2pt,linecolor=blue](7,-0.5)(7,17.5)
\psline[linewidth=2pt,linecolor=blue](8,-0.5)(8,17.5)
\psline[linewidth=2pt,linecolor=blue](9,-0.5)(9,17.5)
\psline[linewidth=2pt,linecolor=blue](10,-0.5)(10,17.5)
\psline[linewidth=2pt,linecolor=blue](11,-0.5)(11,17.5)
\psline[linewidth=2pt,linecolor=blue](12,-0.5)(12,17.5)
\psline[linewidth=2pt,linecolor=blue](13,-0.5)(13,17.5)
\psline[linewidth=2pt,linecolor=blue](14,-0.5)(14,17.5)
\psline[linewidth=2pt,linecolor=blue](15,-0.5)(15,17.5)
\psline[linewidth=2pt,linecolor=blue](16,-0.5)(16,17.5)
\psline[linewidth=2pt,linecolor=blue](17,-0.5)(17,17.5)
\psline[linewidth=2pt,linecolor=blue](-0.5,-0.5)(17.5,17.5)
\psline[linewidth=2pt,linecolor=blue](0.5,-0.5)(17.5,16.5)
\psline[linewidth=2pt,linecolor=blue](1.5,-0.5)(17.5,15.5)
\psline[linewidth=2pt,linecolor=blue](2.5,-0.5)(17.5,14.5)
\psline[linewidth=2pt,linecolor=blue](3.5,-0.5)(17.5,13.5)
\psline[linewidth=2pt,linecolor=blue](4.5,-0.5)(17.5,12.5)
\psline[linewidth=2pt,linecolor=blue](5.5,-0.5)(17.5,11.5)
\psline[linewidth=2pt,linecolor=blue](6.5,-0.5)(17.5,10.5)
\psline[linewidth=2pt,linecolor=blue](7.5,-0.5)(17.5,9.5)
\psline[linewidth=2pt,linecolor=blue](8.5,-0.5)(17.5,8.5)
\psline[linewidth=2pt,linecolor=blue](9.5,-0.5)(17.5,7.5)
\psline[linewidth=2pt,linecolor=blue](10.5,-0.5)(17.5,6.5)
\psline[linewidth=2pt,linecolor=blue](11.5,-0.5)(17.5,5.5)
\psline[linewidth=2pt,linecolor=blue](12.5,-0.5)(17.5,4.5)
\psline[linewidth=2pt,linecolor=blue](13.5,-0.5)(17.5,3.5)
\psline[linewidth=2pt,linecolor=blue](14.5,-0.5)(17.5,2.5)
\psline[linewidth=2pt,linecolor=blue](15.5,-0.5)(17.5,1.5)
\psline[linewidth=2pt,linecolor=blue](16.5,-0.5)(17.5,0.5)
\psline[linewidth=2pt,linecolor=blue](-0.5,0.5)(16.5,17.5)
\psline[linewidth=2pt,linecolor=blue](-0.5,1.5)(15.5,17.5)
\psline[linewidth=2pt,linecolor=blue](-0.5,2.5)(14.5,17.5)
\psline[linewidth=2pt,linecolor=blue](-0.5,3.5)(13.5,17.5)
\psline[linewidth=2pt,linecolor=blue](-0.5,4.5)(12.5,17.5)
\psline[linewidth=2pt,linecolor=blue](-0.5,5.5)(11.5,17.5)
\psline[linewidth=2pt,linecolor=blue](-0.5,6.5)(10.5,17.5)
\psline[linewidth=2pt,linecolor=blue](-0.5,7.5)(9.5,17.5)
\psline[linewidth=2pt,linecolor=blue](-0.5,8.5)(8.5,17.5)
\psline[linewidth=2pt,linecolor=blue](-0.5,9.5)(7.5,17.5)
\psline[linewidth=2pt,linecolor=blue](-0.5,10.5)(6.5,17.5)
\psline[linewidth=2pt,linecolor=blue](-0.5,11.5)(5.5,17.5)
\psline[linewidth=2pt,linecolor=blue](-0.5,12.5)(4.5,17.5)
\psline[linewidth=2pt,linecolor=blue](-0.5,13.5)(3.5,17.5)
\psline[linewidth=2pt,linecolor=blue](-0.5,14.5)(2.5,17.5)
\psline[linewidth=2pt,linecolor=blue](-0.5,15.5)(1.5,17.5)
\psline[linewidth=2pt,linecolor=blue](-0.5,16.5)(0.5,17.5)
\multirput{0}(0,0)(0,1){18}{%
  \multirput{0}(0,0)(1,0){18}{%
     \pscircle*[linecolor=white]{8pt}
     \pscircle[linewidth=1pt,linecolor=black] {8pt}
   }
}
\rput{0}(0,0){{\bf 4}}
\rput{0}(0,1){{\bf 1}}
\rput{0}(0,2){{\bf 3}}
\rput{0}(0,3){{\bf 2}}
\rput{0}(0,4){{\bf }}
\rput{0}(0,5){{\bf }}
\rput{0}(0,6){{\bf }}
\rput{0}(0,7){{\bf }}
\rput{0}(0,8){{\bf }}
\rput{0}(0,9){{\bf }}
\rput{0}(0,10){{\bf }}
\rput{0}(0,11){{\bf }}
\rput{0}(0,12){{\bf }}
\rput{0}(0,13){{\bf 1}}
\rput{0}(0,14){{\bf 4}}
\rput{0}(0,15){{\bf 2}}
\rput{0}(0,16){{\bf 3}}
\rput{0}(0,17){{\bf 1}}
\rput{0}(1,0){{\bf 2}}
\rput{0}(1,1){{\bf 3}}
\rput{0}(1,2){{\bf 2}}
\rput{0}(1,3){{\bf }}
\rput{0}(1,4){{\bf }}
\rput{0}(1,5){{\bf }}
\rput{0}(1,6){{\bf }}
\rput{0}(1,7){{\bf }}
\rput{0}(1,8){{\bf }}
\rput{0}(1,9){{\bf }}
\rput{0}(1,10){{\bf }}
\rput{0}(1,11){{\bf }}
\rput{0}(1,12){{\bf 1}}
\rput{0}(1,13){{\bf 4}}
\rput{0}(1,14){{\bf 2}}
\rput{0}(1,15){{\bf 3}}
\rput{0}(1,16){{\bf 1}}
\rput{0}(1,17){{\bf 4}}
\rput{0}(2,0){{\bf 3}}
\rput{0}(2,1){{\bf 1}}
\rput{0}(2,2){{\bf }}
\rput{0}(2,3){{\bf }}
\rput{0}(2,4){{\bf }}
\rput{0}(2,5){{\bf }}
\rput{0}(2,6){{\bf }}
\rput{0}(2,7){{\bf }}
\rput{0}(2,8){{\bf }}
\rput{0}(2,9){{\bf }}
\rput{0}(2,10){{\bf }}
\rput{0}(2,11){{\bf 1}}
\rput{0}(2,12){{\bf 4}}
\rput{0}(2,13){{\bf 2}}
\rput{0}(2,14){{\bf 3}}
\rput{0}(2,15){{\bf 1}}
\rput{0}(2,16){{\bf 4}}
\rput{0}(2,17){{\bf 2}}
\rput{0}(3,0){{\bf 1}}
\rput{0}(3,1){{\bf }}
\rput{0}(3,2){{\bf }}
\rput{0}(3,3){{\bf }}
\rput{0}(3,4){{\bf }}
\rput{0}(3,5){{\bf }}
\rput{0}(3,6){{\bf }}
\rput{0}(3,7){{\bf }}
\rput{0}(3,8){{\bf }}
\rput{0}(3,9){{\bf }}
\rput{0}(3,10){{\bf 1}}
\rput{0}(3,11){{\bf 3}}
\rput{0}(3,12){{\bf 2}}
\rput{0}(3,13){{\bf 3}}
\rput{0}(3,14){{\bf 1}}
\rput{0}(3,15){{\bf 4}}
\rput{0}(3,16){{\bf 2}}
\rput{0}(3,17){{\bf 3}}
\rput{0}(4,0){{\bf }}
\rput{0}(4,1){{\bf }}
\rput{0}(4,2){{\bf }}
\rput{0}(4,3){{\bf }}
\rput{0}(4,4){{\bf }}
\rput{0}(4,5){{\bf }}
\rput{0}(4,6){{\bf }}
\rput{0}(4,7){{\bf }}
\rput{0}(4,8){{\bf }}
\rput{0}(4,9){{\bf 1}}
\rput{0}(4,10){{\bf 3}}
\rput{0}(4,11){{\bf 2}}
\rput{0}(4,12){{\bf 4}}
\rput{0}(4,13){{\bf 1}}
\rput{0}(4,14){{\bf 4}}
\rput{0}(4,15){{\bf 2}}
\rput{0}(4,16){{\bf 3}}
\rput{0}(4,17){{\bf 1}}
\rput{0}(5,0){{\bf }}
\rput{0}(5,1){{\bf }}
\rput{0}(5,2){{\bf }}
\rput{0}(5,3){{\bf }}
\rput{0}(5,4){{\bf }}
\rput{0}(5,5){{\bf }}
\rput{0}(5,6){{\bf }}
\rput{0}(5,7){{\bf }}
\rput{0}(5,8){{\bf 1}}
\rput{0}(5,9){{\bf 3}}
\rput{0}(5,10){{\bf 2}}
\rput{0}(5,11){{\bf 4}}
\rput{0}(5,12){{\bf 1}}
\rput{0}(5,13){{\bf 3}}
\rput{0}(5,14){{\bf 2}}
\rput{0}(5,15){{\bf 3}}
\rput{0}(5,16){{\bf 1}}
\rput{0}(5,17){{\bf }}
\rput{0}(6,0){{\bf }}
\rput{0}(6,1){{\bf }}
\rput{0}(6,2){{\bf }}
\rput{0}(6,3){{\bf }}
\rput{0}(6,4){{\bf }}
\rput{0}(6,5){{\bf }}
\rput{0}(6,6){{\bf }}
\rput{0}(6,7){{\bf 1}}
\rput{0}(6,8){{\bf 3}}
\rput{0}(6,9){{\bf 2}}
\rput{0}(6,10){{\bf 4}}
\rput{0}(6,11){{\bf 1}}
\rput{0}(6,12){{\bf 3}}
\rput{0}(6,13){{\bf 2}}
\rput{0}(6,14){{\bf 4}}
\rput{0}(6,15){{\bf 1}}
\rput{0}(6,16){{\bf }}
\rput{0}(6,17){{\bf }}
\rput{0}(7,0){{\bf }}
\rput{0}(7,1){{\bf }}
\rput{0}(7,2){{\bf }}
\rput{0}(7,3){{\bf }}
\rput{0}(7,4){{\bf }}
\rput{0}(7,5){{\bf }}
\rput{0}(7,6){{\bf 2}}
\rput{0}(7,7){{\bf 3}}
\rput{0}(7,8){{\bf 2}}
\rput{0}(7,9){{\bf 4}}
\rput{0}(7,10){{\bf 1}}
\rput{0}(7,11){{\bf 3}}
\rput{0}(7,12){{\bf 2}}
\rput{0}(7,13){{\bf 4}}
\rput{0}(7,14){{\bf 1}}
\rput{0}(7,15){{\bf }}
\rput{0}(7,16){{\bf }}
\rput{0}(7,17){{\bf }}
\rput{0}(8,0){{\bf }}
\rput{0}(8,1){{\bf }}
\rput{0}(8,2){{\bf }}
\rput{0}(8,3){{\bf }}
\rput{0}(8,4){{\bf }}
\rput{0}(8,5){{\bf 2}}
\rput{0}(8,6){{\bf 3}}
\rput{0}(8,7){{\bf 1}}
\rput{0}(8,8){{\bf 4}}
\rput{0}(8,9){{\bf 1}}
\rput{0}(8,10){{\bf 3}}
\rput{0}(8,11){{\bf 2}}
\rput{0}(8,12){{\bf 4}}
\rput{0}(8,13){{\bf 1}}
\rput{0}(8,14){{\bf }}
\rput{0}(8,15){{\bf }}
\rput{0}(8,16){{\bf }}
\rput{0}(8,17){{\bf }}
\rput{0}(9,0){{\bf }}
\rput{0}(9,1){{\bf }}
\rput{0}(9,2){{\bf }}
\rput{0}(9,3){{\bf }}
\rput{0}(9,4){{\bf 2}}
\rput{0}(9,5){{\bf 3}}
\rput{0}(9,6){{\bf 1}}
\rput{0}(9,7){{\bf 4}}
\rput{0}(9,8){{\bf 2}}
\rput{0}(9,9){{\bf 3}}
\rput{0}(9,10){{\bf 2}}
\rput{0}(9,11){{\bf 4}}
\rput{0}(9,12){{\bf 1}}
\rput{0}(9,13){{\bf }}
\rput{0}(9,14){{\bf }}
\rput{0}(9,15){{\bf }}
\rput{0}(9,16){{\bf }}
\rput{0}(9,17){{\bf }}
\rput{0}(10,0){{\bf }}
\rput{0}(10,1){{\bf }}
\rput{0}(10,2){{\bf }}
\rput{0}(10,3){{\bf 2}}
\rput{0}(10,4){{\bf 3}}
\rput{0}(10,5){{\bf 1}}
\rput{0}(10,6){{\bf 4}}
\rput{0}(10,7){{\bf 2}}
\rput{0}(10,8){{\bf 3}}
\rput{0}(10,9){{\bf 1}}
\rput{0}(10,10){{\bf 4}}
\rput{0}(10,11){{\bf 1}}
\rput{0}(10,12){{\bf }}
\rput{0}(10,13){{\bf }}
\rput{0}(10,14){{\bf }}
\rput{0}(10,15){{\bf }}
\rput{0}(10,16){{\bf }}
\rput{0}(10,17){{\bf }}
\rput{0}(11,0){{\bf }}
\rput{0}(11,1){{\bf }}
\rput{0}(11,2){{\bf 2}}
\rput{0}(11,3){{\bf 3}}
\rput{0}(11,4){{\bf 1}}
\rput{0}(11,5){{\bf 4}}
\rput{0}(11,6){{\bf 2}}
\rput{0}(11,7){{\bf 3}}
\rput{0}(11,8){{\bf 1}}
\rput{0}(11,9){{\bf 4}}
\rput{0}(11,10){{\bf 2}}
\rput{0}(11,11){{\bf }}
\rput{0}(11,12){{\bf }}
\rput{0}(11,13){{\bf }}
\rput{0}(11,14){{\bf }}
\rput{0}(11,15){{\bf }}
\rput{0}(11,16){{\bf }}
\rput{0}(11,17){{\bf }}
\rput{0}(12,0){{\bf }}
\rput{0}(12,1){{\bf 2}}
\rput{0}(12,2){{\bf 4}}
\rput{0}(12,3){{\bf 1}}
\rput{0}(12,4){{\bf 4}}
\rput{0}(12,5){{\bf 2}}
\rput{0}(12,6){{\bf 3}}
\rput{0}(12,7){{\bf 1}}
\rput{0}(12,8){{\bf 4}}
\rput{0}(12,9){{\bf 2}}
\rput{0}(12,10){{\bf }}
\rput{0}(12,11){{\bf }}
\rput{0}(12,12){{\bf }}
\rput{0}(12,13){{\bf }}
\rput{0}(12,14){{\bf }}
\rput{0}(12,15){{\bf }}
\rput{0}(12,16){{\bf }}
\rput{0}(12,17){{\bf }}
\rput{0}(13,0){{\bf 2}}
\rput{0}(13,1){{\bf 4}}
\rput{0}(13,2){{\bf 1}}
\rput{0}(13,3){{\bf 3}}
\rput{0}(13,4){{\bf 2}}
\rput{0}(13,5){{\bf 3}}
\rput{0}(13,6){{\bf 1}}
\rput{0}(13,7){{\bf 4}}
\rput{0}(13,8){{\bf 2}}
\rput{0}(13,9){{\bf }}
\rput{0}(13,10){{\bf }}
\rput{0}(13,11){{\bf }}
\rput{0}(13,12){{\bf }}
\rput{0}(13,13){{\bf }}
\rput{0}(13,14){{\bf }}
\rput{0}(13,15){{\bf }}
\rput{0}(13,16){{\bf }}
\rput{0}(13,17){{\bf }}
\rput{0}(14,0){{\bf 4}}
\rput{0}(14,1){{\bf 1}}
\rput{0}(14,2){{\bf 3}}
\rput{0}(14,3){{\bf 2}}
\rput{0}(14,4){{\bf 4}}
\rput{0}(14,5){{\bf 1}}
\rput{0}(14,6){{\bf 4}}
\rput{0}(14,7){{\bf 2}}
\rput{0}(14,8){{\bf }}
\rput{0}(14,9){{\bf }}
\rput{0}(14,10){{\bf }}
\rput{0}(14,11){{\bf }}
\rput{0}(14,12){{\bf }}
\rput{0}(14,13){{\bf }}
\rput{0}(14,14){{\bf }}
\rput{0}(14,15){{\bf }}
\rput{0}(14,16){{\bf }}
\rput{0}(14,17){{\bf 2}}
\rput{0}(15,0){{\bf 1}}
\rput{0}(15,1){{\bf 3}}
\rput{0}(15,2){{\bf 2}}
\rput{0}(15,3){{\bf 4}}
\rput{0}(15,4){{\bf 1}}
\rput{0}(15,5){{\bf 3}}
\rput{0}(15,6){{\bf 2}}
\rput{0}(15,7){{\bf }}
\rput{0}(15,8){{\bf }}
\rput{0}(15,9){{\bf }}
\rput{0}(15,10){{\bf }}
\rput{0}(15,11){{\bf }}
\rput{0}(15,12){{\bf }}
\rput{0}(15,13){{\bf }}
\rput{0}(15,14){{\bf }}
\rput{0}(15,15){{\bf }}
\rput{0}(15,16){{\bf 2}}
\rput{0}(15,17){{\bf 4}}
\rput{0}(16,0){{\bf 3}}
\rput{0}(16,1){{\bf 2}}
\rput{0}(16,2){{\bf 4}}
\rput{0}(16,3){{\bf 1}}
\rput{0}(16,4){{\bf 3}}
\rput{0}(16,5){{\bf 2}}
\rput{0}(16,6){{\bf }}
\rput{0}(16,7){{\bf }}
\rput{0}(16,8){{\bf }}
\rput{0}(16,9){{\bf }}
\rput{0}(16,10){{\bf }}
\rput{0}(16,11){{\bf }}
\rput{0}(16,12){{\bf }}
\rput{0}(16,13){{\bf }}
\rput{0}(16,14){{\bf }}
\rput{0}(16,15){{\bf 1}}
\rput{0}(16,16){{\bf 4}}
\rput{0}(16,17){{\bf 1}}
\rput{0}(17,0){{\bf 2}}
\rput{0}(17,1){{\bf 4}}
\rput{0}(17,2){{\bf 1}}
\rput{0}(17,3){{\bf 3}}
\rput{0}(17,4){{\bf 2}}
\rput{0}(17,5){{\bf }}
\rput{0}(17,6){{\bf }}
\rput{0}(17,7){{\bf }}
\rput{0}(17,8){{\bf }}
\rput{0}(17,9){{\bf }}
\rput{0}(17,10){{\bf }}
\rput{0}(17,11){{\bf }}
\rput{0}(17,12){{\bf }}
\rput{0}(17,13){{\bf }}
\rput{0}(17,14){{\bf 1}}
\rput{0}(17,15){{\bf 4}}
\rput{0}(17,16){{\bf 2}}
\rput{0}(17,17){{\bf 3}}
\multirput{0}(17.5,-0.5)(0,1){18}{\psline[linewidth=2pt,linecolor=black]{->}(0.2,-0.2)(-0.1,0.1)}
\uput[0](17.7,-0.7){D1}
\uput[0](17.7,0.3){D2}
\uput[0](17.7,1.3){D3}
\uput[0](17.7,2.3){D4}
\uput[0](17.7,3.3){D5}
\uput[0](17.7,4.3){D6}
\uput[0](17.7,5.3){D7}
\uput[0](17.7,6.3){D8}
\uput[0](17.7,7.3){D9}
\uput[0](17.7,8.3){D10}
\uput[0](17.7,9.3){D11}
\uput[0](17.7,10.3){D12}
\uput[0](17.7,11.3){D13}
\uput[0](17.7,12.3){D14}
\uput[0](17.7,13.3){D15}
\uput[0](17.7,14.3){D16}
\uput[0](17.7,15.3){D17}
\uput[0](17.7,16.3){D18}
\endpspicture
\caption{ \label{prop.12k-6.fig1}
The 4-coloring of $T(18,18)$ after Step~1 in the case $L=4k-2$.
}
\end{figure}

On D3 and D$(12k-7)$, we color all vertices $1$ or $2$; on D4 and 
D$(12k-8)$, we color all vertices $3$ and $4$, and finally, on D5 and
D$(12k-9)$, we color all vertices $1$ and $2$. In every case, there is a 
unique color choice for each vertex. The resulting coloring is depicted on 
Figure~\ref{prop.12k-6.fig1}. The partial degree of $f$ is $\deg f|_R = 8$.

\medskip
\noindent
{\bf Step 2.}
For $k>2$, we find that there are $12k-15$ counter-diagonals to be colored
and we need to sequentially color all of them but nine. (Notice that this is
why this algorithm does not work for $k=1$.) This can be 
achieved by performing the following procedure: suppose that we
have already colored counter-diagonals D$j$ and D$(12k-j-4)$ ($j\geq 5$) 
using colors $1$ and $2$. Then, we color D$(j+1)$ and D$(12k-j-5)$ using
colors $3$ and $4$, and D$(j+2)$ and D$(12k-j-6)$ using colors $1$ and $2$. 
As in Step~1, for each vertex there is a unique choice.

%
%
%
%

This procedure is repeated $3(k-2)$ times, so we add $12(k-2)$
counter-diagonals, and there are only nine counter-diagonals not yet
colored. Indeed, the last colored counter-diagonals D$(6k-7)$ and
D$(6k+3)$ have colors $1$ and $2$, the same as it was at the end of Step~1.

Each of these $3(k-2)$ steps adds $4$ to the degree of the coloring.
Thus, the partial degree of $f$ is $\deg f|_R = 8 + 12(k-2)$.

\medskip
\noindent
{\bf Step 3.}
On D$(6k-6)$, the vertices $(3k-3,3k-3)$ and $(9k-6,9k-6)$ only admit  
a single color (which is $3$ for one of them, and $4$ for the other one).
The rest of the vertices on D$(6k-6)$ are colored $1$ and $2$ (again, there
is a unique choice for each vertex).

On D$(6k+2)$, there are two vertices $(3k+1,3k+1)$ and $(9k-2,9k-2)$ 
admitting a single color (again $3$ or $4$). The other vertices on D$(6k+2)$ 
are colored $1$ or $2$ (again, the choice for each vertex is unique).

On $D(6k-5)$ there are four vertices which admit a single color $\in\{3,4\}$:
vertices $(3k-2,3k-3)$ and $(3k-3,3k-2)$ should be colored $c_1$, while
$(9k-5,9k-6)$ and $(9k-6,9k-5)$ should be colored $c_2\neq c_1$.
The other vertices satisfying $3k-1\leq x \leq 9k-4$ are colored $c_2$,
and the rest of the vertices are colored $c_1$.  

Finally, on D$(6k+1)$, we also find another four vertices admitting a single
color chosen from the set $\{3,4\}$:
vertices $(3k+1,3k)$ and $(3k,3k+1)$ should be colored $c_1$, while
$(9k-2,9k-3)$ and $(9k-3,9k-2)$ should be colored $c_2\neq c_1$.
The other vertices satisfying $3k+2\leq x \leq 9k-4$ are colored $c_2$,
and the rest of the vertices are colored $c_1$.  

The contribution to the partial degree of the
new triangles is $-4$; the partial degree of $f$ is given by
$\deg f|_R = 4 + 12(k-2)$.

%
%
\begin{figure}[htb]
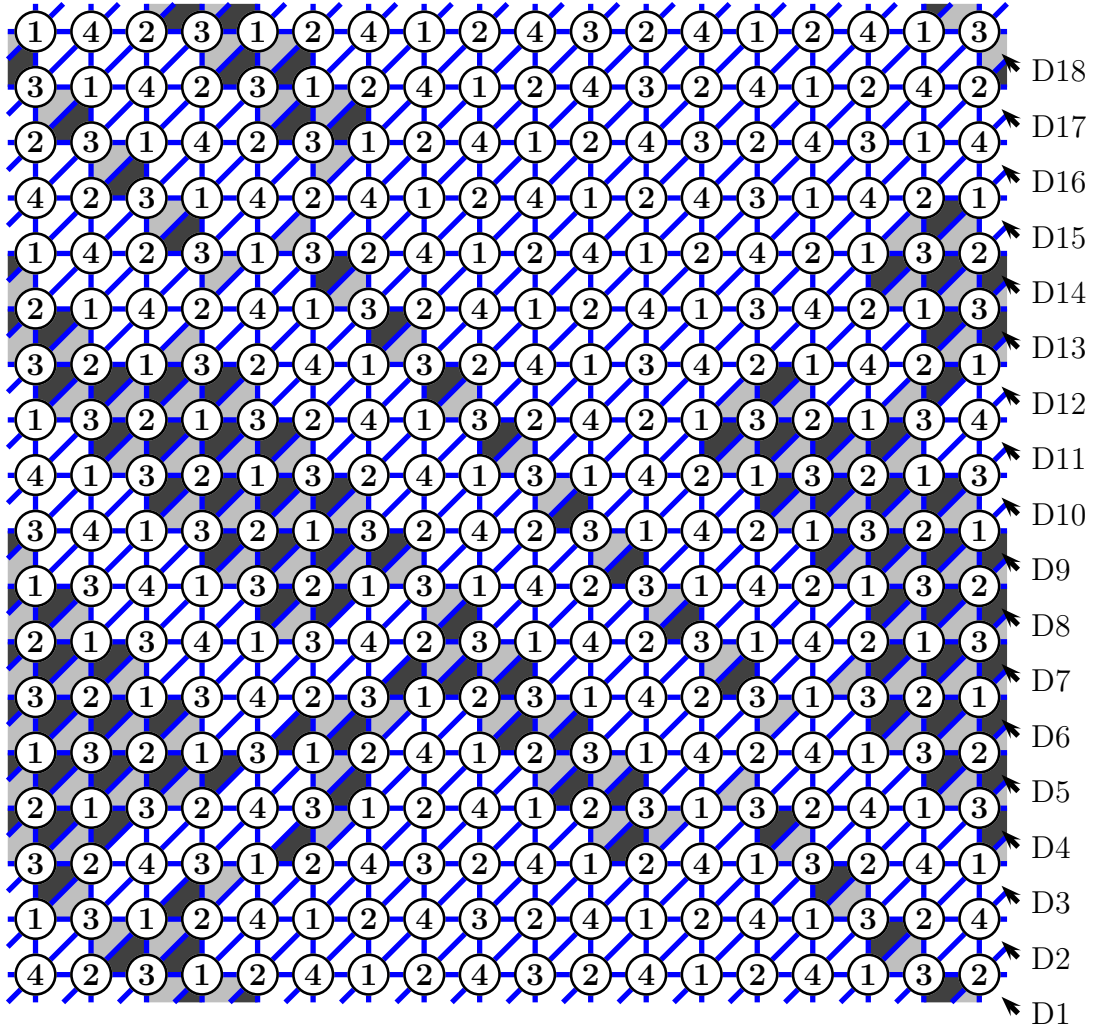

\centering
%
%
\psset{xunit=21pt}
\psset{yunit=21pt}
\psset{labelsep=5pt}
\pspicture(-0.5,-0.9)(18.6,17.5)
\psline*[linewidth=2pt,linecolor=lightgray](0,1)(1,1)(1,2)(0,1)
\psline*[linewidth=2pt,linecolor=darkgray](0,1)(0,2)(1,2)(0,1)
\psline*[linewidth=2pt,linecolor=lightgray](0,2)(1,2)(1,3)(0,2)
\psline*[linewidth=2pt,linecolor=darkgray](0,2)(0,3)(1,3)(0,2)
\psline*[linewidth=2pt,linecolor=lightgray](0,3)(1,3)(1,4)(0,3)
\psline*[linewidth=2pt,linecolor=darkgray](0,3)(0,4)(1,4)(0,3)
\psline*[linewidth=2pt,linecolor=lightgray](0,4)(1,4)(1,5)(0,4)
\psline*[linewidth=2pt,linecolor=darkgray](0,4)(0,5)(1,5)(0,4)
\psline*[linewidth=2pt,linecolor=lightgray](0,5)(1,5)(1,6)(0,5)
\psline*[linewidth=2pt,linecolor=darkgray](0,5)(0,6)(1,6)(0,5)
\psline*[linewidth=2pt,linecolor=lightgray](0,6)(1,6)(1,7)(0,6)
\psline*[linewidth=2pt,linecolor=darkgray](0,6)(0,7)(1,7)(0,6)
\psline*[linewidth=2pt,linecolor=lightgray](0,10)(1,10)(1,11)(0,10)
\psline*[linewidth=2pt,linecolor=darkgray](0,10)(0,11)(1,11)(0,10)
\psline*[linewidth=2pt,linecolor=lightgray](0,11)(1,11)(1,12)(0,11)
\psline*[linewidth=2pt,linecolor=darkgray](0,11)(0,12)(1,12)(0,11)
\psline*[linewidth=2pt,linecolor=darkgray](0,15)(1,15)(1,16)(0,15)
\psline*[linewidth=2pt,linecolor=lightgray](0,15)(0,16)(1,16)(0,15)
\psline*[linewidth=2pt,linecolor=darkgray](1,0)(2,0)(2,1)(1,0)
\psline*[linewidth=2pt,linecolor=lightgray](1,0)(1,1)(2,1)(1,0)
\psline*[linewidth=2pt,linecolor=darkgray](1,2)(1,3)(2,3)(1,2)
\psline*[linewidth=2pt,linecolor=lightgray](1,3)(2,3)(2,4)(1,3)
\psline*[linewidth=2pt,linecolor=darkgray](1,3)(1,4)(2,4)(1,3)
\psline*[linewidth=2pt,linecolor=lightgray](1,4)(2,4)(2,5)(1,4)
\psline*[linewidth=2pt,linecolor=darkgray](1,4)(1,5)(2,5)(1,4)
\psline*[linewidth=2pt,linecolor=lightgray](1,5)(2,5)(2,6)(1,5)
\psline*[linewidth=2pt,linecolor=darkgray](1,5)(1,6)(2,6)(1,5)
\psline*[linewidth=2pt,linecolor=lightgray](1,9)(2,9)(2,10)(1,9)
\psline*[linewidth=2pt,linecolor=darkgray](1,9)(1,10)(2,10)(1,9)
\psline*[linewidth=2pt,linecolor=lightgray](1,10)(2,10)(2,11)(1,10)
\psline*[linewidth=2pt,linecolor=darkgray](1,10)(1,11)(2,11)(1,10)
\psline*[linewidth=2pt,linecolor=darkgray](1,14)(2,14)(2,15)(1,14)
\psline*[linewidth=2pt,linecolor=lightgray](1,14)(1,15)(2,15)(1,14)
\psline*[linewidth=2pt,linecolor=darkgray](2,0)(3,0)(3,1)(2,0)
\psline*[linewidth=2pt,linecolor=lightgray](2,0)(2,1)(3,1)(2,0)
\psline*[linewidth=2pt,linecolor=darkgray](2,1)(3,1)(3,2)(2,1)
\psline*[linewidth=2pt,linecolor=lightgray](2,3)(3,3)(3,4)(2,3)
\psline*[linewidth=2pt,linecolor=darkgray](2,3)(2,4)(3,4)(2,3)
\psline*[linewidth=2pt,linecolor=lightgray](2,4)(3,4)(3,5)(2,4)
\psline*[linewidth=2pt,linecolor=darkgray](2,4)(2,5)(3,5)(2,4)
\psline*[linewidth=2pt,linecolor=lightgray](2,8)(3,8)(3,9)(2,8)
\psline*[linewidth=2pt,linecolor=darkgray](2,8)(2,9)(3,9)(2,8)
\psline*[linewidth=2pt,linecolor=lightgray](2,9)(3,9)(3,10)(2,9)
\psline*[linewidth=2pt,linecolor=darkgray](2,9)(2,10)(3,10)(2,9)
\psline*[linewidth=2pt,linecolor=lightgray](2,10)(3,10)(3,11)(2,10)
\psline*[linewidth=2pt,linecolor=darkgray](2,10)(2,11)(3,11)(2,10)
\psline*[linewidth=2pt,linecolor=lightgray](2,11)(3,11)(3,12)(2,11)
\psline*[linewidth=2pt,linecolor=darkgray](2,13)(3,13)(3,14)(2,13)
\psline*[linewidth=2pt,linecolor=lightgray](2,13)(2,14)(3,14)(2,13)
\psline*[linewidth=2pt,linecolor=darkgray](2,17)(3,17)(3,17.5)(2.5,17.5)(2,17)
\psline*[linewidth=2pt,linecolor=darkgray](3,0)(3,-0.5)(2.5,-0.5)(3,0)
\psline*[linewidth=2pt,linecolor=lightgray](2,17)(2,17.5)(2.5,17.5)(2,17)
\psline*[linewidth=2pt,linecolor=lightgray](2,0)(2,-0.5)(2.5,-0.5)(3,0)(2,0)
\psline*[linewidth=2pt,linecolor=lightgray](3,1)(3,2)(4,2)(3,1)
\psline*[linewidth=2pt,linecolor=darkgray](3,3)(3,4)(4,4)(3,3)
\psline*[linewidth=2pt,linecolor=lightgray](3,7)(4,7)(4,8)(3,7)
\psline*[linewidth=2pt,linecolor=darkgray](3,7)(3,8)(4,8)(3,7)
\psline*[linewidth=2pt,linecolor=lightgray](3,8)(4,8)(4,9)(3,8)
\psline*[linewidth=2pt,linecolor=darkgray](3,8)(3,9)(4,9)(3,8)
\psline*[linewidth=2pt,linecolor=lightgray](3,9)(4,9)(4,10)(3,9)
\psline*[linewidth=2pt,linecolor=darkgray](3,9)(3,10)(4,10)(3,9)
\psline*[linewidth=2pt,linecolor=lightgray](3,10)(4,10)(4,11)(3,10)
\psline*[linewidth=2pt,linecolor=darkgray](3,10)(3,11)(4,11)(3,10)
\psline*[linewidth=2pt,linecolor=lightgray](3,12)(3,13)(4,13)(3,12)
\psline*[linewidth=2pt,linecolor=darkgray](3,16)(4,16)(4,17)(3,16)
\psline*[linewidth=2pt,linecolor=lightgray](3,16)(3,17)(4,17)(3,16)
\psline*[linewidth=2pt,linecolor=darkgray](3,17)(4,17)(4,17.5)(3.5,17.5)(3,17)
\psline*[linewidth=2pt,linecolor=darkgray](4,0)(4,-0.5)(3.5,-0.5)(4,0)
\psline*[linewidth=2pt,linecolor=lightgray](3,17)(3,17.5)(3.5,17.5)(3,17)
\psline*[linewidth=2pt,linecolor=lightgray](3,0)(3,-0.5)(3.5,-0.5)(4,0)(3,0)
\psline*[linewidth=2pt,linecolor=darkgray](4,2)(5,2)(5,3)(4,2)
\psline*[linewidth=2pt,linecolor=darkgray](4,4)(5,4)(5,5)(4,4)
\psline*[linewidth=2pt,linecolor=lightgray](4,6)(5,6)(5,7)(4,6)
\psline*[linewidth=2pt,linecolor=darkgray](4,6)(4,7)(5,7)(4,6)
\psline*[linewidth=2pt,linecolor=lightgray](4,7)(5,7)(5,8)(4,7)
\psline*[linewidth=2pt,linecolor=darkgray](4,7)(4,8)(5,8)(4,7)
\psline*[linewidth=2pt,linecolor=lightgray](4,8)(5,8)(5,9)(4,8)
\psline*[linewidth=2pt,linecolor=darkgray](4,8)(4,9)(5,9)(4,8)
\psline*[linewidth=2pt,linecolor=lightgray](4,9)(5,9)(5,10)(4,9)
\psline*[linewidth=2pt,linecolor=darkgray](4,9)(4,10)(5,10)(4,9)
\psline*[linewidth=2pt,linecolor=lightgray](4,13)(5,13)(5,14)(4,13)
\psline*[linewidth=2pt,linecolor=darkgray](4,15)(5,15)(5,16)(4,15)
\psline*[linewidth=2pt,linecolor=lightgray](4,15)(4,16)(5,16)(4,15)
\psline*[linewidth=2pt,linecolor=darkgray](4,16)(5,16)(5,17)(4,16)
\psline*[linewidth=2pt,linecolor=lightgray](4,16)(4,17)(5,17)(4,16)
\psline*[linewidth=2pt,linecolor=lightgray](5,2)(5,3)(6,3)(5,2)
\psline*[linewidth=2pt,linecolor=darkgray](5,3)(6,3)(6,4)(5,3)
\psline*[linewidth=2pt,linecolor=lightgray](5,3)(5,4)(6,4)(5,3)
\psline*[linewidth=2pt,linecolor=darkgray](5,4)(6,4)(6,5)(5,4)
\psline*[linewidth=2pt,linecolor=lightgray](5,4)(5,5)(6,5)(5,4)
\psline*[linewidth=2pt,linecolor=darkgray](5,6)(5,7)(6,7)(5,6)
\psline*[linewidth=2pt,linecolor=lightgray](5,7)(6,7)(6,8)(5,7)
\psline*[linewidth=2pt,linecolor=darkgray](5,7)(5,8)(6,8)(5,7)
\psline*[linewidth=2pt,linecolor=lightgray](5,8)(6,8)(6,9)(5,8)
\psline*[linewidth=2pt,linecolor=darkgray](5,8)(5,9)(6,9)(5,8)
\psline*[linewidth=2pt,linecolor=lightgray](5,12)(6,12)(6,13)(5,12)
\psline*[linewidth=2pt,linecolor=darkgray](5,12)(5,13)(6,13)(5,12)
\psline*[linewidth=2pt,linecolor=lightgray](5,14)(5,15)(6,15)(5,14)
\psline*[linewidth=2pt,linecolor=darkgray](5,15)(6,15)(6,16)(5,15)
\psline*[linewidth=2pt,linecolor=lightgray](5,15)(5,16)(6,16)(5,15)
\psline*[linewidth=2pt,linecolor=lightgray](6,4)(6,5)(7,5)(6,4)
\psline*[linewidth=2pt,linecolor=darkgray](6,5)(7,5)(7,6)(6,5)
\psline*[linewidth=2pt,linecolor=lightgray](6,7)(7,7)(7,8)(6,7)
\psline*[linewidth=2pt,linecolor=darkgray](6,7)(6,8)(7,8)(6,7)
\psline*[linewidth=2pt,linecolor=lightgray](6,11)(7,11)(7,12)(6,11)
\psline*[linewidth=2pt,linecolor=darkgray](6,11)(6,12)(7,12)(6,11)
\psline*[linewidth=2pt,linecolor=darkgray](7,5)(8,5)(8,6)(7,5)
\psline*[linewidth=2pt,linecolor=lightgray](7,5)(7,6)(8,6)(7,5)
\psline*[linewidth=2pt,linecolor=darkgray](7,6)(8,6)(8,7)(7,6)
\psline*[linewidth=2pt,linecolor=lightgray](7,6)(7,7)(8,7)(7,6)
\psline*[linewidth=2pt,linecolor=lightgray](7,10)(8,10)(8,11)(7,10)
\psline*[linewidth=2pt,linecolor=darkgray](7,10)(7,11)(8,11)(7,10)
\psline*[linewidth=2pt,linecolor=darkgray](8,4)(9,4)(9,5)(8,4)
\psline*[linewidth=2pt,linecolor=lightgray](8,4)(8,5)(9,5)(8,4)
\psline*[linewidth=2pt,linecolor=darkgray](8,5)(9,5)(9,6)(8,5)
\psline*[linewidth=2pt,linecolor=lightgray](8,5)(8,6)(9,6)(8,5)
\psline*[linewidth=2pt,linecolor=lightgray](8,9)(9,9)(9,10)(8,9)
\psline*[linewidth=2pt,linecolor=darkgray](8,9)(8,10)(9,10)(8,9)
\psline*[linewidth=2pt,linecolor=darkgray](9,3)(10,3)(10,4)(9,3)
\psline*[linewidth=2pt,linecolor=lightgray](9,3)(9,4)(10,4)(9,3)
\psline*[linewidth=2pt,linecolor=darkgray](9,4)(10,4)(10,5)(9,4)
\psline*[linewidth=2pt,linecolor=lightgray](9,4)(9,5)(10,5)(9,4)
\psline*[linewidth=2pt,linecolor=darkgray](9,8)(10,8)(10,9)(9,8)
\psline*[linewidth=2pt,linecolor=lightgray](9,8)(9,9)(10,9)(9,8)
\psline*[linewidth=2pt,linecolor=darkgray](10,2)(11,2)(11,3)(10,2)
\psline*[linewidth=2pt,linecolor=lightgray](10,2)(10,3)(11,3)(10,2)
\psline*[linewidth=2pt,linecolor=darkgray](10,3)(11,3)(11,4)(10,3)
\psline*[linewidth=2pt,linecolor=lightgray](10,3)(10,4)(11,4)(10,3)
\psline*[linewidth=2pt,linecolor=darkgray](10,7)(11,7)(11,8)(10,7)
\psline*[linewidth=2pt,linecolor=lightgray](10,7)(10,8)(11,8)(10,7)
\psline*[linewidth=2pt,linecolor=lightgray](11,2)(11,3)(12,3)(11,2)
\psline*[linewidth=2pt,linecolor=darkgray](11,6)(12,6)(12,7)(11,6)
\psline*[linewidth=2pt,linecolor=lightgray](11,6)(11,7)(12,7)(11,6)
\psline*[linewidth=2pt,linecolor=lightgray](12,3)(13,3)(13,4)(12,3)
\psline*[linewidth=2pt,linecolor=darkgray](12,5)(13,5)(13,6)(12,5)
\psline*[linewidth=2pt,linecolor=lightgray](12,5)(12,6)(13,6)(12,5)
\psline*[linewidth=2pt,linecolor=lightgray](12,9)(13,9)(13,10)(12,9)
\psline*[linewidth=2pt,linecolor=darkgray](12,9)(12,10)(13,10)(12,9)
\psline*[linewidth=2pt,linecolor=lightgray](12,10)(13,10)(13,11)(12,10)
\psline*[linewidth=2pt,linecolor=lightgray](13,2)(14,2)(14,3)(13,2)
\psline*[linewidth=2pt,linecolor=darkgray](13,2)(13,3)(14,3)(13,2)
\psline*[linewidth=2pt,linecolor=lightgray](13,4)(13,5)(14,5)(13,4)
\psline*[linewidth=2pt,linecolor=lightgray](13,8)(14,8)(14,9)(13,8)
\psline*[linewidth=2pt,linecolor=darkgray](13,8)(13,9)(14,9)(13,8)
\psline*[linewidth=2pt,linecolor=lightgray](13,9)(14,9)(14,10)(13,9)
\psline*[linewidth=2pt,linecolor=darkgray](13,9)(13,10)(14,10)(13,9)
\psline*[linewidth=2pt,linecolor=lightgray](13,10)(14,10)(14,11)(13,10)
\psline*[linewidth=2pt,linecolor=darkgray](13,10)(13,11)(14,11)(13,10)
\psline*[linewidth=2pt,linecolor=lightgray](14,1)(15,1)(15,2)(14,1)
\psline*[linewidth=2pt,linecolor=darkgray](14,1)(14,2)(15,2)(14,1)
\psline*[linewidth=2pt,linecolor=lightgray](14,5)(15,5)(15,6)(14,5)
\psline*[linewidth=2pt,linecolor=lightgray](14,7)(15,7)(15,8)(14,7)
\psline*[linewidth=2pt,linecolor=darkgray](14,7)(14,8)(15,8)(14,7)
\psline*[linewidth=2pt,linecolor=lightgray](14,8)(15,8)(15,9)(14,8)
\psline*[linewidth=2pt,linecolor=darkgray](14,8)(14,9)(15,9)(14,8)
\psline*[linewidth=2pt,linecolor=lightgray](14,9)(15,9)(15,10)(14,9)
\psline*[linewidth=2pt,linecolor=darkgray](14,9)(14,10)(15,10)(14,9)
\psline*[linewidth=2pt,linecolor=lightgray](15,0)(16,0)(16,1)(15,0)
\psline*[linewidth=2pt,linecolor=darkgray](15,0)(15,1)(16,1)(15,0)
\psline*[linewidth=2pt,linecolor=lightgray](15,4)(16,4)(16,5)(15,4)
\psline*[linewidth=2pt,linecolor=darkgray](15,4)(15,5)(16,5)(15,4)
\psline*[linewidth=2pt,linecolor=lightgray](15,5)(16,5)(16,6)(15,5)
\psline*[linewidth=2pt,linecolor=darkgray](15,5)(15,6)(16,6)(15,5)
\psline*[linewidth=2pt,linecolor=lightgray](15,6)(16,6)(16,7)(15,6)
\psline*[linewidth=2pt,linecolor=darkgray](15,6)(15,7)(16,7)(15,6)
\psline*[linewidth=2pt,linecolor=lightgray](15,7)(16,7)(16,8)(15,7)
\psline*[linewidth=2pt,linecolor=darkgray](15,7)(15,8)(16,8)(15,7)
\psline*[linewidth=2pt,linecolor=lightgray](15,8)(16,8)(16,9)(15,8)
\psline*[linewidth=2pt,linecolor=darkgray](15,8)(15,9)(16,9)(15,8)
\psline*[linewidth=2pt,linecolor=lightgray](15,9)(16,9)(16,10)(15,9)
\psline*[linewidth=2pt,linecolor=darkgray](15,9)(15,10)(16,10)(15,9)
\psline*[linewidth=2pt,linecolor=lightgray](15,10)(16,10)(16,11)(15,10)
\psline*[linewidth=2pt,linecolor=lightgray](15,12)(16,12)(16,13)(15,12)
\psline*[linewidth=2pt,linecolor=darkgray](15,12)(15,13)(16,13)(15,12)
\psline*[linewidth=2pt,linecolor=lightgray](15,13)(16,13)(16,14)(15,13)
\psline*[linewidth=2pt,linecolor=lightgray](16,3)(17,3)(17,4)(16,3)
\psline*[linewidth=2pt,linecolor=darkgray](16,3)(16,4)(17,4)(16,3)
\psline*[linewidth=2pt,linecolor=lightgray](16,4)(17,4)(17,5)(16,4)
\psline*[linewidth=2pt,linecolor=darkgray](16,4)(16,5)(17,5)(16,4)
\psline*[linewidth=2pt,linecolor=lightgray](16,5)(17,5)(17,6)(16,5)
\psline*[linewidth=2pt,linecolor=darkgray](16,5)(16,6)(17,6)(16,5)
\psline*[linewidth=2pt,linecolor=lightgray](16,6)(17,6)(17,7)(16,6)
\psline*[linewidth=2pt,linecolor=darkgray](16,6)(16,7)(17,7)(16,6)
\psline*[linewidth=2pt,linecolor=lightgray](16,7)(17,7)(17,8)(16,7)
\psline*[linewidth=2pt,linecolor=darkgray](16,7)(16,8)(17,8)(16,7)
\psline*[linewidth=2pt,linecolor=lightgray](16,8)(17,8)(17,9)(16,8)
\psline*[linewidth=2pt,linecolor=darkgray](16,8)(16,9)(17,9)(16,8)
\psline*[linewidth=2pt,linecolor=darkgray](16,10)(16,11)(17,11)(16,10)
\psline*[linewidth=2pt,linecolor=lightgray](16,11)(17,11)(17,12)(16,11)
\psline*[linewidth=2pt,linecolor=darkgray](16,11)(16,12)(17,12)(16,11)
\psline*[linewidth=2pt,linecolor=lightgray](16,12)(17,12)(17,13)(16,12)
\psline*[linewidth=2pt,linecolor=darkgray](16,12)(16,13)(17,13)(16,12)
\psline*[linewidth=2pt,linecolor=lightgray](16,13)(17,13)(17,14)(16,13)
\psline*[linewidth=2pt,linecolor=darkgray](16,13)(16,14)(17,14)(16,13)
\psline*[linewidth=2pt,linecolor=lightgray](16,17)(17,17)(17,17.5)(16.5,17.5)(16,17)
\psline*[linewidth=2pt,linecolor=lightgray](17,0)(17,-0.5)(16.5,-0.5)(17,0)
\psline*[linewidth=2pt,linecolor=darkgray](16,17)(16,17.5)(16.5,17.5)(16,17)
\psline*[linewidth=2pt,linecolor=darkgray](16,0)(16,-0.5)(16.5,-0.5)(17,0)(16,0)
\psline*[linewidth=2pt,linecolor=lightgray](17,2)(17.5,2)(17.5,2.5)(17,2)
\psline*[linewidth=2pt,linecolor=lightgray](0,2)(-0.5,2)(-0.5,2.5)(0,3,)(0,2)
\psline*[linewidth=2pt,linecolor=darkgray](17,2)(17,3)(17.5,3)(17.5,2.5)(17,2)
\psline*[linewidth=2pt,linecolor=darkgray](0,3)(-0.5,3)(-0.5,2.5)(0,3)
\psline*[linewidth=2pt,linecolor=lightgray](17,3)(17.5,3)(17.5,3.5)(17,3)
\psline*[linewidth=2pt,linecolor=lightgray](0,3)(-0.5,3)(-0.5,3.5)(0,4,)(0,3)
\psline*[linewidth=2pt,linecolor=darkgray](17,3)(17,4)(17.5,4)(17.5,3.5)(17,3)
\psline*[linewidth=2pt,linecolor=darkgray](0,4)(-0.5,4)(-0.5,3.5)(0,4)
\psline*[linewidth=2pt,linecolor=lightgray](17,4)(17.5,4)(17.5,4.5)(17,4)
\psline*[linewidth=2pt,linecolor=lightgray](0,4)(-0.5,4)(-0.5,4.5)(0,5,)(0,4)
\psline*[linewidth=2pt,linecolor=darkgray](17,4)(17,5)(17.5,5)(17.5,4.5)(17,4)
\psline*[linewidth=2pt,linecolor=darkgray](0,5)(-0.5,5)(-0.5,4.5)(0,5)
\psline*[linewidth=2pt,linecolor=lightgray](17,5)(17.5,5)(17.5,5.5)(17,5)
\psline*[linewidth=2pt,linecolor=lightgray](0,5)(-0.5,5)(-0.5,5.5)(0,6,)(0,5)
\psline*[linewidth=2pt,linecolor=darkgray](17,5)(17,6)(17.5,6)(17.5,5.5)(17,5)
\psline*[linewidth=2pt,linecolor=darkgray](0,6)(-0.5,6)(-0.5,5.5)(0,6)
\psline*[linewidth=2pt,linecolor=lightgray](17,6)(17.5,6)(17.5,6.5)(17,6)
\psline*[linewidth=2pt,linecolor=lightgray](0,6)(-0.5,6)(-0.5,6.5)(0,7,)(0,6)
\psline*[linewidth=2pt,linecolor=darkgray](17,6)(17,7)(17.5,7)(17.5,6.5)(17,6)
\psline*[linewidth=2pt,linecolor=darkgray](0,7)(-0.5,7)(-0.5,6.5)(0,7)
\psline*[linewidth=2pt,linecolor=lightgray](17,7)(17.5,7)(17.5,7.5)(17,7)
\psline*[linewidth=2pt,linecolor=lightgray](0,7)(-0.5,7)(-0.5,7.5)(0,8,)(0,7)
\psline*[linewidth=2pt,linecolor=darkgray](17,7)(17,8)(17.5,8)(17.5,7.5)(17,7)
\psline*[linewidth=2pt,linecolor=darkgray](0,8)(-0.5,8)(-0.5,7.5)(0,8)
\psline*[linewidth=2pt,linecolor=lightgray](17,11)(17.5,11)(17.5,11.5)(17,11)
\psline*[linewidth=2pt,linecolor=lightgray](0,11)(-0.5,11)(-0.5,11.5)(0,12,)(0,11)
\psline*[linewidth=2pt,linecolor=darkgray](17,11)(17,12)(17.5,12)(17.5,11.5)(17,11)
\psline*[linewidth=2pt,linecolor=darkgray](0,12)(-0.5,12)(-0.5,11.5)(0,12)
\psline*[linewidth=2pt,linecolor=lightgray](17,12)(17.5,12)(17.5,12.5)(17,12)
\psline*[linewidth=2pt,linecolor=lightgray](0,12)(-0.5,12)(-0.5,12.5)(0,13,)(0,12)
\psline*[linewidth=2pt,linecolor=darkgray](17,12)(17,13)(17.5,13)(17.5,12.5)(17,12)
\psline*[linewidth=2pt,linecolor=darkgray](0,13)(-0.5,13)(-0.5,12.5)(0,13)
\psline*[linewidth=2pt,linecolor=darkgray](17,16)(17.5,16)(17.5,16.5)(17,16)
\psline*[linewidth=2pt,linecolor=darkgray](0,16)(-0.5,16)(-0.5,16.5)(0,17,)(0,16)
\psline*[linewidth=2pt,linecolor=lightgray](17,16)(17,17)(17.5,17)(17.5,16.5)(17,16)
\psline*[linewidth=2pt,linecolor=lightgray](0,17)(-0.5,17)(-0.5,16.5)(0,17)
\psline[linewidth=2pt,linecolor=blue](-0.5,0)(17.5,0)
\psline[linewidth=2pt,linecolor=blue](-0.5,1)(17.5,1)
\psline[linewidth=2pt,linecolor=blue](-0.5,2)(17.5,2)
\psline[linewidth=2pt,linecolor=blue](-0.5,3)(17.5,3)
\psline[linewidth=2pt,linecolor=blue](-0.5,4)(17.5,4)
\psline[linewidth=2pt,linecolor=blue](-0.5,5)(17.5,5)
\psline[linewidth=2pt,linecolor=blue](-0.5,6)(17.5,6)
\psline[linewidth=2pt,linecolor=blue](-0.5,7)(17.5,7)
\psline[linewidth=2pt,linecolor=blue](-0.5,8)(17.5,8)
\psline[linewidth=2pt,linecolor=blue](-0.5,9)(17.5,9)
\psline[linewidth=2pt,linecolor=blue](-0.5,10)(17.5,10)
\psline[linewidth=2pt,linecolor=blue](-0.5,11)(17.5,11)
\psline[linewidth=2pt,linecolor=blue](-0.5,12)(17.5,12)
\psline[linewidth=2pt,linecolor=blue](-0.5,13)(17.5,13)
\psline[linewidth=2pt,linecolor=blue](-0.5,14)(17.5,14)
\psline[linewidth=2pt,linecolor=blue](-0.5,15)(17.5,15)
\psline[linewidth=2pt,linecolor=blue](-0.5,16)(17.5,16)
\psline[linewidth=2pt,linecolor=blue](-0.5,17)(17.5,17)
\psline[linewidth=2pt,linecolor=blue](0,-0.5)(0,17.5)
\psline[linewidth=2pt,linecolor=blue](1,-0.5)(1,17.5)
\psline[linewidth=2pt,linecolor=blue](2,-0.5)(2,17.5)
\psline[linewidth=2pt,linecolor=blue](3,-0.5)(3,17.5)
\psline[linewidth=2pt,linecolor=blue](4,-0.5)(4,17.5)
\psline[linewidth=2pt,linecolor=blue](5,-0.5)(5,17.5)
\psline[linewidth=2pt,linecolor=blue](6,-0.5)(6,17.5)
\psline[linewidth=2pt,linecolor=blue](7,-0.5)(7,17.5)
\psline[linewidth=2pt,linecolor=blue](8,-0.5)(8,17.5)
\psline[linewidth=2pt,linecolor=blue](9,-0.5)(9,17.5)
\psline[linewidth=2pt,linecolor=blue](10,-0.5)(10,17.5)
\psline[linewidth=2pt,linecolor=blue](11,-0.5)(11,17.5)
\psline[linewidth=2pt,linecolor=blue](12,-0.5)(12,17.5)
\psline[linewidth=2pt,linecolor=blue](13,-0.5)(13,17.5)
\psline[linewidth=2pt,linecolor=blue](14,-0.5)(14,17.5)
\psline[linewidth=2pt,linecolor=blue](15,-0.5)(15,17.5)
\psline[linewidth=2pt,linecolor=blue](16,-0.5)(16,17.5)
\psline[linewidth=2pt,linecolor=blue](17,-0.5)(17,17.5)
\psline[linewidth=2pt,linecolor=blue](-0.5,-0.5)(17.5,17.5)
\psline[linewidth=2pt,linecolor=blue](0.5,-0.5)(17.5,16.5)
\psline[linewidth=2pt,linecolor=blue](1.5,-0.5)(17.5,15.5)
\psline[linewidth=2pt,linecolor=blue](2.5,-0.5)(17.5,14.5)
\psline[linewidth=2pt,linecolor=blue](3.5,-0.5)(17.5,13.5)
\psline[linewidth=2pt,linecolor=blue](4.5,-0.5)(17.5,12.5)
\psline[linewidth=2pt,linecolor=blue](5.5,-0.5)(17.5,11.5)
\psline[linewidth=2pt,linecolor=blue](6.5,-0.5)(17.5,10.5)
\psline[linewidth=2pt,linecolor=blue](7.5,-0.5)(17.5,9.5)
\psline[linewidth=2pt,linecolor=blue](8.5,-0.5)(17.5,8.5)
\psline[linewidth=2pt,linecolor=blue](9.5,-0.5)(17.5,7.5)
\psline[linewidth=2pt,linecolor=blue](10.5,-0.5)(17.5,6.5)
\psline[linewidth=2pt,linecolor=blue](11.5,-0.5)(17.5,5.5)
\psline[linewidth=2pt,linecolor=blue](12.5,-0.5)(17.5,4.5)
\psline[linewidth=2pt,linecolor=blue](13.5,-0.5)(17.5,3.5)
\psline[linewidth=2pt,linecolor=blue](14.5,-0.5)(17.5,2.5)
\psline[linewidth=2pt,linecolor=blue](15.5,-0.5)(17.5,1.5)
\psline[linewidth=2pt,linecolor=blue](16.5,-0.5)(17.5,0.5)
\psline[linewidth=2pt,linecolor=blue](-0.5,0.5)(16.5,17.5)
\psline[linewidth=2pt,linecolor=blue](-0.5,1.5)(15.5,17.5)
\psline[linewidth=2pt,linecolor=blue](-0.5,2.5)(14.5,17.5)
\psline[linewidth=2pt,linecolor=blue](-0.5,3.5)(13.5,17.5)
\psline[linewidth=2pt,linecolor=blue](-0.5,4.5)(12.5,17.5)
\psline[linewidth=2pt,linecolor=blue](-0.5,5.5)(11.5,17.5)
\psline[linewidth=2pt,linecolor=blue](-0.5,6.5)(10.5,17.5)
\psline[linewidth=2pt,linecolor=blue](-0.5,7.5)(9.5,17.5)
\psline[linewidth=2pt,linecolor=blue](-0.5,8.5)(8.5,17.5)
\psline[linewidth=2pt,linecolor=blue](-0.5,9.5)(7.5,17.5)
\psline[linewidth=2pt,linecolor=blue](-0.5,10.5)(6.5,17.5)
\psline[linewidth=2pt,linecolor=blue](-0.5,11.5)(5.5,17.5)
\psline[linewidth=2pt,linecolor=blue](-0.5,12.5)(4.5,17.5)
\psline[linewidth=2pt,linecolor=blue](-0.5,13.5)(3.5,17.5)
\psline[linewidth=2pt,linecolor=blue](-0.5,14.5)(2.5,17.5)
\psline[linewidth=2pt,linecolor=blue](-0.5,15.5)(1.5,17.5)
\psline[linewidth=2pt,linecolor=blue](-0.5,16.5)(0.5,17.5)
\multirput{0}(0,0)(0,1){18}{%
  \multirput{0}(0,0)(1,0){18}{%
     \pscircle*[linecolor=white]{8pt}
     \pscircle[linewidth=1pt,linecolor=black] {8pt}
   }
}
\rput{0}(0,0){{\bf 4}}
\rput{0}(0,1){{\bf 1}}
\rput{0}(0,2){{\bf 3}}
\rput{0}(0,3){{\bf 2}}
\rput{0}(0,4){{\bf 1}}
\rput{0}(0,5){{\bf 3}}
\rput{0}(0,6){{\bf 2}}
\rput{0}(0,7){{\bf 1}}
\rput{0}(0,8){{\bf 3}}
\rput{0}(0,9){{\bf 4}}
\rput{0}(0,10){{\bf 1}}
\rput{0}(0,11){{\bf 3}}
\rput{0}(0,12){{\bf 2}}
\rput{0}(0,13){{\bf 1}}
\rput{0}(0,14){{\bf 4}}
\rput{0}(0,15){{\bf 2}}
\rput{0}(0,16){{\bf 3}}
\rput{0}(0,17){{\bf 1}}
\rput{0}(1,0){{\bf 2}}
\rput{0}(1,1){{\bf 3}}
\rput{0}(1,2){{\bf 2}}
\rput{0}(1,3){{\bf 1}}
\rput{0}(1,4){{\bf 3}}
\rput{0}(1,5){{\bf 2}}
\rput{0}(1,6){{\bf 1}}
\rput{0}(1,7){{\bf 3}}
\rput{0}(1,8){{\bf 4}}
\rput{0}(1,9){{\bf 1}}
\rput{0}(1,10){{\bf 3}}
\rput{0}(1,11){{\bf 2}}
\rput{0}(1,12){{\bf 1}}
\rput{0}(1,13){{\bf 4}}
\rput{0}(1,14){{\bf 2}}
\rput{0}(1,15){{\bf 3}}
\rput{0}(1,16){{\bf 1}}
\rput{0}(1,17){{\bf 4}}
\rput{0}(2,0){{\bf 3}}
\rput{0}(2,1){{\bf 1}}
\rput{0}(2,2){{\bf 4}}
\rput{0}(2,3){{\bf 3}}
\rput{0}(2,4){{\bf 2}}
\rput{0}(2,5){{\bf 1}}
\rput{0}(2,6){{\bf 3}}
\rput{0}(2,7){{\bf 4}}
\rput{0}(2,8){{\bf 1}}
\rput{0}(2,9){{\bf 3}}
\rput{0}(2,10){{\bf 2}}
\rput{0}(2,11){{\bf 1}}
\rput{0}(2,12){{\bf 4}}
\rput{0}(2,13){{\bf 2}}
\rput{0}(2,14){{\bf 3}}
\rput{0}(2,15){{\bf 1}}
\rput{0}(2,16){{\bf 4}}
\rput{0}(2,17){{\bf 2}}
\rput{0}(3,0){{\bf 1}}
\rput{0}(3,1){{\bf 2}}
\rput{0}(3,2){{\bf 3}}
\rput{0}(3,3){{\bf 2}}
\rput{0}(3,4){{\bf 1}}
\rput{0}(3,5){{\bf 3}}
\rput{0}(3,6){{\bf 4}}
\rput{0}(3,7){{\bf 1}}
\rput{0}(3,8){{\bf 3}}
\rput{0}(3,9){{\bf 2}}
\rput{0}(3,10){{\bf 1}}
\rput{0}(3,11){{\bf 3}}
\rput{0}(3,12){{\bf 2}}
\rput{0}(3,13){{\bf 3}}
\rput{0}(3,14){{\bf 1}}
\rput{0}(3,15){{\bf 4}}
\rput{0}(3,16){{\bf 2}}
\rput{0}(3,17){{\bf 3}}
\rput{0}(4,0){{\bf 2}}
\rput{0}(4,1){{\bf 4}}
\rput{0}(4,2){{\bf 1}}
\rput{0}(4,3){{\bf 4}}
\rput{0}(4,4){{\bf 3}}
\rput{0}(4,5){{\bf 4}}
\rput{0}(4,6){{\bf 1}}
\rput{0}(4,7){{\bf 3}}
\rput{0}(4,8){{\bf 2}}
\rput{0}(4,9){{\bf 1}}
\rput{0}(4,10){{\bf 3}}
\rput{0}(4,11){{\bf 2}}
\rput{0}(4,12){{\bf 4}}
\rput{0}(4,13){{\bf 1}}
\rput{0}(4,14){{\bf 4}}
\rput{0}(4,15){{\bf 2}}
\rput{0}(4,16){{\bf 3}}
\rput{0}(4,17){{\bf 1}}
\rput{0}(5,0){{\bf 4}}
\rput{0}(5,1){{\bf 1}}
\rput{0}(5,2){{\bf 2}}
\rput{0}(5,3){{\bf 3}}
\rput{0}(5,4){{\bf 1}}
\rput{0}(5,5){{\bf 2}}
\rput{0}(5,6){{\bf 3}}
\rput{0}(5,7){{\bf 2}}
\rput{0}(5,8){{\bf 1}}
\rput{0}(5,9){{\bf 3}}
\rput{0}(5,10){{\bf 2}}
\rput{0}(5,11){{\bf 4}}
\rput{0}(5,12){{\bf 1}}
\rput{0}(5,13){{\bf 3}}
\rput{0}(5,14){{\bf 2}}
\rput{0}(5,15){{\bf 3}}
\rput{0}(5,16){{\bf 1}}
\rput{0}(5,17){{\bf 2}}
\rput{0}(6,0){{\bf 1}}
\rput{0}(6,1){{\bf 2}}
\rput{0}(6,2){{\bf 4}}
\rput{0}(6,3){{\bf 1}}
\rput{0}(6,4){{\bf 2}}
\rput{0}(6,5){{\bf 3}}
\rput{0}(6,6){{\bf 4}}
\rput{0}(6,7){{\bf 1}}
\rput{0}(6,8){{\bf 3}}
\rput{0}(6,9){{\bf 2}}
\rput{0}(6,10){{\bf 4}}
\rput{0}(6,11){{\bf 1}}
\rput{0}(6,12){{\bf 3}}
\rput{0}(6,13){{\bf 2}}
\rput{0}(6,14){{\bf 4}}
\rput{0}(6,15){{\bf 1}}
\rput{0}(6,16){{\bf 2}}
\rput{0}(6,17){{\bf 4}}
\rput{0}(7,0){{\bf 2}}
\rput{0}(7,1){{\bf 4}}
\rput{0}(7,2){{\bf 3}}
\rput{0}(7,3){{\bf 2}}
\rput{0}(7,4){{\bf 4}}
\rput{0}(7,5){{\bf 1}}
\rput{0}(7,6){{\bf 2}}
\rput{0}(7,7){{\bf 3}}
\rput{0}(7,8){{\bf 2}}
\rput{0}(7,9){{\bf 4}}
\rput{0}(7,10){{\bf 1}}
\rput{0}(7,11){{\bf 3}}
\rput{0}(7,12){{\bf 2}}
\rput{0}(7,13){{\bf 4}}
\rput{0}(7,14){{\bf 1}}
\rput{0}(7,15){{\bf 2}}
\rput{0}(7,16){{\bf 4}}
\rput{0}(7,17){{\bf 1}}
\rput{0}(8,0){{\bf 4}}
\rput{0}(8,1){{\bf 3}}
\rput{0}(8,2){{\bf 2}}
\rput{0}(8,3){{\bf 4}}
\rput{0}(8,4){{\bf 1}}
\rput{0}(8,5){{\bf 2}}
\rput{0}(8,6){{\bf 3}}
\rput{0}(8,7){{\bf 1}}
\rput{0}(8,8){{\bf 4}}
\rput{0}(8,9){{\bf 1}}
\rput{0}(8,10){{\bf 3}}
\rput{0}(8,11){{\bf 2}}
\rput{0}(8,12){{\bf 4}}
\rput{0}(8,13){{\bf 1}}
\rput{0}(8,14){{\bf 2}}
\rput{0}(8,15){{\bf 4}}
\rput{0}(8,16){{\bf 1}}
\rput{0}(8,17){{\bf 2}}
\rput{0}(9,0){{\bf 3}}
\rput{0}(9,1){{\bf 2}}
\rput{0}(9,2){{\bf 4}}
\rput{0}(9,3){{\bf 1}}
\rput{0}(9,4){{\bf 2}}
\rput{0}(9,5){{\bf 3}}
\rput{0}(9,6){{\bf 1}}
\rput{0}(9,7){{\bf 4}}
\rput{0}(9,8){{\bf 2}}
\rput{0}(9,9){{\bf 3}}
\rput{0}(9,10){{\bf 2}}
\rput{0}(9,11){{\bf 4}}
\rput{0}(9,12){{\bf 1}}
\rput{0}(9,13){{\bf 2}}
\rput{0}(9,14){{\bf 4}}
\rput{0}(9,15){{\bf 1}}
\rput{0}(9,16){{\bf 2}}
\rput{0}(9,17){{\bf 4}}
\rput{0}(10,0){{\bf 2}}
\rput{0}(10,1){{\bf 4}}
\rput{0}(10,2){{\bf 1}}
\rput{0}(10,3){{\bf 2}}
\rput{0}(10,4){{\bf 3}}
\rput{0}(10,5){{\bf 1}}
\rput{0}(10,6){{\bf 4}}
\rput{0}(10,7){{\bf 2}}
\rput{0}(10,8){{\bf 3}}
\rput{0}(10,9){{\bf 1}}
\rput{0}(10,10){{\bf 4}}
\rput{0}(10,11){{\bf 1}}
\rput{0}(10,12){{\bf 2}}
\rput{0}(10,13){{\bf 4}}
\rput{0}(10,14){{\bf 1}}
\rput{0}(10,15){{\bf 2}}
\rput{0}(10,16){{\bf 4}}
\rput{0}(10,17){{\bf 3}}
\rput{0}(11,0){{\bf 4}}
\rput{0}(11,1){{\bf 1}}
\rput{0}(11,2){{\bf 2}}
\rput{0}(11,3){{\bf 3}}
\rput{0}(11,4){{\bf 1}}
\rput{0}(11,5){{\bf 4}}
\rput{0}(11,6){{\bf 2}}
\rput{0}(11,7){{\bf 3}}
\rput{0}(11,8){{\bf 1}}
\rput{0}(11,9){{\bf 4}}
\rput{0}(11,10){{\bf 2}}
\rput{0}(11,11){{\bf 3}}
\rput{0}(11,12){{\bf 4}}
\rput{0}(11,13){{\bf 1}}
\rput{0}(11,14){{\bf 2}}
\rput{0}(11,15){{\bf 4}}
\rput{0}(11,16){{\bf 3}}
\rput{0}(11,17){{\bf 2}}
\rput{0}(12,0){{\bf 1}}
\rput{0}(12,1){{\bf 2}}
\rput{0}(12,2){{\bf 4}}
\rput{0}(12,3){{\bf 1}}
\rput{0}(12,4){{\bf 4}}
\rput{0}(12,5){{\bf 2}}
\rput{0}(12,6){{\bf 3}}
\rput{0}(12,7){{\bf 1}}
\rput{0}(12,8){{\bf 4}}
\rput{0}(12,9){{\bf 2}}
\rput{0}(12,10){{\bf 1}}
\rput{0}(12,11){{\bf 4}}
\rput{0}(12,12){{\bf 1}}
\rput{0}(12,13){{\bf 2}}
\rput{0}(12,14){{\bf 4}}
\rput{0}(12,15){{\bf 3}}
\rput{0}(12,16){{\bf 2}}
\rput{0}(12,17){{\bf 4}}
\rput{0}(13,0){{\bf 2}}
\rput{0}(13,1){{\bf 4}}
\rput{0}(13,2){{\bf 1}}
\rput{0}(13,3){{\bf 3}}
\rput{0}(13,4){{\bf 2}}
\rput{0}(13,5){{\bf 3}}
\rput{0}(13,6){{\bf 1}}
\rput{0}(13,7){{\bf 4}}
\rput{0}(13,8){{\bf 2}}
\rput{0}(13,9){{\bf 1}}
\rput{0}(13,10){{\bf 3}}
\rput{0}(13,11){{\bf 2}}
\rput{0}(13,12){{\bf 3}}
\rput{0}(13,13){{\bf 4}}
\rput{0}(13,14){{\bf 3}}
\rput{0}(13,15){{\bf 2}}
\rput{0}(13,16){{\bf 4}}
\rput{0}(13,17){{\bf 1}}
\rput{0}(14,0){{\bf 4}}
\rput{0}(14,1){{\bf 1}}
\rput{0}(14,2){{\bf 3}}
\rput{0}(14,3){{\bf 2}}
\rput{0}(14,4){{\bf 4}}
\rput{0}(14,5){{\bf 1}}
\rput{0}(14,6){{\bf 4}}
\rput{0}(14,7){{\bf 2}}
\rput{0}(14,8){{\bf 1}}
\rput{0}(14,9){{\bf 3}}
\rput{0}(14,10){{\bf 2}}
\rput{0}(14,11){{\bf 1}}
\rput{0}(14,12){{\bf 4}}
\rput{0}(14,13){{\bf 2}}
\rput{0}(14,14){{\bf 1}}
\rput{0}(14,15){{\bf 4}}
\rput{0}(14,16){{\bf 1}}
\rput{0}(14,17){{\bf 2}}
\rput{0}(15,0){{\bf 1}}
\rput{0}(15,1){{\bf 3}}
\rput{0}(15,2){{\bf 2}}
\rput{0}(15,3){{\bf 4}}
\rput{0}(15,4){{\bf 1}}
\rput{0}(15,5){{\bf 3}}
\rput{0}(15,6){{\bf 2}}
\rput{0}(15,7){{\bf 1}}
\rput{0}(15,8){{\bf 3}}
\rput{0}(15,9){{\bf 2}}
\rput{0}(15,10){{\bf 1}}
\rput{0}(15,11){{\bf 4}}
\rput{0}(15,12){{\bf 2}}
\rput{0}(15,13){{\bf 1}}
\rput{0}(15,14){{\bf 4}}
\rput{0}(15,15){{\bf 3}}
\rput{0}(15,16){{\bf 2}}
\rput{0}(15,17){{\bf 4}}
\rput{0}(16,0){{\bf 3}}
\rput{0}(16,1){{\bf 2}}
\rput{0}(16,2){{\bf 4}}
\rput{0}(16,3){{\bf 1}}
\rput{0}(16,4){{\bf 3}}
\rput{0}(16,5){{\bf 2}}
\rput{0}(16,6){{\bf 1}}
\rput{0}(16,7){{\bf 3}}
\rput{0}(16,8){{\bf 2}}
\rput{0}(16,9){{\bf 1}}
\rput{0}(16,10){{\bf 3}}
\rput{0}(16,11){{\bf 2}}
\rput{0}(16,12){{\bf 1}}
\rput{0}(16,13){{\bf 3}}
\rput{0}(16,14){{\bf 2}}
\rput{0}(16,15){{\bf 1}}
\rput{0}(16,16){{\bf 4}}
\rput{0}(16,17){{\bf 1}}
\rput{0}(17,0){{\bf 2}}
\rput{0}(17,1){{\bf 4}}
\rput{0}(17,2){{\bf 1}}
\rput{0}(17,3){{\bf 3}}
\rput{0}(17,4){{\bf 2}}
\rput{0}(17,5){{\bf 1}}
\rput{0}(17,6){{\bf 3}}
\rput{0}(17,7){{\bf 2}}
\rput{0}(17,8){{\bf 1}}
\rput{0}(17,9){{\bf 3}}
\rput{0}(17,10){{\bf 4}}
\rput{0}(17,11){{\bf 1}}
\rput{0}(17,12){{\bf 3}}
\rput{0}(17,13){{\bf 2}}
\rput{0}(17,14){{\bf 1}}
\rput{0}(17,15){{\bf 4}}
\rput{0}(17,16){{\bf 2}}
\rput{0}(17,17){{\bf 3}}
\multirput{0}(17.5,-0.5)(0,1){18}{\psline[linewidth=2pt,linecolor=black]{->}(0.2,-0.2)(-0.1,0.1)}
\uput[0](17.7,-0.7){D1}
\uput[0](17.7,0.3){D2}
\uput[0](17.7,1.3){D3}
\uput[0](17.7,2.3){D4}
\uput[0](17.7,3.3){D5}
\uput[0](17.7,4.3){D6}
\uput[0](17.7,5.3){D7}
\uput[0](17.7,6.3){D8}
\uput[0](17.7,7.3){D9}
\uput[0](17.7,8.3){D10}
\uput[0](17.7,9.3){D11}
\uput[0](17.7,10.3){D12}
\uput[0](17.7,11.3){D13}
\uput[0](17.7,12.3){D14}
\uput[0](17.7,13.3){D15}
\uput[0](17.7,14.3){D16}
\uput[0](17.7,15.3){D17}
\uput[0](17.7,16.3){D18}
\endpspicture
\caption{ \label{prop.12k-6.fig3}
The 4-coloring of $T(18,18)$ after Step~4 in the case $L=4k-2$.
}
\end{figure}
%
%

\medskip
\noindent
{\bf Step 4.}
There are only five counter-diagonals to be colored. 
All vertices on D$(6k-4)$ are colored $1$ or $2$ using the following simple
rule: the vertex $(x,y)$ is colored $1$ (resp.\ $2$) if the vertex 
$(x,y-1)$ is colored $4$ (resp.\  $3$). In particular, those vertices with 
$3k-1 \leq x \leq 9k-5$ are colored alike. 

On D$(6k)$ we find two vertices admitting a single color in the set $\{1,2\}$:
$(3k-1,3k-1)$ and ($9k-2,9k-4)$ taking respectively, colors $c_1$ and $c_2$.
The vertices satisfying $3k\leq x \leq 9k-4$ are colored $c_1$, and the others
are colored $c_2$.

On D$(6k-3)$ we find two vertices $(3k-1,3k-2)$ and $(9k-4,9k-5)$ 
that admit a single color from the set $\{3,4\}$. The other vertices are 
colored $1$ and $2$ (there is a unique choice for each vertex).

On D$(6k-2)$ there are four vertices admitting a single color from the set
$\{3,4\}$: the vertices $(3k,3k-2)$ and $(3k-1,3k-1)$ are colored $c_1$,
while $(9k-3,9k-5)$ and $(9k-4,9k-4)$ are colored $c_2\neq c_1$. Those
vertices satisfying $3k+1\leq x \leq 9k-2$ are colored $c_2$, and the rest
are colored $c_1$.  

The last counter-diagonal D$(6k-1)$ contains seven vertices that admit a 
single color: $(3k+1,3k-2)$, $(3k,3k-1)$, $(3k-1,3k)$, $(9k-1,9k-6)$,
$(9k-2,9k-5)$, $(9k-3,9k-4)$, and $(9k-4,9k-3)$. The other vertices are
colored $3$ and $4$ (there is a unique choice for each vertex). 

The resulting coloring is depicted in
Figure~\ref{prop.12k-6.fig3}. The contribution to the partial degree of the
new triangles is $2$; the partial degree of $f$ is given by
$\deg f|_R = 6 + 12(k-2) \equiv 6 \pmod{12}$.

\medskip

The above argument proves the base case of our induction. Now we have to
find a four-coloring of the triangulation $T(12k-6,3)$ with $k\geq 2$
such that it has the same top-row coloring $c_{3}$ as $f$ 
(see Figure~\ref{prop.12k-6.fig3}). We proceed as for the previous cases: 
let $t = \lfloor \tfrac{3k-6}{2} \rfloor$; the 4-coloring we need
is defined as follows for $k$ even:
\begin{eqnarray*}
 c_0 = c_3 & = &  [1423]^{t+1} 1241 2432 4124 1 [3241]^t 3 \\
       c_2 & = & 3[1423]^{t+1} 1241 243  2413   [2413]^t 42 \;=\; 
                 3[1423]^{t+1} 1241 243         [2413]^{t+1} 42  \\ 
       c_1 & = &  [2314]^{t+1} 2312 4124 3 2413 [2413]^t 4 \;=\; 
                  [2314]^{t+1} 2312 4124 3 [2413]^{t+1} 4\,.
\end{eqnarray*}
If $k$ is odd, then we have:
\begin{eqnarray*}
 c_0 =c_3 & = & [1423]^{t+1} 1421 3213 4132 13 [2413]^{t+1} \\
      c_2 & = & [3142]^{t+1} 3142 1321 3413    [2413]^{t+1} 42 \;=\; 
                [3142]^{t+2}      1321 3413    [2413]^{t+1} 42 \\
      c_1 & = & [2314]^{t+1} 2314 2132 1341 3  [2413]^{t+1} 4 \;=\; 
                [2314]^{t+2}      2132 1341 3  [2413]^{t+1} 4 \,. 
\end{eqnarray*}
Again, it is easy to verify that this gives a proper 4-coloring of $T(3L,3)$,
and by Proposition~\ref{prop.T_Lx3}, it has zero degree.
This completes the proof of the theorem. \qed

\medskip

Theorems~\ref{theo.main} and~\ref{theo.asym} imply that WSK is non
ergodic on any triangulation $T(3L,3M)$ with $3\leq L\leq M$. 
Proposition~\ref{prop.T_Lx3} together with Fisk's theorem implies that
WSK is ergodic on any triangulation $T(3,3L)$. The triangulations 
$T(6,3L)$ are special in the sense that WSK is ergodic depending on 
the value of $L$. In particular, WSK is not ergodic for any $T(6,6p)$ with 
odd $p$, because of Theorem~\ref{theo.main} [or Theorem~\ref{theo_L=6}]
and Lemma~\ref{lemma.tech}. 

By direct computer enumeration of the $299 146 792$ proper four-colorings of
$T(6,9)$, we have checked that all of them have zero degree. 
We have also checked with a computer that we can transform any of these
colorings into the three-coloring by a {\em finite} number of K--changes.
Therefore we have obtained a computer--assisted proof of the following
Theorem:

\begin{proposition} \label{prop_tri_L=6x9}
$\Kc(T(6,9),4) = 1$
\end{proposition}
 
\medskip

\noindent
{\bf Remark.} Fisk's Theorem~\ref{theo_Fisk} can be used to prove the
ergodicity of the WSK on $T(6,9)$ directly from the fact that all
colorings have zero degree. 

%
%
\section{Summary and open problems} \label{sec.summary}

We have considered the question of the ergodicity of the 
Wang--Swendsen--Koteck\'y dynamics for the zero-temperature 
4--state Potts antiferromagnet on triangulations $T(3L,3M)$ of the torus.
This dynamics is equivalent (for the zero-temperature case only) to that
of the Kempe chains studied in Combinatorics. We have obtained two
main results:

1) For the wider family of the even triangulations of the torus (which 
contains the triangulations $T(3L,3M)$ as a proper subset), we find that the
degree of a 4--coloring modulo 12 is invariant under Kempe changes. 

2) For any triangulation $T(3L,3M)$ of the torus with $3\le L \le M$, 
there are at least two Kempe equivalence classes for 4 colors. In other
words, the Wang--Swendsen--Koteck\'y dynamics with 4 colors on these
triangulations is non-ergodic. For $L=2$, we can only show that this 
dynamics is non-ergodic for $M=2p$ with odd $p$. 

In addition to their intrinsic mathematical interest, these results have 
a great practical importance in Statistical Mechanics. The 
triangular-lattice 4--state Potts antiferromagnet is believed to have
a zero temperature critical point \cite[and references therein]{transfer3}.
But we {\em cannot}\/ study the critical properties of this model
using WSK dynamics because of the non-ergodicity of the algorithm. (This also
holds for the single-site Metropolis dynamics, as it corresponds 
to a particular subset of moves of the WSK dynamics.) 
Indeed, one can simulate the 4--state Potts antiferromagnet at zero temperature
using the WSK algorithm on planar graphs (e.g., a triangular grid with
free boundary conditions); but surface effects cannot be eliminated, and one
has to go to much larger lattice sizes to attain high--precision results.
It is therefore important to devise a new Monte Carlo algorithm for
this model which is ergodic at zero temperature.

There are other open problems related to the ergodicity of the Kempe 
dynamics. The case of four-colors on triangulations of the torus is rather
special, as we can make use of concepts borrowed from Algebraic Topology.
However these techniques cannot be applied to the cases of $q=5,6$ colors,
and the ergodicity of the corresponding WSK dynamics is still an open problem.

Finally, let us mention that {\em at zero temperature}, the 4--state Potts
model on the triangular lattice is essentially equivalent to the 
3--state Potts model on the kagom\'e lattice. We have found that the
WSK dynamics for this model also fails to be ergodic on most 
kagom\'e lattices when embedded on a torus. The details will be published 
elsewhere.

%
%
\section*{Acknowledgments}

We are indebted to Alan Sokal for his participation on the early stages
of this work, and his encouragement and useful suggestions later on.
We also wish to thank Eduardo J.S. Villase\~nor for useful discussions. 

J.S.\ is grateful for the kind hospitality of
the Physics Department of New York University and the Mathematics 
Department of University College London, where part of this work was done; and  
also thanks the Isaac Newton Institute for Mathematical Sciences,
University of Cambridge, for hospitality during the programme on
Combinatorics and Statistical Mechanics (January--June 2008).

The authors' research was supported in part by 
the ARRS (Slovenia) Research Program P1--0297, by an NSERC Discovery Grant,
 and by the Canada Research Chair program (B.M.), 
by U.S.\ National Science Foundation grants
PHY--0116590 and PHY--0424082,  
and by Spanish MEC grants MTM2005--08618 and MTM2008--03020 (J.S.).

%
%

\end{document}